\documentclass[11pt,           
english,     
]{article}

\usepackage[english,german,strings]{babel}     
\usepackage{amsmath}           
\usepackage[utf8]{inputenc}   
\usepackage[T1]{fontenc}      
\usepackage[sort]{cite}        
\usepackage{fancyhdr}          
\usepackage[a4paper]{geometry} 
\usepackage{xspace}            
\usepackage{tikz}              
\usetikzlibrary{arrows,positioning}

\usepackage{nchairx} 
\usepackage[expansion=false    
]{microtype}        
\usepackage[%
final=true,         
pdfpagelabels,		
hypertexnames=false 
]{hyperref}

\makeatletter
\newcommand{\chairxauthorbibfont}{\textsc}
\newcommand{\chairxtitlebibfont}{\textit}
\newcommand{\chairxseriesbibfont}{\textit}

\newcommand{\bibnote}[2]{\nocite{#1}\@namedef{#1chairxnote}{#2}}
\makeatother

\newcommand{\matR}[2]{\field{R}^{{#1} \times {#2}}}

\newcommand{\nablaCan}{\nabla^{\mathrm{can1}}}
\newcommand{\nablaCanSecond}{\nabla^{\mathrm{can2}}}
\newcommand{\nablaHor}{\nabla^{\mathrm{Hor}}}
\newcommand{\nablaHorAlpha}{\nabla^{\mathrm{Hor}, \alpha}}
\newcommand{\nablaAlpha}{\nabla^{\alpha}}

\newcommand{\liegroup}[1]{\mathrm{{#1}}}

\newcommand{\keywords}[1]{\noindent \textbf{Keywords: {#1}}} 
\newcommand{\mathematicssubjectlassification}[1]{\noindent\textbf{2020 Mathematics Subject Classification: {#1}}}

\begin{document}

\selectlanguage{english}

\title{Covariant Derivatives on Homogeneous Spaces - \\
	Horizontal Lifts and Parallel Transport}

\author{Markus Schlarb}
\date{Department of Mathematics, \\
	Julius-Maximilians-Universit{\"a}t W{\"u}rzburg, \\
	Germany \\ \medskip
	\texttt{markus.schlarb@mathematik.uni-wuerzburg.de}\\~\\
	\today}

\maketitle

\begin{abstract}
	We consider invariant covariant derivatives
	on reductive homogeneous spaces corresponding to the well-known
	invariant affine connections.
	These invariant covariant derivatives are expressed in terms of 
	horizontally lifted vector fields on the Lie group.
	This point of view allows for a characterization
	of parallel vector fields along curves.
	Moreover, metric invariant covariant derivatives
	on a reductive homogeneous space equipped with an
	invariant pseudo-Riemannian metric are characterized.
	As a by-product, a new proof for the existence
	of invariant covariant derivatives on reductive homogeneous spaces
	and their the one-to-one correspondence to certain bilinear maps
	is obtained.
\end{abstract}

\keywords{geodesic equation, horizontal lifts, invariant covariant derivatives,
	parallel vector fields along curves, reductive homogeneous spaces}

\medskip

\mathematicssubjectlassification{%
53B05, %Linear and affine connections
53C30, %Differential geometry of homogeneous manifolds
53C22 %Geodesics in global differential geometry
}

\section{Introduction}
\label{sec:introduction}

Reductive homogeneous spaces play a role in a wide range of applications
from mathematical physics to an engineering context.
Without going into details, geodesics and parallel transport are
certainly of interest.
These notions can be defined with respect to 
invariant covariant derivatives
which correspond to the
well-known invariant affine connections
from the literature.
In fact, the existence of invariant affine
connections on a reductive homogeneous space $G / H$ with a fixed
reductive decomposition $\liealg{g} = \liealg{h} \oplus \liealg{m}$
and their one-to-one-correspondence to $\Ad(H)$-invariant bilinear maps
$\liealg{m} \times \liealg{m} \to \liealg{m}$ were proven
in~\cite{nomizu:1954}.

The initial motivation for this text was
to derive a characterization
of parallel vector fields along curves
generalizing~\cite[Lem. 1]{jurdjevic.markina.leite:2023}
to an arbitrary reductive homogeneous space
equipped with some invariant covariant derivative.
In order to obtain such a characterization, given in
Corollary~\ref{corollary:parallel_vectorfields_ODE_on_m} below,
we express an arbitrary invariant covariant
derivative on $G / H$
associated to an $\Ad(H)$-invariant bilinear
map $\liealg{m} \times \liealg{m} \to \liealg{m}$ in terms
of horizontally lifted vector fields on $G$.
This expression generalizes
formulas for the Levi-Civita covariant derivative on $G / H$
in terms of horizontally lifted vector fields
from the literature, where
$G$ is equipped with a bi-invariant metric and $G / H$ is
endowed with a pseudo-Riemannian metric such that $G \to G / H$
is a pseudo-Riemannian submersion.
Indeed, in the proof of~\cite[Lem. 1]{jurdjevic.markina.leite:2023},
a formula for the Levi-Civita covariant derivative on
a pseudo-Riemannian symmetric space 
in terms of horizontally lifted vector fields 
is obtained.
Moreover, a formula for the Levi-Civita covariant
derivative in terms of horizontal vector fields 
is derived in~\cite[Sec. 4.2]{rabenoro.pennec:2023}
for certain homogeneous spaces of compact Lie groups equipped
with bi-invariant metrics.
Here we also mention the recent work~\cite{xu:2022},
where similar questions are independently discussed
in the context of spray geometry.

We now give an overview of this text.
We start with introducing some notations in
Section~\ref{sec:notation_terminology}.
Moreover, in Section~\ref{sec:background},
we recall some facts on reductive homogeneous spaces
and discuss the
principal connections
defined by reductive
decompositions.
After this preparation, we come to
Section~\ref{sec:covariant_derivatives}, where
invariant covariant derivatives are investigated in detail.
In Subsection~\ref{subsec:covariant_derivatives},
we show that an invariant covariant derivative is uniquely
determined by evaluating it on certain fundamental vector fields
of the left action
$G \times G / H \ni (g, g^{\prime} \cdot H) \mapsto (g g^{\prime}) \cdot H \in G / H$.
Afterwards, we express an invariant covariant derivative
$\nablaAlpha$
corresponding to an $\Ad(H)$-invariant bilinear map
$\alpha \colon \liealg{m} \times \liealg{m} \to \liealg{m}$ 
in terms of horizontally lifted vector fields on $G$
as follows.
For two vector fields $X$ and $Y$ on $G / H$
whose horizontal lifts are the vector fields on $G$ denoted by
$\overline{X}$ and $\overline{Y}$, respectively,
we express the horizontal lift $\overline{\nablaAlpha_X Y}$
of $\nablaAlpha_X Y$
in terms of $\overline{X}$ and $\overline{Y}$.
The exact expression for $\overline{\nablaAlpha_X Y}$ is obtained in
Theorem~\ref{theorem:covariant_derivative}.
As a by-product, a new proof for the existence of invariant
covariant derivatives associated to $\Ad(H)$-invariant bilinear maps
$\liealg{m} \times \liealg{m} \to \liealg{m}$ is obtained.
Moreover, the formula from
Theorem~\ref{theorem:covariant_derivative}
is used to derive the curvature of $\nablaAlpha$
in Subsection~\ref{subsec:torsion_curvature}
In addition, we characterize invariant metric covariant derivatives
if $G / H$ is equipped with an invariant pseudo-Riemannian metric.
In Subsection~\ref{subsec:parallel_vector_fields},
we turn our attention to vector fields along curves.
In particular, the expression of $\nablaAlpha$ in terms
of horizontally lifted vector fields from Theorem~\ref{theorem:covariant_derivative}
allows for
characterizing parallel vector fields
along curves on $G / H$
in terms of an ODE on $\liealg{m}$.
In addition, we obtain
a geodesic equation for
the reductive homogeneous space $G / H$ equipped
with an invariant covariant derivative.
If this geodesic equation is specialized to a Lie group
endowed with some left-invariant metric, the well-known geodesic
equation from~\cite[Ap. 2]{arnold:1978} is obtained.
Finally, we discuss the canonical invariant covariant derivatives of first
and second kind which correspond to the 
canonical affine connections of first and second kind
from~\cite[Sec. 10]{nomizu:1954}.

\section{Notations and Terminology}
\label{sec:notation_terminology}

We start with introducing the notation and terminology
that is used throughout
this text.

\begin{notation}
	We follow the convention
	in~\cite[Chap. 2]{oneill:1983}.
	A scalar product is defined as a non-degenerated symmetric bilinear form.
	An inner product is a positive definite symmetric bilinear form.
\end{notation}

Next we introduce some notations concerning differential geometry.
Let $M$ be a smooth (finite-dimensional) manifold.
We denote by $T M$ and $T^* M$ the tangent and cotangent bundle
of $M$, respectively.
For a smooth map $f \colon M \to N$ between manifolds $M$ and $N$,
the tangent map of $f$ is denoted by
$T f \colon TM \to TN$.	
We write $\Cinfty(M)$ for the algebra of smooth real-valued functions on $M$.

Let $E \to M$ be a vector bundle over $M$ with typical fiber $V$.
The smooth sections of $E$ are denoted by $\Secinfty(E)$. 
We write $\End(E) \cong E^* \tensor E$ for the endomorphism bundle of $E$.
Moreover, we denote by $E^{\tensor k}$, $\Sym^k E$ and $\Anti^k E$
the $k$-th tensor power, the $k$-th symmetrized tensor power and the
$k$-th anti-symmetrized tensor power of $E$.
If $T \in \Secinfty\big((T^* M)^{\tensor k} \tensor (T M)^{\tensor \ell} \big)$
is a tensor field on $M$ and $X \in \Secinfty(T M)$ is a vector field
on, we write
$\Lie_X T$ for the Lie derivative.
The pull-back of a smooth function $x \colon N \to \field{R}$
by the smooth map $f \colon M \to N$ is denoted by
$f^* x = x \circ f \colon M \to \field{R}$.
More generally, if $\omega \in \Secinfty\big(\Anti^k(T^* N)\big) \tensor V$
is a differential form taking values in a finite dimensional
$\field{R}$-vector space $V$, its pull-back by $f$ is denoted
by $f^* \omega$.

Concerning the regularity of curves on manifolds,
we use the following convention.

\begin{notation}
	Whenever $c \colon I \to M$ denotes a curve in a manifold $M$
	defined on an interval $I \subseteq \field{R}$,
	we assume for simplicity
	that $c$ is smooth if not indicated otherwise.
	If $I$ is not open, 
	we assume that $c$ can be extended to smooth
	curve defined on
	an open interval $J \subseteq \field{R}$ containing $I$.
	Moreover, we implicitly assume that $0$ is contained in $I$
	if we write $0 \in I$.
\end{notation}

\begin{notation}
	If not indicated otherwise, we use Einstein summation
	convention.
\end{notation}

\section{Background on Reductive Homogeneous Spaces}
\label{sec:background}

In this section, we introduce some more notations and
recall some well-known facts
concerning Lie groups and reductive homogeneous spaces.
Moreover, the principal connection on
the $H$-principal fiber bundle $G \to G / H$
obtained by a reductive decomposition is discussed in detail.

\subsection{Lie groups}

We start with introducing some notations and well-known facts
concerning Lie groups
and Lie algebras.
Let $G$ be a Lie group and denote its Lie algebra by $\liealg{g}$.
The identity of $G$ is usually denoted by $e$.
We write 
\begin{equation}
	\ell_g \colon G \to G,
	\quad 
	h \mapsto \ell_g(h) = g h 
\end{equation}
for the left translation by $g \in G$ and the right translation
by $g \in G$ is denoted by
\begin{equation}
	r_g \colon G \to G, 
	\quad
	h \mapsto r_g(h) = h g .
\end{equation}
The conjugation by an element $g \in G$ is given by
\begin{equation}
	\Conj_{g} \colon G \to G,
	\quad
	h \mapsto \Conj_g(h) 
	=
	(\ell_g \circ r_{g^{-1}})(h)
	=
	(r_{g^{-1}} \circ \ell_g)(h)
	=
	g h g^{-1} 
\end{equation}
and the adjoint representation of $G$ is defined as
\begin{equation}
	\Ad \colon G \to \liegroup{GL}(\liealg{g}),
	\quad
	g \mapsto \Ad_g 
	= \big(\xi \mapsto \Ad_g(\xi) = T_e \Conj_g \xi \big) .
\end{equation}
Moreover, we denote the adjoint representation of $\liealg{g}$ by
\begin{equation}
	\ad \colon \liealg{g} \to \liealg{gl}(\liealg{g}),
	\quad
	\xi \mapsto \big(\eta \mapsto \ad_{\xi}(\eta) = [\xi, \eta] \big) .
\end{equation}
Next we recall~\cite[Def. 19.7]{gallier.quaintance:2020}.
A vector field $X \in \Secinfty(T G)$ is called left-invariant
or right-invariant
if for all $g, k \in G$
\begin{equation}
	T_k \ell_g X(k) = X (\ell_g(k))
	\quad
	\text{ or } \quad
	T_k r_g X(k) = X(r_g(k)),
\end{equation}
respectively,
holds.
For $\xi \in \liealg{g}$, we denote by $\xi^L \in \Secinfty(T G)$
and $\xi^R \in \Secinfty(T G)$
the corresponding left and right-invariant vector fields, respectively,
which are given by 
\begin{equation}
	\xi^L(g) = T_e \ell_g \xi 
	\quad \text{ and } \quad
	\xi^R(g) = T_e r_g \xi,
	\quad g \in G .
\end{equation} 	
We write
\begin{equation}
	\exp \colon \liealg{g} \to G .
\end{equation}
for the exponential map of $G$.

\subsection{Reductive Homogeneous Spaces}
Next we recall some well-known facts on reductive homogeneous spaces
and introduce the notation that is used throughout this text.
We refer to~\cite[Sec. 23.4]{gallier.quaintance:2020}
or~\cite[Chap. 11]{oneill:1983}
for details.

Let $G$ be a Lie group and let $\liealg{g}$ be its Lie algebra.
Moreover, let $H \subseteq G$ a closed subgroup
whose Lie algebra is denoted by
$\liealg{h} \subseteq \liealg{g}$.
We consider the homogeneous space $G / H$.
Then
\begin{equation}
	\label{equation:action_on_homogeneous_space}
	\tau \colon G \times G / H \to G / H, 
	\quad 
	(g, g^{\prime} \cdot H) \mapsto (g g^{\prime}) \cdot H
\end{equation}
is a smooth action of $G$ on $G / H$ from the left,
where $g \cdot H \in G / H$ denotes the coset defined by $g \in G$.
Borrowing the notation from~\cite[p. 676]{gallier.quaintance:2020},
for fixed $g \in G$, the associated diffeomorphism is denote by 
\begin{equation}
	\tau_g \colon G / H \to G / H, 
	\quad 
	g^{\prime} \cdot H 
	\mapsto \tau_g(g^{\prime} \cdot H) 
	=
	(g g^{\prime}) \cdot H .
\end{equation}
In addition, we write
\begin{equation}
	\pr \colon G \to G /H, \quad
	g \mapsto \pr(g) = g \cdot H
\end{equation}
for the canonical projection.

Since reductive homogeneous spaces play a central role in this text,
we recall their definition
from~\cite[Def. 23.8]{gallier.quaintance:2020},
see also~\cite[Sec. 7]{nomizu:1954}
or~\cite[Def. 21, Chap. 11]{oneill:1983}.

\begin{definition}
	Let $G$ be a Lie group and $\liealg{g}$ be its Lie algebra.
	Moreover, let
	$H \subseteq G$ be a closed subgroup and denote its Lie
	algebra by $\liealg{h} \subseteq \liealg{g}$.
	Then the homogeneous space $G / H$ is called reductive 
	if there exists a subspace
	$\liealg{m} \subseteq \liealg{g}$ 
	such that $\liealg{g} = \liealg{h} \oplus \liealg{m}$ is fulfilled and
	\begin{equation}
		\Ad_h(\liealg{m}) \subseteq \liealg{m}
	\end{equation}
	holds for all $h \in H$.
\end{definition}
Following~\cite[Prop. 23.22]{gallier.quaintance:2020}, we 
recall a well-known property of the isotropy representation of
a reductive homogeneous space.
This is the next lemma.

\begin{lemma}
	\label{lemma:isotropy_representation_equivalence}
	The isotropy representation 
	of a reductive homogeneous space $G / H$
	with reductive decomposition
	$\liealg{g} = \liealg{h} \oplus \liealg{m}$
	\begin{equation}
		H \ni h 
		\mapsto T_{\pr(e)} \tau_h 
		\in \liegroup{GL}\big(T_ {\pr(e)} G / H \big)
	\end{equation}
	is equivalent to the representation
	\begin{equation}
		H \to \liegroup{GL}(\liealg{m}),
		\quad
		h \mapsto \Ad_h\at{\liealg{m}} = \big(X \mapsto \Ad_h(X) \big),
	\end{equation}
	i.e.
	\begin{equation}
		T_{\pr(e)}  \tau_h  \circ  T_e \pr\at{\liealg{m}} 
		=
		T_e \pr \circ \Ad_h\at{\liealg{m}}
	\end{equation}
	is fulfilled for all $h \in H$.
\end{lemma}

\begin{notation}
	Let $\liealg{g} = \liealg{h} \oplus \liealg{m}$ be
	a reductive decomposition of $\liealg{g}$.
	Then the projection onto $\liealg{m}$ 
	whose kernel is given by $\liealg{h}$ is denoted by
	$\pr_{\liealg{m}} \colon \liealg{g} \to \liealg{m}$.
	We write 
	$\pr_{\liealg{h}} \colon \liealg{g} \to \liealg{h}$ for the projection whose kernel is given by $\liealg{m}$.
	Moreover, we write for $\xi \in \liealg{g}$
	\begin{equation}
		\xi_{\liealg{m}} = \pr_{\liealg{m}}(\xi)
		\quad  \text{ and } \quad
		\xi_{\liealg{h}} = \pr_{\liealg{h}}(\xi) .
	\end{equation}
\end{notation}
A scalar product $\langle \cdot, \cdot \rangle \colon \liealg{m} \times \liealg{m} \to \field{R}$
is called $\Ad(H)$-invariant if
\begin{equation}
	\big\langle \Ad_h(X), \Ad_h(Y) \big\rangle
	=
	\big\langle X, Y \big\rangle
\end{equation}
holds for all $h \in H$ and $X, Y \in \liealg{m}$,
see e.g~\cite[p. 301]{oneill:1983}
or~\cite[Sec. 23.4]{gallier.quaintance:2020}
for the positive definite case.
Reformulating and adapting~\cite[Def. 23.5]{gallier.quaintance:2020},
we call a pseudo-Riemannian metric
$\langle \! \langle \!  \cdot, \cdot \rangle \! \rangle
\in \Secinfty\big( \Sym^2 T^* (G / H)\big)$
invariant if
\begin{equation}
	\langle \! \langle v_p, w_p \rangle \! \rangle_{p}
	=
	\langle \! \langle T_p \tau_g v_p , T_p \tau_g w_p \rangle \! \rangle_{\tau_g(p)},
	\quad
	p \in G / H, \quad v_p, w_p \in T_p (G  / H)
\end{equation}
holds for all $g \in G$.
In the next lemma which is taken 
from~\cite[Chap. 11, Prop. 22]{oneill:1983}, see
also~\cite[Prop. 23.22]{gallier.quaintance:2020}
for the Riemannian case,
invariant metrics on $G / H$ are related to
$\Ad(H)$-invariant scalar products on $\liealg{m}$.

\begin{lemma}
	By requiring the linear isomorphism
	$T_e \pr\at{m} \colon \liealg{m} \to T_{\pr(e)} (G / H)$
	to be an isometry,
	there is a one-to-one correspondence between
	$\Ad(H)$-invariant scalar products on $\liealg{m}$ and
	invariant pseudo-Riemannian metrics on $G / H$.
\end{lemma}
Naturally reductive homogeneous spaces are special
reductive homogeneous spaces.
We recall their definition from~\cite[Chap. 11, Def. 23]{oneill:1983}.

\begin{definition}
	\label{definition:naturally_reductive_space}
	Let $G / H$ be a reductive homogeneous space equipped with
	an invariant pseudo-Riemannian metric corresponding to the
	$\Ad(H)$-invariant scalar product
	$\langle \cdot, \cdot \rangle \colon \liealg{m} \times \liealg{m} \to \field{R}$.
	Then $G /H$ is called a naturally reductive homogeneous space if
	\begin{equation}
		\big\langle [X, Y]_{\liealg{m}}, Z \big\rangle 
		= 
		\big\langle X, [Y, Z]_{\liealg{m}} \big\rangle
	\end{equation}
	holds for all $X, Y, Z \in \liealg{m}$.
\end{definition}
The following lemma can be considered as a generalization
of~\cite[Prop. 23.29 (1)-(2)]{gallier.quaintance:2020}
to pseudo-Riemannian metrics and
Lie groups which are not necessarily connected.

\begin{lemma}
	\label{lemma:normal_naturally_reductive_is_naturally_reductive}
	Let $G$ be a Lie group and denote by $\liealg{g}$ its Lie algebra.
	Moreover, let $G$ be equipped with a bi-invariant metric
	and let
	$\langle \cdot, \cdot \rangle \colon \liealg{g} \times \liealg{g}
	\to \field{R}$
	be the corresponding $\Ad(G)$-invariant scalar product.
	Moreover, let $H \subseteq G$ be a closed subgroup such that
	its Lie algebra 
	$\liealg{h} \subseteq \liealg{g}$
	is non-degenerated with respect to $\langle \cdot, \cdot \rangle$.
	Then $G / H$ is a reductive homogeneous space with reductive 
	decomposition $\liealg{g} = \liealg{h} \oplus \liealg{m}$, 
	where $\liealg{m} = \liealg{h}^{\perp}$
	is the orthogonal complement of $\liealg{h}$
	with respect to $\langle \cdot, \cdot \rangle$.
	Moreover, if $G / H$ is equipped with the invariant metric
	corresponding to the scalar product on $\liealg{m}$ that is obtained by restricting $\langle \cdot, \cdot \rangle$ to $\liealg{m}$,
	the reductive homogeneous space $G / H$ is naturally reductive.
	\begin{proof}
		The claim can be proven analogously to the proof 
		of~\cite[Prop. 23.29 (1)-(2)]{gallier.quaintance:2020}
		by taking the assumption
		$\liealg{h} \oplus \liealg{h}^{\perp} 
		=
		\liealg{h} \oplus \liealg{m}
		= \liealg{g}$
		into account.
	\end{proof}
\end{lemma}

\begin{remark}
	\label{remark:normal_naturally_recductive}
	Inspired by the terminology
	in~\cite[Sec. 23.6, p. 710]{gallier.quaintance:2020},
	we refer to the naturally reductive homogeneous spaces from
	Lemma~\ref{lemma:normal_naturally_reductive_is_naturally_reductive}
	as normal naturally reductive homogeneous spaces.
\end{remark}

We end this subsection with considering another special class of
reductive homogeneous spaces.
To this end, we state the following
definition which can be found in~\cite[p. 209]{helgason:1978}.

\begin{definition}
	Let $G$ be a connected Lie group and let $H$ be a closed subgroup.
	Then $(G, H)$ is called a symmetric pair if there exists a
	smooth involutive
	automorphism $\sigma \colon G \to G$, i.e.
	an automorphism of Lie groups fulfilling $\sigma^2 = \sigma$,
	such that $(H_{\sigma})_0 \subseteq H \subseteq H_{\sigma}$ holds.
	Here $H_{\sigma}$ denotes the set of fixed points of $\sigma$ and
	$(H_{\sigma})_0$ denotes the connected component of $H_{\sigma}$
	containing the identity $e \in G$.
\end{definition}

Inspired by the terminology used
in~\cite[Def. 23.13]{gallier.quaintance:2020}, we refer to the 
triple $(G, H, \sigma)$ as symmetric pair, as well, where $(G, H)$ is a symmetric pair
with respect to the involutive automorphism
$\sigma \colon G \to G$.
These symmetric pairs lead to reductive homogeneous spaces
which are called symmetric homogeneous spaces if a certain ``canonical''
reductive decomposition is chosen, see
e.g.~\cite[Sec. 14]{nomizu:1954}.
Note that the definition 
in~\cite[Sec. 14]{nomizu:1954}
does not require an invariant pseudo-Riemannian metric on $G / H$.

The next lemma,
see e.g.~\cite[Sec. 14]{nomizu:1954},
shows that a symmetric homogeneous space is a reductive homogeneous
space with respect to the so-called canonical
reductive decomposition.
Here we also refer
to~\cite[Prop. 23.33]{gallier.quaintance:2020} for a proof.

\begin{lemma}
	\label{lemma:symmetric_pair_canonical_reductive_decomposition}
	Let $(G , H, \sigma)$ be a symmetric pair
	and define the subspaces of $\liealg{g}$ by
	\begin{equation}
		\liealg{h} = \{ X \in \liealg{g} \mid T_e \sigma X = X\}
		\subseteq \liealg{g}
		\quad  \text{ and } \quad
		\liealg{m} = \{X \in \liealg{g} \mid T_e \sigma X = - X \} 
		\subseteq \liealg{g} .
	\end{equation}
	Then $\liealg{g} = \liealg{h} \oplus \liealg{m}$ is a
	reductive decomposition of $\liealg{g}$ turning $G / H$
	into a reductive homogeneous space.
	Moreover, the inclusion
	\begin{equation}
		[\liealg{m} , \liealg{m}] \subseteq \liealg{h}
	\end{equation}
	is fulfilled.
\end{lemma}
Next we define symmetric homogeneous spaces and canonical reductive
decompositions following~\cite[Sec. 14]{nomizu:1954}.

\begin{definition}
	Let $(G, H, \sigma)$ be a symmetric pair.
	Then the reductive decomposition
	$\liealg{g} = \liealg{h} \oplus \liealg{m}$
	from Lemma~\ref{lemma:symmetric_pair_canonical_reductive_decomposition}
	is called canonical reductive decomposition.
	Moreover, the reductive homogeneous space $G / H$ with
	the reductive decomposition
	$\liealg{g} = \liealg{h} \oplus \liealg{m}$ 
	is called symmetric homogeneous space.
\end{definition}
For pseudo-Riemannian symmetric spaces we state the next remark
following~\cite[Chap. 11, p. 317]{oneill:1983}, see
also~\cite[Sec. 23.8]{gallier.quaintance:2020}
for the Riemannian case.

\begin{remark}
	\label{remark:pseudo-Riemannian-symmetric-naturally-reductive}
	Let $(G, H, \sigma)$ be symmetric pair and let $G / H$ be the
	associated symmetric homogeneous space with canonical reductive
	decomposition $\liealg{g} = \liealg{h} \oplus \liealg{m}$.
	Let $G / H$ be equipped with an invariant pseudo-Riemannian metric
	and let
	$\langle \cdot, \cdot \rangle \colon 
	\liealg{m} \times \liealg{m} \to \field{R}$
	be the associated $\Ad(H)$-invariant scalar product.
	Then $G / H$ is a naturally reductive homogeneous space
	since $[\liealg{m}, \liealg{m}] \subseteq \liealg{h}$ implies
	that the condition on the scalar product
	$\langle \cdot, \cdot \rangle$ from
	Definition~\ref{definition:naturally_reductive_space}
	is always satisfied.
	In the sequel, we refer to symmetric homogeneous spaces equipped 
	with an invariant pseudo-Riemannian metric as pseudo-Riemannian
	symmetric homogeneous space
	or pseudo-Riemannian symmetric spaces, for short.
\end{remark}

\subsection{Reductive Decompositions and Principal Connections}

In this section, we consider $G$ as a $H$-principal fiber bundle over
$G / H$ and discuss certain principal connections
on $\pr \colon G \to G / H$.
For general properties of principal fiber bundles and connections,
we refer to~\cite[Sec. 18-19]{michor:2008}
and~\cite[Sec. 1.1-1.3]{rudolph.schmidt:2017}.

Let $G$ be a Lie group and $H \subseteq G$ be a closed subgroup.
It is well-known that
$\pr \colon G \to G / H$ is a $H$-principle
fiber bundle, see e.g. \cite[Sec. 18.15]{michor:2008},
where the base is the homogeneous space $G / H$.
The $H$-principal action on $G$ is denoted by 
\begin{equation}
	\label{equation:H_principal-action_on_G_mod_H}
	\racts \colon G \times H \to G, 
	\quad 
	(g, h) \mapsto g \racts h
	= r_h(g) 
	= \ell_g(h) 
	= g h,
\end{equation}
if not indicated otherwise.
We now assume that $G / H$ is a reductive homogeneous space
and the reductive decomposition 
$\liealg{g} = \liealg{h} \oplus \liealg{m}$
is fixed.
This reductive decomposition
can be used to obtain a
principal connection on $\pr \colon G \to G / H$,
see~\cite[Thm. 11.1]{kobayashi.nomizu:1963a}.
Although this fact is well-known, we provide a detailed proof
in order to keep this text more self-contained.
To this end,
we recall a well-known fiber-wise expression for the vertical bundle
of $\pr \colon G \to G / H$ which follows for example 
from~\cite[Sec. 18.18]{michor:2008}.
We have for fixed $g \in G$
by~\cite[Lem. 1.3.1]{rudolph.schmidt:2017},
see also~\cite[Sec. 18.8]{michor:2008}
\begin{equation}
	\label{equation:vertical_bundle_reductive_space_description}
	\Ver(G)_g 
	=
	\big\{ \tfrac{\D}{\D t} \big( g \racts \exp(t \eta) \big) \at{t = 0} \mid \eta \in \liealg{h} \big\}
%	= 
%	\big\{ \tfrac{\D}{\D t} \ell_g (\exp(t \eta)) \at{t = 0} \mid \eta \in \liealg{h}  \big\}
	= 
	(T_e \ell_g )\liealg{h} .
\end{equation}
The next proposition provides explicit formulas for the principal connection
and the associated principal connection one-form on $G \to G / H$ defined by a
reductive decomposition.

\begin{proposition}
	\label{proposition:principal_connection_reductive_homogeneous_space}
	Consider $\pr \colon G \to G / H$ as a $H$-principal fiber bundle,
	where $G / H$ is a reductive homogeneous space with a fixed reductive
	decomposition $\liealg{g} = \liealg{h} \oplus \liealg{m}$
	and define $\Hor(G) \subseteq T G$ fiber-wise by 
	\begin{equation}
		\Hor(G)_g = (T_e \ell_g) \liealg{m}, \quad g \in G.
	\end{equation}
	Then $\Hor(G)$ is a subbundle of $TG$ defining a horizontal bundle
	on $T G$, i.e. a complement of the vertical bundle
	$\Ver(G) = \ker (T \pr) \subseteq TG$
	which yields a principal connection on $\pr \colon G \to G / H$. 
	This principal connection $\mathcal{P} \in \Secinfty\big(\End(T G) \big)$ 
	corresponding to $\Hor(G)$ is given by
	\begin{equation}
		\label{equation:principal_connection_G_mod_H_defined_by_reductive_split}
		\mathcal{P}\at{g}(v_g) = T_e \ell_g \circ \pr_{\liealg{h}} \circ (T_e \ell_g)^{-1} v_g,
		\quad 
		g \in G, \quad v_g \in T_g G .
	\end{equation}
	The corresponding connection one-form $\omega \in \Secinfty(T^* G) \tensor \liealg{h}$ reads 
	\begin{equation}
		\omega\at{g}(v_g) = \pr_{\liealg{h}} \circ (T_e \ell_g)^{-1} v_g 
	\end{equation}
	for $g \in G$ and $v_g \in T_g G$.	
	\begin{proof}
		Although, this statement is well-known,
		see e.g.~\cite[Thm. 11.1]{kobayashi.nomizu:1963a},
		we provide a proof, nevertheless.
		Indeed, $\Hor(G)$ is a complement of the vertical bundle
		$\Ver(G) = \ker(T \pr) \subseteq T G$
		due to 
		$\liealg{g} = \liealg{h} \oplus \liealg{m}$
		implying $T G = \Ver(G) \oplus \Hor(G)$
		as desired.
		Moreover, $\mathcal{P}$ defined
		by~\eqref{equation:principal_connection_G_mod_H_defined_by_reductive_split}
		is clearly a smooth
		endomorphism of the vector bundle $T G \to G$, i.e.
		$\mathcal{P} \in \Secinfty\big(\End(T G) \big)$.
		In addition, $\mathcal{P}^2 = \mathcal{P}$ is obviously fulfilled.
		Moreover, one has
		$\image(\mathcal{P}) 
		=
		\ker( T \pr) 
		=
		\Ver(G)$
		and
		$\ker(\mathcal{P}) = \Hor(G)$
		showing that $\mathcal{P}$
		is the connection corresponding to the horizontal bundle $\Hor(G)$. 
		
		We now show that $\omega$ is the connection one-form
		corresponding to $\mathcal{P}$ by using the
		correspondence from~\cite[Sec. 19.1]{michor:2008}.
		Let $\eta \in \liealg{h}$ and denote by $\eta_G$
		the corresponding fundamental vector field, i.e.
		we have for $g \in G$
		\begin{equation*}
			\eta_G(g) 
			= \tfrac{\D}{\D t}
			\big( g \racts \exp(t \eta) \big)\at{t = 0}
			= \tfrac{\D}{\D t} \ell_g\big( \exp(t \eta) \big) \at{t = 0}
			= T_e \ell_g \eta.
		\end{equation*} 
		By this notation, one obtains for $v_g \in T_g G$
		\begin{equation*}
			\big(\omega\at{g}(v_g)\big)_G(g) 
			=
			T_e \ell_g \big( \omega\at{g}(v_g) \big)
			=
			T_e \ell_g \big( \pr_{\liealg{h}} \circ (T_e \ell_g)^{-1} (v_g)\big)
			= 
			\mathcal{P}\at{g}(v_g) .
		\end{equation*}
		Moreover, we have
		\begin{equation*}
			\omega\at{g}\big(\eta_G(g) \big) 
			=
			\big( \pr_{\liealg{h}} \circ (T_e \ell_g)^{-1} \big) T_e \ell_g \eta
			=
			\pr_{\liealg{h}}(\eta)
			=
			\eta
		\end{equation*}
		for all $\eta \in \liealg{h}$
		proving that $\omega \in \Secinfty(T^* G) \tensor \liealg{h}$
		is the connection one-from corresponding
		to $\mathcal{P} \in \Secinfty\big( \End(T G) \big)$.
		
		It remains to show that $\mathcal{P}$ is a principal connection.
		By~\cite[Sec. 19.1]{michor:2008}
		this is equivalent to showing that
		$\omega$ has the equivariance property
		\begin{equation*}
			\big((\cdot \racts h)^* \omega\big)\at{g}(v_g)
			=
			\Ad_{h^{-1}}\big(\omega\at{g}(v_g) \big)
		\end{equation*}
		for all $h \in H$, $g \in G$ and $v_g \in T_g G$,
		where $(\cdot \racts h)^* \omega$ denotes the pull-back of $\omega$
		by $(\cdot \racts h) \colon P \ni p \mapsto p \racts h \in P$.
		Since
		$\liealg{g} = \liealg{h} \oplus \liealg{m}$ is a reductive
		decomposition,
		we obtain for $h \in H$ and $\xi \in \liealg{g}$
		\begin{equation}
			\label{equation:proposition_principal_connection_reductive_homogeneous_space_pr_h_Ad_commute}
			(\pr_{\liealg{h}} \circ \Ad_h)(\xi) 
			=
			\pr_{\liealg{h}} \big(\Ad_h(\xi_{\liealg{h}}) 
			+ \Ad_h(\xi_{\liealg{m}}) \big)
			=
			\Ad_{h}(\xi_{\liealg{h}}) 
			=
			(\Ad_h \circ \pr_{\liealg{h}})(\xi) .
		\end{equation}
		Using~\eqref{equation:proposition_principal_connection_reductive_homogeneous_space_pr_h_Ad_commute}
		and the chain-rule,
		we compute for $h \in H$, $g \in G$ and $v_g \in T_g G$
		\begin{equation*}
			\begin{split}
				\big((\cdot \racts h)^* \omega\big)\at{g}(v_g)
				&=
				\omega\at{g \racts h}\big(T_g (\cdot \racts h) v_g \big) \\
				&=
				\omega\at{gh} \big(T_g r_h v_g \big) \\
				&=
				\big(\pr_{\liealg{h}} \circ (T_e \ell_{g h})^{-1} \big) T_g r_h v_g \\
				&= 
				\pr_{\liealg{h}} \circ T_{g h} \ell_{h^{-1} g^{-1}} \circ  T_g r_h v_g \\
				&=
				\pr_{\liealg{h}} \circ T_g (\ell_h^{-1} \circ \ell_g^{-1} \circ r_h) v_g \\
				&= 
				\pr_{\liealg{h}} \circ T_g (\ell_{h^{-1}} \circ r_h \circ \ell_{g^{-1}} ) v_g \\
				&=
				\pr_{\liealg{h}} \circ T_e (\ell_{h^{-1}} \circ r_h) \circ T_g \ell_{g^{-1}}  v_g \\
				%		 			&= 
				%	 				\pr_{\liealg{h}} \circ T_e (\ell_{h^{-1}} \circ r_h) \circ (T_e \ell_g)^{-1} v_g \\
				&=
				\pr_{\liealg{h}} \circ T_e \Conj_{h^{-1}} \circ (T_e \ell_g)^{-1} v_g \\
				&= 
				\pr_{\liealg{h}} \circ \Ad_{h^{-1}} \circ  (T_e \ell_g)^{-1} v_g \\
				&=
				\Ad_{h^{-1}} \circ \pr_{\liealg{h}} \circ (T_e \ell_g)^{-1} v_g \\
				&= 
				\Ad_{h^{-1}}\big(\omega\at{g}(v_g) \big).
			\end{split}
		\end{equation*}
		Hence $\omega$ is the connection one-form corresponding to the
		principal $\mathcal{P}$.
	\end{proof}
\end{proposition}
By~\cite[Prop. 23.23]{gallier.quaintance:2020},
adapted to the pseudo-Riemannian case,
we obtain the following remark concerning pseudo-Riemannian
reductive homogeneous spaces.
\begin{remark}
	\label{remark:principal_connection_normal_pseudo_Riemannian_reductive_homogeneous_space}
	Let $G / H$ be a reductive homogeneous
	space with reductive decomposition
	$\liealg{g} = \liealg{h} \oplus \liealg{m}$
	endowed with an invariant pseudo-Riemannian metric
	corresponding to the $\Ad(H)$-invariant scalar product
	$\langle \cdot, \cdot \rangle \colon \liealg{m} \times \liealg{m} \to \field{R}$.
	By~\cite[Prop. 23.23]{gallier.quaintance:2020},
	which clearly extends to the pseudo-Riemannian case, 
	the scalar product
	$\langle \cdot, \cdot \rangle$ on $\liealg{m}$
	can be extended to a scalar product
	$\langle \cdot, \cdot \rangle_{\liealg{g}}$
	on $\liealg{g}$ such that
	$\liealg{m} = \liealg{h}^{\perp}$ is fulfilled.
	Then $\pr \colon G \to G / H$ becomes
	a pseudo-Riemannian submersion
	by~\cite[Prop. 23.23]{gallier.quaintance:2020},
	where $G$ is equipped with the left-invariant metric
	defined by $\langle \cdot, \cdot \rangle_{\liealg{g}}$.
	Obviously, the horizontal bundle defined by
	$\Hor(G) = \Ver(G)^{\perp}$
	yields the connection on $G$ which coincides with the
	principal connection from
	Proposition~\ref{proposition:principal_connection_reductive_homogeneous_space}
	defined by the reductive decomposition
	$\liealg{g} = \liealg{h} \oplus \liealg{m}$.
\end{remark}

\section{Invariant Covariant Derivatives}
\label{sec:covariant_derivatives}

In this section,
we consider the invariant covariant derivatives
on a reductive homogenoeus space $G / H$
that correspond to the invariant affine connections investigated
in~\cite{nomizu:1954}.
These invariant covariant derivatives are expressed in terms of
horizontally lifted vector fields
yielding another proof for their existence.
In particular, this expression is used to 
characterize parallel vector fields along curves
in terms of an ODE on $\liealg{m}$.

Throughout this subsection, we use the following notation.

\begin{notation}
	If not indicated otherwise,
	we denote by $G / H$ a reductive homogeneous space with a fixed
	reductive decomposition $\liealg{g} = \liealg{h} \oplus \liealg{m}$.
\end{notation}

\subsection{Invariant Covariant Derivatives}
\label{subsec:covariant_derivatives}

We start with introducing the notion of an invariant
covariant derivative on a reductive homogeneous space.
In view of the one-to-one correspondence of covariant derivatives
and affine connections, see
Remark~\ref{remark:affine_connection_covariant_derivatives} below,
the next definition can be seen as a reformulation
of~\cite[Eq. (2.3) and Sec. 8, p. 43]{nomizu:1954}.

\begin{definition}
	\label{definition:invariant_covariant_derivative_on_G_H_reductive}
	Let $G / H$ be a homogeneous space.
	Then a covariant derivative
	\begin{equation}
		\nabla \colon \Secinfty\big(T (G / H)\big) \times \Secinfty\big(T (G / H) \big) \to \Secinfty\big( T(G / H)\big)
	\end{equation}
	on $G / H$
	is called $G$-invariant, or invariant for short,
	if 
	\begin{equation}
		\nabla_X Y = (\tau_{g^{-1}})_* \big( \nabla_{ (\tau_g)_* X}  (\tau_g)_* Y\big)
	\end{equation}
	holds for all $g \in G$ and $X, Y \in \Secinfty\big(T (G / H)\big)$,
	where
	$(\tau_g)_* X$ denotes the well-known push-forward
	of $X$ by the diffeomorphism
	$\tau_g \colon G / H \to G / H$ given by
	$(\tau_g)_* X =  T  \tau_g \circ X \circ \tau_{g^{-1}}$.
\end{definition}
Obviously, for a fixed $g \in G$ the push-forward
$(\tau_g)_* X  
=
T \tau_g  \circ X \circ \tau_{g^{-1}}$
of a vector field
$X \in \Secinfty\big(T (G / H) \big)$ by
$\tau_g \colon G / H \to G / H$ is point-wise given by
\begin{equation}
	\label{equation:push_forward_covariant_derivative_on_G_H}
	\big((\tau_g)_* X \big) (\pr(k))
	=
	T_{\tau_{g^{-1}(\pr(k))}} \tau_g  X \big(\tau_{g^{-1}}(\pr(k)) \big), \quad \pr(k) \in G / H  .
\end{equation}

In the next remark, we relate the notion of affine connections
from~\cite{nomizu:1954}
to covariant derivatives.

\begin{remark}
	\label{remark:affine_connection_covariant_derivatives}
	Let $M$ be a manifold and denote by
	$\End_{\Cinfty(M)}\big(\Secinfty(T M)\big)$ the
	endomorphisms of the $\Cinfty(M)$-module $\Secinfty(T M)$.
	An affine connection is defined in~\cite{nomizu:1954}
	as a map
	\begin{equation}
		\label{equation:remark_affine_connection_covariant_derivative}
		t \colon \Secinfty(T M) \ni X \mapsto t (X) \in \End_{\Cinfty(M)}\big(\Secinfty(T M) \big)
	\end{equation}
	such that
	\begin{equation}
		\label{equation:remark_affine_connection_covariant_derivative_defining_properties}
		t(X_1 + X_1) = t(X_1) + t(X_2)
		\quad
		\text{ and }
		\quad
		t(f X)(Y) = ft(X)(Y) + ( \Lie_X f ) t(X)(Y)
	\end{equation}
	holds
	for all $X_1, X_2, X, Y \in \Secinfty(T M)$.
	As pointed out in~\cite[Sec. 4.5]{berestovskii.nikonorov:2020},
	an affine connection $t \colon  \Secinfty(T M) \to \End_{\Cinfty(M)}\big(\Secinfty(T M) \big)$
	defines a covariant derivative
	$\nabla \colon \Secinfty(T M) \times \Secinfty(T M) \to \Secinfty(T M)$
	on $TM$ by
	\begin{equation}
		\label{equation:affine_connection_covariant_derivative_relation}
		\nabla_X Y = t(Y)(X),
		\quad X, Y \in \Secinfty( T M) .
	\end{equation}
	Obviously, the converse is also true.
	Given a covariant derivative $\nabla$ on $TM$,
	Equation~\eqref{equation:affine_connection_covariant_derivative_relation}
	yields
	an affine connection.
\end{remark}

In the sequel, we discuss the invariant covariant derivatives on $G / H$
corresponding to the invariant affine connections on $G / H$
from~\cite[Thm. 8.1]{nomizu:1954}.
This correspondence is made precise in
Proposition~\ref{proposition:one-to-one-correspondence-invariant-affine-connection-covariant-derivatives}, below.

We first recall the notion of an $\Ad(H)$-invariant bilinear map
from~\cite[Sec. 8]{nomizu:1954}. 

\begin{definition}
	\label{definition:Ad_H_invariant_bilinear_map}
	Let $G / H$ be a reductive homogeneous space with reductive decomposition $\liealg{g} = \liealg{h} \oplus \liealg{m}$.
	Then the bilinear map
	$\alpha \colon \liealg{m} \times \liealg{m} \to \liealg{m}$ 
	is called $\Ad(H)$-invariant if
	\begin{equation}
		\Ad_h\big(\alpha(X, Y) \big) 
		=
		\alpha\big( \Ad_h(X), \Ad_h(Y) \big)
	\end{equation}
	holds for all $X, Y \in \liealg{m}$ and $h  \in H$.
	More generally,
	for $\ell \in \field{N}$,
	we call a $\ell$-linear 
	map $\alpha \colon \liealg{m}^{\ell} \to \liealg{m}$
	$\Ad(H)$-invariant if
	\begin{equation}
		\Ad_h\big(\alpha(X_1, \ldots, X_\ell) \big) 
		=
		\alpha\big( \Ad_h(X_1), \ldots, \Ad_h(X_\ell) \big)		
	\end{equation}
	holds for all $X_1, \ldots, X_\ell \in \liealg{m}$ and $h \in H$.
\end{definition}

\begin{remark}
	As we have already pointed out in the introduction,
	the one-to-one correspondence
	between invariant affine connections
	and $\Ad(H)$-invariant bilinear maps
	$\liealg{m} \times \liealg{m} \to \liealg{m}$
	is well-known by~\cite[Thm. 8.1]{nomizu:1954}.
	Nevertheless, the discussion in this text differs from
	the discussion in~\cite{nomizu:1954}.
	Inspired by~\cite[Sec. 23.6]{gallier.quaintance:2020},
	we consider invariant covariant derivatives evaluated
	at the fundamental
	vector fields of the action $\tau \colon G \times G / H \to G / H$
	at the point $\pr(e)$ which already determines them uniquely.
	Moreover, we express invariant covariant derivatives on $G / H$
	in terms of
	horizontally lifted vector fields on $G$.
	Beside yielding another proof for the existence of an
	invariant covariant derivative associated with an
	$\Ad(H)$-invariant bilinear map
	$\liealg{m} \times \liealg{m} \to \liealg{m}$,
	this point of view allows in particular for an easy characterization
	of parallel vector fields, see
	Subsection~\ref{subsec:parallel_vector_fields}
	below.
\end{remark}

\subsubsection{Invariant Covariant Derivatives evaluated on Fundamental Vector Fields}

Before we continue with
considering invariant covariant derivatives, we
take a closer look on the fundamental vector fields on $G / H$ associated
with the action $\tau \colon G \times G / H \to G / H$
from~\eqref{equation:action_on_homogeneous_space}.	
Let $X \in \liealg{g}$.
The fundamental vector field $X_{G / H} \in \Secinfty\big(T (G / H)\big)$
associated with $X$
is defined by
\begin{equation}
	X_{G / H}(\pr(g)) = \tfrac{\D}{\D t} \tau_{\exp(t X)}(\pr(g)) \at{ t = 0}
\end{equation}
for $\pr(g)  \in G / H$ with $g \in G$.
In the next lemma, we state some properties of $X_{G / H}$.
Note that its third claim is well-known.

\begin{lemma}
	\label{lemma:reductive_homogeneous_space_fundamental_vfs}
	Let $G / H$ be a reductive homogeneous space with reductive decomposition $\liealg{g} = \liealg{h} \oplus \liealg{m}$.
	Moreover, let
	$X \in \liealg{m}$, let
	$\{A_1, \ldots, A_N\} \subseteq \liealg{m}$ 
	be a basis of $\liealg{m}$
	and let $\{A^1, \ldots, A^N\} \subseteq \liealg{m}^*$
	be its dual basis.
	Let $A_i^L \in \Secinfty\big( \Hor( G)\big)$
	denote the left-invariant vector field on $G$
	defined by $A_i$ for $i \in \{1, \ldots, N\}$.
	Then the following assertions are fulfilled:
	\begin{enumerate}
		\item
		\label{item:lemma_reductive_homogeneous_space_fundamental_vfs_horizontal_lifts}
		The horizontal lift of $X_{G / H}$ is given by
		\begin{equation}
			\overline{X_{G / H}}(g) = A^i \big( \Ad_{g^{-1}}(X)\big) A_i^L (g)
		\end{equation}
		for all $g \in G$.
		\item
		\label{lemma:reductive_homogeneous_space_fundamental_vfs_lie_derivative_horizontal_lift_at_identity}
		Let $Y \in \liealg{m}$ and define
		the smooth functions
		$y^j \colon G \ni g \mapsto y^j(g) = A^j\big(\Ad_{g^{-1}}(Y) \big) \in \field{R}$,
		where $j \in \{1, \ldots, N\}$.
		Then one has
		\begin{equation}
			\big(\Lie_{\overline{X_{G / H}}} y^j \big)(e) A_j^L(e)
			=
			- [X, Y]_{\liealg{m}} .
		\end{equation}
		\item
		\label{item:lemma_reductive_homogeneous_space_fundamental_vfs_push_forward_by_canonical_action}
		One has
		\begin{equation}
			(\tau_g)_* X_{G / H}(\pr(k)) = \big( \Ad_g(X)\big)_{G / H}(\pr(k))
		\end{equation}
		for all $g \in G$ and $\pr(k) \in G / H$.
	\end{enumerate}
	\begin{proof}
		We first show Claim~\ref{item:lemma_reductive_homogeneous_space_fundamental_vfs_horizontal_lifts}.
		To this end, we compute for $g \in G$
		\begin{equation}
			\begin{split}
				\label{equation:lemma_reductive_homogeneous_space_fundamental_vfs_horizontal_lift_right_invariant}
				X_{G / H}(\pr(g)) 
%				&=
%				\tfrac{\D}{\D t} \tau_{\exp(t X)}(\pr(g)) \at{t = 0} \\
				&=
				\tfrac{\D}{\D t} (\tau_{\exp(t X)} \circ \pr)(g) \at{t = 0} \\
				&=
				\tfrac{\D}{\D t} (\pr \circ \ell_{\exp(t X)} )(g) \at{t = 0} \\
				&=
				\tfrac{\D}{\D t} \pr \circ r_g (\exp(t X)) \at{t = 0} \\
				&=
				T_g \pr \circ T_e r_g X \\
				&=
				T_g \pr \circ X^R(g)
			\end{split}
		\end{equation}
		showing that $X_{G / H}$ is $\pr$-related to
		the right-invariant vector field $X^R$.
		Next we express $X^R$ in terms of left-invariant vector fields.
		Let $g \in G$.
		We now compute
		\begin{equation}
			\label{equation:lemma_reductive_homogeneous_space_fundamental_vfs_right_left_invariant_relation}
			\begin{split}
				X^R(g) 
				&=
				\big( T_e \ell_g \circ (T_e \ell_g)^{-1} \big) T_e r_g X \\
				&=
				T_e \ell_g \circ T_g \ell_{g^{-1}} \circ T_e r_g X  \\
				&= 
				T_e \ell_g \circ T_e (\ell_{g^{-1}} \circ r_g) X \\
				&=
				T_e \ell_g \circ T_e \Conj_{g^{-1}} X \\
				&= 
				T_e \ell_g \Ad_{g^{-1}} (X) \\
				&= 
				\big(\Ad_{g^{-1}}(X)\big)^L(g) .
			\end{split}		
		\end{equation}
		Let $\mathcal{P}$ be the principal connection from
		Proposition~\ref{proposition:principal_connection_reductive_homogeneous_space}.
		Then the horizontal lift of $X_{G / H}$ is given
		by $\overline{X_{G / H}} = (\id_{T G} - \mathcal{P}) \circ X^R$
		because of $X_{G / H} \circ \pr = T \pr \circ X^R$
		according
		to~\eqref{equation:lemma_reductive_homogeneous_space_fundamental_vfs_horizontal_lift_right_invariant}.
		Using~\eqref{equation:lemma_reductive_homogeneous_space_fundamental_vfs_right_left_invariant_relation}
		and
		$\pr_{\liealg{m}}(Y) 
		=
		A^i\big(\pr_{\liealg{m}}(Y)\big) A_i = A^i(Y) A_i$
		for all $Y \in \liealg{g}$, we have for $g \in G$
		\begin{equation*}
			\begin{split}
				\overline{X_{G / H}}(g) 
				&=
				\big( \id_{T G} - \mathcal{P}\big) \circ X^R(g) \\
				&=
				(\id_{T G} - \mathcal{P}) \big(\Ad_{g^{-1}}(X) \big)^L(g) \\
				&=
				\big( T_e \ell_g \circ \pr_ {\liealg{m}} \circ (T_e \ell_g)^{-1} \big) 
				T_e \ell_g \Ad_{g^{-1}}(X) \\
				&=
				T_e \ell_g \big( A^i(\Ad_{g^{-1}}(X)) A_i  \big) \\
				&=
				A^i \big(\Ad_{g^{-1}}(X) \big) A_i^L(g) .
			\end{split}
		\end{equation*}
		Next we show
		Claim~\ref{lemma:reductive_homogeneous_space_fundamental_vfs_lie_derivative_horizontal_lift_at_identity}.
		The curve $\gamma \colon \field{R} \ni t \mapsto \exp(t X) \in G$ fulfills $\gamma(0) = e$ and $\dot{\gamma}(0) = X$.
		Therefore we compute, again by
		$\pr_{\liealg{m}}(Y) = A^i(Y) A_i$
		for all $Y \in \liealg{g}$
		\begin{equation*}
			\begin{split}
				\big(\Lie_{\overline{X_{G / H}}} y^j \big)(e) A_j^L(e)
				&=
				\big( \tfrac{\D}{\D t} y^j(\gamma(t)) \at{t = 0}\big) A_j  \\
				&=
				\big( \tfrac{\D}{\D t} A^j
				\big(\Ad_{\exp(t X)^{-1}} (Y) \big) 
				\at{t = 0} \big) A_j \\
				&=
				\big(  
				A^j\big(\tfrac{\D}{\D t} \Ad_{\exp(- t X)} (Y) \at{t = 0} \big)  \big) A_j  \\
				&=
				A^j\big(-[X, Y] \big) A_j \\
				&=
				- [X, Y]_{\liealg{m}}
			\end{split}
		\end{equation*}
		as desired.
		
		Although a proof of Claim~\ref{item:lemma_reductive_homogeneous_space_fundamental_vfs_push_forward_by_canonical_action}
		can be found for example
		in~\cite[Prop. 23.20]{gallier.quaintance:2020},
		following this reference,
		we repeat it here for the reader's convenience.
		We compute for $g, k \in G$
		\begin{equation*}
			\begin{split}
				\big(( \tau_g)_* X_{G / H}\big)(\pr(k))
				&=
				\big( T_{\tau_{g^{-1}}(\pr(k))} \tau_g \big) X_{G / H} \big(\tau_g^{-1} (\pr(k)) \big) \\
				&=
				\big( T_{\tau_{g^{-1}}(\pr(k))} \tau_g \big)
				\tfrac{\D}{\D t} \tau_{\exp(t X)} \big( \tau_{g^{-1}}(\pr(k))\big) \at{t = 0} \\
				&=
				\tfrac{\D}{\D t} \tau_{g \exp(t X) g^{-1}}(\pr(k)) \at{t = 0} \\
				&=
				\tfrac{\D}{\D t} \tau_{\exp(t \Ad_g(X))}(\pr(k)) \at{t = 0} \\
				&= 
				\big(\Ad_g(X)\big)_{G / H}(\pr(k)),
			\end{split}
		\end{equation*}
		where~\eqref{equation:push_forward_covariant_derivative_on_G_H}
		is used in the first equality and we also exploited
		$\Conj_g \circ \exp = \exp \circ \Ad_g$.
		This yields the desired result.
	\end{proof}
\end{lemma}
It is well-known that there is a one-to-one correspondence between
$\Ad(H)$-invariant tensors on $\liealg{m}$ and invariant tensor fields
on $G / H$, see e.g.~\cite[Chap. 11, p. 312]{oneill:1983}.
In the sequel, we need the following lemma which can be regarded
as a special case of this assertion. 
In order to keep this text as self-contained as possible, we include
proof which is inspired by the proof
of~\cite[Prop. 23.23]{gallier.quaintance:2020}
and~\cite[Chap. 11, Prop. 22]{oneill:1983}.

\begin{lemma}
	\label{lemma:one-to-one-correspondence_invariant_tensor_fields-Ad_H_invariant_bilinear_map}
	Let $G / H$ be a reductive homogeneous space 
	with reductive decomposition $\liealg{g} = \liealg{h} \oplus \liealg{m}$.
	There is a one-to-one-correspondence between
	$\Ad(H)$-invariant $\ell$-linear maps
	$d \colon \liealg{m}^{\ell} \to \liealg{m}$
	with $\ell \in \field{N}$
	and tensor fields
	$D \in \Secinfty\big(T^* (G / H)^{\tensor \ell} \tensor T(G / H) \big)$
	on $G / H$ fulfilling
	\begin{equation}
		\label{equation:lemma_one-to-one-correspondence_invariant_tensor_fields-Ad_H_invariant_bilinear_map_invariance}
		(\tau_g)_* \big( D(X_1, \dots, X_{\ell}) \big) 
		=
		D\big((\tau_g)_* X_1, \ldots,  (\tau_g)_* X_{\ell} \big)
	\end{equation}
	for all $X_1, \ldots, X_{\ell} \in \Secinfty(T M)$ and $g \in G$
	by requiring
	\begin{equation}
		\label{equation:lemma_one-to-one-correspondence_invariant_tensor_fields-Ad_H_invariant_bilinear_map_correspondence}
		D\at{\pr(e)}\big(T_e \pr X_1, \ldots,  T_e \pr X_{\ell} \big) 
		=
		T_e \pr \big( d(X_1, \ldots, X_{\ell}) \big)
	\end{equation}
	for all $X_1, \ldots, X_{\ell} \in \liealg{m}$.
	\begin{proof}
		We use ideas that can be found in~\cite[Prop. 23.23]{gallier.quaintance:2020},
		see also~\cite[Chap. 11, Prop. 22]{oneill:1983}.
		In this proof, we write $o = \pr(e) = e \cdot H \in G / H$
		for the coset defined by $e \in G$.
		Let
		$D \in \Secinfty\big( T^* (G / H)^{\tensor \ell} \tensor T (G / H) \big)$
		be a tensor field
		satisfying~\eqref{equation:lemma_one-to-one-correspondence_invariant_tensor_fields-Ad_H_invariant_bilinear_map_invariance}.
		Using
		Lemma~\ref{lemma:isotropy_representation_equivalence},
		i.e. $T_{\pr(e)} \tau_h \circ T_e \pr\at{\liealg{m}}
		= T_o \tau_h \circ T_e \pr\at{\liealg{m}}
		= T_e \pr \circ \Ad_h \at{\liealg{m}}$
		for all $h \in H$,
		and Lemma~\ref{lemma:reductive_homogeneous_space_fundamental_vfs},
		Claim~\ref{item:lemma_reductive_homogeneous_space_fundamental_vfs_push_forward_by_canonical_action},
		i.e. $(\tau_h)_* X_{G / H} = \big( \Ad_h(X)\big)_{G / H}$
		for all $X \in \liealg{m}$
		as well as $T_e \pr X = X_{G / H}(o)$
		we compute for $X_1, \ldots, X_{\ell} \in \liealg{m}$ and $h \in H$
		\begin{equation}
			\label{equation:lemma_one-to-one-correspondence_invariant_tensor_fields-Ad_H_invariant_bilinear_map_computation_defines_Ad_H_invariant_map}
			\begin{split}
				&(T_e \pr \circ  \Ad_h ) d(X_1 , \ldots, X_{\ell} )\big) \\
				&=
				(T_{o} \tau_h  \circ T_e \pr)  d(X_1, \ldots, X_{\ell}) \\
				&=
				(T_{o} \tau_h) \big( D\at{o} \big(T_e \pr X_1, \ldots T_e \pr X_{\ell} \big) \big) \\
				&=
				(T_{o} \tau_h) D\at{o}\big( (X_1)_{G / H}(o),
				\ldots, (X_{\ell})_{G / H}(o) \big) \\
				&=
				(T_{\tau_{h^{-1}}(o)}  \tau_h )
				\big(D \big( (X_1)_{G / H}, \ldots,  (X_{\ell})_{G / H}\big)(\tau_{h^{-1}}(o)) \\
				&=
				(\tau_h)_* \big( D \big( (X_1)_{G / H}, \ldots,  (X_{\ell})_{G / H}) \big) \at{o} \\
				&=
				D\big( (\tau_h)_* (X_1)_{G / H}, \ldots, (\tau_h)_* 
				(X_{\ell})_{G / H} \big) \at{o}  \\
				&=
				D\at{o}\big( (\Ad_h(X_1))_{G / H}(o), \ldots, (\Ad_h(X_{\ell}))_{G / H}(o) \big) \\
				&=
				T_e \pr \big( d\big( \Ad_h(X_1), \ldots,  \Ad_h(X_{\ell}) \big) \big) ,
			\end{split}
		\end{equation}
		where we exploited $\tau_h(o) = o$
		for all $h \in H$.
		Thus the tensor field $D$ fulfilling~\eqref{equation:lemma_one-to-one-correspondence_invariant_tensor_fields-Ad_H_invariant_bilinear_map_invariance}
		yields an $\Ad(H)$-invariant $\ell$-linear map $d$
		via~\eqref{equation:lemma_one-to-one-correspondence_invariant_tensor_fields-Ad_H_invariant_bilinear_map_correspondence}
		since $T_e \pr \at{\liealg{m}} \colon \liealg{m} \to T_{\pr(e)} (G / H)$ is a linear isomorphism.
		Conversely, assume that $d \colon \liealg{m}^{\ell} \to \liealg{m}$
		is an $\Ad(H)$-invariant $\ell$-linear map.
		Then~\eqref{equation:lemma_one-to-one-correspondence_invariant_tensor_fields-Ad_H_invariant_bilinear_map_correspondence}
		defines a unique invariant tensor field on $G / H$
		fulfilling~\eqref{equation:lemma_one-to-one-correspondence_invariant_tensor_fields-Ad_H_invariant_bilinear_map_invariance}
		by setting for $v_{1, o}, \ldots v_{\ell, o} \in T_o (G / H)$
		\begin{equation}
			\label{equation:lemma_one-to-one-correspondence_invariant_tensor_fields-Ad_H_invariant_bilinear_map_defi_at_pr_e}
			D\at{o}(v_{1, o}, \ldots, v_{\ell, o}) 
			=
			T_e \pr \Big( d \Big( \big(T_e \pr\at{\liealg{m}} \big)^{-1} v_{1, o},
			\ldots, 
			\big(T_e \pr\at{\liealg{m}} \big)^{-1} v_{\ell, o} \Big) \Big)
		\end{equation}
		and defining for $g \in G$
		with $p = \pr(g) = \tau_g(o) \in G / H$
		for 
		$v_{1 p}, \ldots, v_{\ell p} \in T_{p} (G / H)$
		\begin{equation}
			\label{equation:lemma_one-to-one-correspondence_invariant_tensor_fields-Ad_H_invariant_bilinear_map_defi_at_general_point}
			\begin{split}
				D\at{p}(v_{1 p}, \ldots, v_{\ell p})
				&=
				T_ {o} \tau_g \Big(
				D\at{o} \big(T_{p} \tau_{g^{-1}} v_{1 p},
				\ldots, 
				(T_{p} \tau_{g^{-1}} v_{\ell p} ) \big) \Big) .
%				 \\
%				&=
%				(T_e \pr d) \Big( \big(T_e \pr\at{\liealg{m}} \big)^{-1} 
%				(T_{p} \tau_{g^{-1}}) v_{1, p},
%				\ldots, 
%				\big(T_e \pr\at{\liealg{m}} \big)^{-1} (T_{p} \tau_{g^{-1}}) v_{\ell, p} \Big)
			\end{split} 
		\end{equation}
		We first show that this yields a well-defined expression.
		Let $k \in G$ be another element with $\pr(g) = p = \pr(k)$, i.e.
		there exists a $h \in H$ with $g = k h$ and
		therefore $k^{-1} = h g^{-1}$ as well as $k = g h^{-1}$ is fulfilled.
		Using the definition of $D$
		in~\eqref{equation:lemma_one-to-one-correspondence_invariant_tensor_fields-Ad_H_invariant_bilinear_map_defi_at_pr_e}
		and~\eqref{equation:lemma_one-to-one-correspondence_invariant_tensor_fields-Ad_H_invariant_bilinear_map_defi_at_general_point},
		we compute
		\begin{equation*}
			\begin{split}
				D\at{p} \big(v_{1 p}, \ldots,  v_{\ell p} \big) 
				&=
				D\at{\pr(k)} \big(v_{1 p}, \ldots,  v_{\ell p} \big)  \\
				&=
				T_ {o} \tau_k \big( D\at{o} \big(T_{p} \tau_{k^{-1}} v_{1 p},
				\ldots, 
				T_{p} \tau_{k^{-1}} v_{\ell p} \big) \big) \\
				&=
				(T_ {o} \tau_g \circ T_{o} \tau_{h^{-1}})
				\big( D\at{o} \big( 
				(T_{o} \tau_h \circ T_{p} \tau_{g^{-1}} ) v_{1 p},
				\ldots,
				(T_{o} \tau_h \circ T_{p} \tau_{g^{-1}} ) v_{\ell p} \big) \big) \\
				&=
				(T_ {o} \tau_g \circ T_{o} \tau_{h^{-1}}) 
				\big( T_{o} \tau_{h}
				\big(
				D\at{o} 	 \big(T_{p} \tau_{g^{-1}} v_{1 p},
				\ldots, (T_{p} \tau_{g^{-1}} v_{\ell p} ) \big) \big) \big) \\
				&=
				D\at{\pr(g)} \big(v_{1 p}, \ldots,  v_{\ell p} \big),
			\end{split}
		\end{equation*}
		where the fourth equality follows by a calculation similar
		to~\eqref{equation:lemma_one-to-one-correspondence_invariant_tensor_fields-Ad_H_invariant_bilinear_map_computation_defines_Ad_H_invariant_map}
		exploiting the $\Ad(H)$-invariance of
		$d \colon \liealg{m}^{\ell} \to \liealg{m}$.
		It remains to proof that $D$ has the desired invariance property.
		To this end, let
		$X_1, \ldots, X_{\ell} \in \Secinfty\big(T (G / H) \big)$
		be vector fields and let $g \in G$.
		Suppressing the ``foot points'' of the tangent maps,
		we compute by the
		definition of $D$ for
		 $q = \pr(k) \in G / H$ represented by some $k \in G$
		\begin{equation*}
			\begin{split}
				&\big( (\tau_g)_* D(X_1, \ldots, X_{\ell}) \big)(q) \\
				&=
				T \tau_g \circ D(X_1, \ldots, X_{\ell}) \circ \tau_{g^{-1}}(q) \\
				&=
				T \tau_g \big( D \at{\tau_{g^{-1}}(q)}
				\big( X (\tau_{g^{-1}}(q)), \ldots,
				X_{\ell}(\tau_{g^{-1}}(q)) \big)\big)  \\
				&=
				T \tau_g 
				\big( 
				T \tau_{g^{-1} k} D\at{o} 
				\big( 
				T \tau_{(g^{-1} k)^{-1}} X_1 (\tau_{g^{-1}}(q)), 
				\ldots,
				T \tau_{(g^{-1} k)^{-1}} X_{\ell}(\tau_{g^{-1}}(q))
				\big) \\
				&=
				T \tau_k D\at{o} 
				\big( T \tau_{k^{-1}} \circ T \tau_{g} \circ X_1 \circ \tau_{g^{-1}}(q) , 
				\ldots,
				T \tau_{k^{-1}} \circ T \tau_{g} \circ X_{\ell} \circ \tau_{g^{-1}}(q) \big) \\
				&=
				T \tau_k D\at{o} \big( T \tau_{k^{-1}} \big( (\tau_g)_* X_1 \big)(q) ,
				\ldots,
				T \tau_{k^{-1}} \big( (\tau_g)_* X_{\ell} \big)(q) \big) \\
				&=
				D\at{q} \big( (\tau_g)_* X_1(q), 
				\ldots, (\tau_g)_* X_{\ell}(q)\big) \\
				&=
				\big( D \big( (\tau_g)_* X_1,
				\ldots,
				(\tau_g)_* X_{\ell} \big)\big)(q) .
			\end{split}
		\end{equation*}
		Thus
		$(\tau_g)_* D(X_1, \ldots, X_{\ell}) 
		=
		D \big( (\tau_g)_* X_1, \ldots, (\tau_g)_* X_{\ell} \big)$
		is shown for all $g \in G$ and
		$X_1, \ldots,  X_{\ell} \in \Secinfty\big(T (G / H) \big)$ as desired.
	\end{proof}
\end{lemma}

In the remainder part of this subsection, we investigate invariant
covariant derivatives and their relation to
$\Ad(H)$-invariant bilinear maps
$\liealg{m} \times \liealg{m} \to \liealg{m}$.
We first show that an invariant covariant derivative
on $G / H$ yields an
$\Ad(H)$-invariant bilinear map
by evaluating it on fundamental vector fields and considering
its value at $\pr(e) \in G / H$.
This is motivated by the discussion
in~\cite[Sec. 23.6]{gallier.quaintance:2020}.

Before we proceed, we point out that the right-hand side
of~\eqref{equation:lemma:G_invariant_covariant_derivative_Ad_H_invariant_bilinear_map_defintion_alpha}
in the next lemma is chosen such that it coincides
with the expression from
Definition~\ref{definition:invariant_covariant_derivative},
below.

\begin{lemma}
	\label{lemma:G_invariant_covariant_derivative_Ad_H_invariant_bilinear_map}
	Let $G / H$ be a reductive homogeneous space with reductive
	decomposition $\liealg{g} = \liealg{h} \oplus \liealg{m}$.
	Moreover, let
	$\nabla \colon \Secinfty\big( T (G / H)\big) \times \Secinfty\big( T (G / H)\big) \to \Secinfty\big( T (G / H)\big)$
	be an invariant covariant derivative.
	Then the following assertions are fulfilled:
	\begin{enumerate}
		\item
		\label{item:lemma_G_invariant_covariant_derivative_Ad_H_invariant_bilinear_map_definition}
		Let $X, Y \in \liealg{m}$.
		Then 
		\begin{equation}
			\label{equation:lemma:G_invariant_covariant_derivative_Ad_H_invariant_bilinear_map_defintion_alpha}
			\nabla_{X_{G / H}} Y_{G / H} \at{\pr(e)}
			=
			T_e \pr \big( - [X, Y]_{\liealg{m}} + \alpha(X, Y) \big)
		\end{equation}
		defines an $\Ad(H)$-invariant bilinear map 
		$\alpha \colon \liealg{m} \times \liealg{m} \to \liealg{m}$.
		\item
		\label{item:lemma_G_invariant_covariant_derivative_Ad_H_invariant_bilinear_map_uniqueness}
		Let $\nabla^1$ and $\nabla^2$ be both
		invariant covariant derivatives on $G / H$.
		Then $\nabla^1_{X_{G / H}} Y_{G / H} \at{ \pr(e)} = \nabla^2_{X_{G / H}} Y_{G / H} \at{\pr(e)}$
		for all $X, Y \in \liealg{m}$
		implies
		$\nabla^1 = \nabla^2$.
	\end{enumerate}
	\begin{proof}
		We first show
		Claim~\ref{item:lemma_G_invariant_covariant_derivative_Ad_H_invariant_bilinear_map_definition}.
		Obviously, the map
		$\liealg{m} \times \liealg{m} \ni (X, Y) \mapsto [X, Y]_{\liealg{m}} \in \liealg{m}$
		is
		bilinear and $\Ad(H)$-invariant. 
		Moreover, by exploiting that
		$T_e \pr\at{\liealg{m}} \colon \liealg{m} \to T_{\pr(e)} (G / H)$ is a
		linear isomorphism,
		Claim~\ref{item:lemma_G_invariant_covariant_derivative_Ad_H_invariant_bilinear_map_definition}
		is equivalent to the assertion that
		\begin{equation*}
			\beta \colon \liealg{m} \times \liealg{m} \to \liealg{m},
			\quad
			(X, Y) \mapsto \big( T_e \pr \at{\liealg{m}}\big)^{-1} 
			\Big( \nabla_{X_{G / H}} Y_{G / H} \at{\pr(e)}
			\Big)
		\end{equation*}
		is an $\Ad(H)$-invariant bilinear map.
		The map $\beta$ is bilinear since
		the covariant derivative
		$\nabla \colon \Secinfty\big(T (G / H) \big) \times \Secinfty\big(T (G / H) \big) \to \Secinfty\big(T (G / H) \big)$
		is $\field{R}$-bilinear.
		Next, let $h \in H$. 
		Using ideas of~\cite[Sec. 8]{nomizu:1954},
		we obtain by
		Lemma~\ref{lemma:reductive_homogeneous_space_fundamental_vfs},
		Claim~\ref{item:lemma_reductive_homogeneous_space_fundamental_vfs_push_forward_by_canonical_action}
		for $X, Y \in \liealg{m}$
		\begin{equation}
			\label{equation:lemma_G_invariant_covariant_derivative_Ad_H_invariant_bilinear_map_first_computation}
			\begin{split}
				\beta\big(\Ad_h(X), \Ad_h(Y) \big)
				&=
				\big( T_e \pr \at{\liealg{m}}\big)^{-1} 
				\Big( \nabla_{\big(\Ad_h(X) \big)_{G / H}} \big(\Ad_h(Y) \big)_{G / H} \at{\pr(e)}
				\Big) \\
				&=
				\big( T_e \pr \at{\liealg{m}}\big)^{-1} 
				\Big(
				\nabla_{(\tau_h)_* X_{G / H}} (\tau_h)_* Y_{G / H} \at{\pr(e)}
				\Big) .
			\end{split}
		\end{equation}
		Moreover, using
		Lemma~\ref{lemma:isotropy_representation_equivalence},
		i.e. $T_{\pr(e)} \tau_h \circ T_e \pr\at{\liealg{m}}
		= T_e \pr \circ \Ad_h \at{\liealg{m}}$,
		and exploiting that
		$T_e \pr \at{\liealg{m}} \colon \liealg{m}  \to T_{\pr(e)} (G / H)$
		is a linear isomorphism, we obtain by the invariance of 
		$\nabla$
		\begin{equation}
			\label{equation:lemma_G_invariant_covariant_derivative_Ad_H_invariant_bilinear_map_second_computation}
			\begin{split}
				T_e \pr \big(\Ad_h\big(\beta(X, Y) \big) \big)
				&=
				T_{\pr(e)} \tau_h \circ T_e \pr \beta(X, Y) \\
				&=
				T_{\pr(e)} \tau_h  \big(\nabla_{X_{G / H}} Y_{G / H} \at{\pr(e)} \big) \\
				&=
				T_{\tau_h^{-1}(\pr(e))} \tau_h 
				\big(\big( \nabla_{X_{G / H}} Y_{G / H} \big)(\tau_{h^{-1}}(\pr(e))) \big) \\
				&=
				\big((\tau_h)_* \nabla_{X_{G / H}} Y_{G / H} \big)(\pr(e)) \\
				&=
				\nabla_{(\tau_h)_* X_{G / H}} ((\tau_h)_* Y_{G / H}) \at{\pr(e)} \\
				&=
				T_e \pr \big( \beta(\Ad_h(X), \Ad_h(Y)) \big),
			\end{split}
		\end{equation}
		where the last equality holds
		by~\eqref{equation:lemma_G_invariant_covariant_derivative_Ad_H_invariant_bilinear_map_first_computation}.
		Obviously,~\eqref{equation:lemma_G_invariant_covariant_derivative_Ad_H_invariant_bilinear_map_second_computation}
		is equivalent to
		\begin{equation*}
			\beta\big(\Ad_h(X), \Ad_h(Y) \big)
			=
			\Ad_h \big( \beta(X, Y) \big)
		\end{equation*}
		for all $h \in H$ and $X, Y \in \liealg{m}$.
		Hence
		Claim~\ref{item:lemma_G_invariant_covariant_derivative_Ad_H_invariant_bilinear_map_definition}
		is proven.
		
		We now show
		Claim~\ref{item:lemma_G_invariant_covariant_derivative_Ad_H_invariant_bilinear_map_uniqueness}.
		Let $\nabla^1, \nabla^2 \colon 
		\Secinfty\big(T(G / H) \big) \times \Secinfty\big(T (G / H) \big) \to \Secinfty\big(T (G / H) \big)$
		be two invariant covariant derivatives.
		Then their difference
		\begin{equation*}
			D(X, Y) = \nabla^1_X Y  - \nabla^2_X Y ,
			\quad
			X, Y \in \Secinfty(T (G / H))
		\end{equation*}
		defines a tensor field
		$D \in \Secinfty\big( T^* (G / H)^{\tensor 2} \tensor T (G / H)\big)$ on $G / H$
		according to~\cite[Prop. 4.13]{lee:2013}.
		Moreover, this tensor field corresponds to an $\Ad(H)$-invariant 
		bilinear map 
		$d \colon \liealg{m} \times \liealg{m} \to \liealg{m}$
		via
		$D\at{\pr(e)} \big(T_e \pr X, T_e \pr Y \big)
		=
		T_e \pr \big( d(X, Y) \big)$
		for all $X, Y \in \liealg{m}$ by Lemma~\ref{lemma:one-to-one-correspondence_invariant_tensor_fields-Ad_H_invariant_bilinear_map}
		because of
		\begin{equation*}
			\begin{split}
				(\tau_g)_* \big( D(X, Y) \big)
				&=
				(\tau_g)_* \big( \nabla^1_X Y -  \nabla^2_X Y \big) \\
				&=
				\nabla^1_{(\tau_g)_* X} (\tau_g)_*Y 
				-
				\nabla^2_{(\tau_g)_* X} (\tau_g)_* Y \\
				&=
				D\big((\tau_g)_* X, (\tau_g)_* Y \big) 
			\end{split}
		\end{equation*}
		for all $X, Y \in \Secinfty\big(T (G / H) \big)$.
		By
		$\nabla^1_{X_{G / H}} Y_{G / H} \at{\pr(e)} 
		=
		\nabla^2_{X_{G / H}} Y_{G / H} \at{\pr(e)}$
		for all $X, Y \in \liealg{m}$, we obtain
		\begin{equation*}
			0 
			=
			\nabla^1_{X_{G / H}} Y_{G / H} \at{\pr(e)} - \nabla^2_{X_{G / H}} Y_{G / H} \at{\pr(e)}
			=
			D(X_{G / H}(e), Y_{G / H}(e))
			= T_e \pr \big( d(X, Y) \big).
		\end{equation*}
		Hence $d(X, Y) = 0$ is fulfilled
		for all $X, Y \in \liealg{m}$.
		This implies
		$D = 0$ as desired.
	\end{proof}
\end{lemma}

\subsubsection{Invariant Covariant Derivatives in Terms Horizontal Lifts}

Lemma~\ref{lemma:G_invariant_covariant_derivative_Ad_H_invariant_bilinear_map},
Claim~\ref{item:lemma_G_invariant_covariant_derivative_Ad_H_invariant_bilinear_map_definition}
shows that an invariant covariant derivative $\nabla$ on $G / H$ 
defines an $\Ad(H)$-invariant bilinear map
$\alpha \colon \liealg{m} \times \liealg{m} \to \liealg{m}$
by~\eqref{equation:lemma:G_invariant_covariant_derivative_Ad_H_invariant_bilinear_map_defintion_alpha}.
Moreover, it shows that $\nabla$ is uniquely determined
by an $\Ad(H)$-invariant bilinear map
$\alpha \colon \liealg{m} \times \liealg{m}  \to \liealg{m}$
by~\eqref{equation:lemma:G_invariant_covariant_derivative_Ad_H_invariant_bilinear_map_defintion_alpha}.

However, it does not show that 
such an invariant covariant derivative $\nabla$
on $G / H$ exists.
In the sequel, we obtain another proof for the existence of
an invariant covariant derivative $\nabla^{\alpha}$ on $G / H$
for a given $\Ad(H)$-invariant bilinear map
$\alpha \colon \liealg{m} \times \liealg{m} \to \liealg{m}$
by expressing $\nablaAlpha$ in terms of
horizontally lifted vector fields on $G$.
To this end, we state some lemmas as preparation.

\begin{lemma}
	\label{lemma:push_forward_left_invariant_vectorfield_function}
	Let $g \in G$. Then the following assertions are fulfilled:
	\begin{enumerate}
		\item
		Let $\overline{X} \in \Secinfty\big( \Hor(G)\big)$ and $g \in G$.
		Then the push-forward of $\overline{X}$ by $\ell_g \colon G \to G$
		is a horizontal
		vector field on $G$,
		i.e.
		$(\ell_g)_* \overline{X} 
		=
		T \ell_g \circ \overline{X} \circ \ell_{g^{-1}} 
		\in \Secinfty\big( \Hor(G) \big)$
		holds.
		\item
		\label{item:lemma:push_forward_left_invariant_vectorfield_function_pushforward_function}
		Let $f \colon G \to \field{R}$ be smooth and let
		$X \in \liealg{m}$. Moreover, let $g \in G$.
		Denoting by
		$(\ell_g)_* X^L
		=
		T \ell_g \circ X^L \circ \ell_{g^{-1}}$
		the push-forward of $X^L \in \Secinfty\big( \Hor(G) \big)$
		by $\ell_g \colon G \to G$,
		one has
		\begin{equation}
			(\ell_g)_* (f X^L) = \big((\ell_{g^{-1}})^* f\big) X^L .
		\end{equation}
	\end{enumerate}
	\begin{proof}
		The first claim is obvious.
		It remains to prove the second claim.
		To this end, we compute for $g, k \in G$
		\begin{equation*}
			\begin{split}
				\big((\ell_g)_* (f X^L) \big)(k)
				&=
				T_{g^{-1} k} \ell_g  \circ \big( f X^L \big) \circ \ell_{g^{-1}}(k) \\
				&=
				T_{g^{-1} k}  \ell_g 
				\big( f(g^{-1} k) X^L(g^{-1} k) \big) \\
				&=
				f\big(\ell_{g^{-1}}(k) \big)
				\big(T_{g^{-1} k} \ell_g  \circ  X^L \circ \ell_{g^{-1}} (k)\big) \\
				&= 
				\big((\ell_{g^{-1}})^* f \big)(k) X^L(k) ,
			\end{split}
		\end{equation*} 
		where exploited that $X^L\in \Secinfty\big(\Hor(G)\big)$
		is a left-invariant vector field.
		This yields the desired result.
	\end{proof}
\end{lemma}

\begin{lemma}
	\label{lemma:nablaHor}
	Let $G / H$ be a reductive homogeneous space
	with reductive decomposition
	$\liealg{g} = \liealg{m} \oplus \liealg{h}$
	and let $\Hor(G) \subseteq TG$ be the horizontal bundle
	from
	Proposition~\ref{proposition:principal_connection_reductive_homogeneous_space}.
	Moreover, let
	$\alpha \colon \liealg{m} \times \liealg{m} \to \liealg{m}$
	be an $\Ad(H)$-invariant billinear map.
	Let $\{ A_1, \ldots, A_N \} \subseteq \liealg{m}$ be a basis
	of $\liealg{m}$ and 
	denote by $A_1^L, \ldots, A_N^L \in \Secinfty\big(\Hor(G) \big)$
	the corresponding left-invariant frame.
	Let $\overline{X}, \overline{Y} \in \Secinfty(\Hor(G))$
	be horizontal vector fields on $G$
	and expand them in the frame $A_1^L, \ldots A_N^L$, i.e.
	$\overline{X} = x^i A_i^L$ and $\overline{Y} = y^j A_j^L$,
	with some uniquely determined smooth functions
	$x^i, y^j \colon G \to \field{R}$, where $i, j \in \{1, \ldots, N\}$.
	Using this notation and Einstein summation convention,
	as usual, we set
	\begin{equation}
		\label{equation:lemma_nablaHor_definition}
		\nablaHorAlpha_{\overline{X}} \overline{Y} 
		= 
		\big( \Lie_{\overline{X}}  y^j \big)A_j^L 
		+  x^i y^j \big( \alpha(A_i, A_j) \big)^L .
	\end{equation}
	Then~\eqref{equation:lemma_nablaHor_definition} defines a map
	$\nablaHorAlpha 
	\colon \Secinfty\big(\Hor(G)\big) \times \Secinfty\big(\Hor(G)\big) 
	\to \Secinfty\big(\Hor(G) \big)$
	fulfilling
	\begin{equation}
		\label{equation:lemma_nablaHor_definition_covariant_derivative_like_properties}
		\nablaHorAlpha_{f \overline{X}} \overline{Y} 
		=
		f \nablaHorAlpha_{\overline{X}} \overline{Y} 
		\quad
		\text{ and } 
		\quad
		\nablaHorAlpha_{\overline{X}} (f \overline{Y})
		=
		\big(\Lie_{\overline{X}} f \big) \overline{Y} 
		+ f \nablaHorAlpha_{\overline{X}} \overline{Y} 
	\end{equation}
	for all $f \in \Cinfty(G)$ and
	$\overline{X}, \overline{Y} \in \Secinfty\big(\Hor(G)\big)$.
	Moreover, $\nablaHorAlpha$ has the following properties:
	\begin{enumerate}
		\item
		\label{item:lemma_nablaHor_invariance}
		For each $g \in G$, the map
		$\nablaHorAlpha$ is
		invariant under $\ell_g \colon G \to G$
		in the sense that
		\begin{equation}
			\nablaHorAlpha_{\overline{X}} \overline{Y} 
			= 
			\big(\ell_{g^{-1}} \big)_*
			\big( 
			\nablaHorAlpha_{ (\ell_g)_* \overline{X}}  
			\big((\ell_g)_* \overline{Y} \big) \big),
			\quad \overline{X}, \overline{Y} \in \Secinfty\big(\Hor(G)\big) .
		\end{equation}
		holds.
		\item
		\label{item:lemma_nablaHor_uniqueness}
		The map
		$\nablaHorAlpha \colon 
		\Secinfty\big(\Hor(G)\big) \times  \Secinfty\big(\Hor(G)\big) \to \Secinfty\big(\Hor(G)\big)$
		fulfills
		\begin{equation}
			\label{equation:lemma_nablaHor_invariance_defines_nabla_Hor_e}
			\nablaHorAlpha_{X^L} Y^L \at{e} 
			=
			\big(\alpha(X, Y)\big)^L(e)
			=
			\alpha(X, Y)
			\quad X, Y \in \liealg{m} .
		\end{equation}
	\end{enumerate}
	\begin{proof}
		We first show that $\nablaHorAlpha$ is well-defined.
		Let $\{ B_1, \ldots, B_N \} \subseteq \liealg{m}$ be
		another basis of $\liealg{m}$.
		Then one has $A_i = a^k_i B_k$
		and therefore $A_i^L = a^k_i (B_k^L)$,
		where $\big(a^k_i\big) \in \matR{N}{N}$
		is some invertible matrix.
		Writing $\overline{X} = x^i A_i^L$ and
		$\overline{Y} = y^j A_j^L$ yields
		$\overline{X} = (x^i a_i^k) B_k^L$ as well as
		$\overline{Y}
		= (y^j a_j^{\ell}) B_{\ell}^L$.
		Using the $\field{R}$-linearity of Lie derivatives
		and the bilinearity of
		$\alpha \colon \liealg{m} \times \liealg{m} \to \liealg{m}$
		we compute
		\begin{equation*}
			\begin{split}
				\nablaHorAlpha_{\overline{X}} \overline{Y}
				&= 
				\big( \Lie_{\overline{X}} (y^j a^{\ell}_j) \big) B_{\ell}^L
				+ 
				(x^i a_i^k) (y^j a_j^{\ell}) \big( \alpha(B_k, B_{\ell}) \big)^L \\
				&=
				\big( \Lie_{\overline{X}} y^j \big) a_j^{\ell} B_{\ell}^L
				+ 
				x^i  y^j \big(  \alpha( a_i^k B_k, a_j^{\ell} B_{\ell}) \big)^L \\
				&=
				\big( \Lie_{\overline{X}}  y^j \big)A_j^L 
				+  x^i y^j \big( \alpha(A_i, A_j) \big)^L .
			\end{split}
		\end{equation*}
		Thus 
		$\nablaHorAlpha \colon 
		\Secinfty\big(\Hor(G)\big) \times \Secinfty\big(\Hor(G)\big) 
		\to \Secinfty\big(\Hor(G)\big)$
		is well-defined.
		By a straightforward computation,
		one verifies that $\nablaHorAlpha$ fulfills~\eqref{equation:lemma_nablaHor_definition_covariant_derivative_like_properties}.
		
		Next we show Claim~\ref{item:lemma_nablaHor_invariance}.
		By
		Lemma~\ref{lemma:push_forward_left_invariant_vectorfield_function},
		we obtain
		\begin{equation}
			\label{equation:lemma_nabla_Hor_push_forward_horizontal_lift}
			(\ell_g)_* \overline{X}
			=
			(\ell_g)_*\big( x^i A_i^L \big)
			=
			\big((\ell_{g^{-1}})^* x^i \big) A_i^L  
		\end{equation}
		for all $g \in G$
		and analogously $(\ell_g)_* \overline{Y}
		=
		\big((\ell_{g^{-1}})^* y^j \big) A_j^L$.
		By~\eqref{equation:lemma_nabla_Hor_push_forward_horizontal_lift}
		and using
		$\ell_g^* A_i^L = A_i^L$ due the left-invariance of $A_i^L$,
		we compute
		\begin{equation*}
			\begin{split}
				&\nablaHorAlpha_{((\ell_g)_* \overline{X})} 
				((\ell_g)_* \overline{Y})  \\
				&=
				\Big(\Lie_{((\ell_{g^{-1}})^* x^i) A_i^L} \big((\ell_{g^{-1}})^* y^j \big) \Big) A_j^L 
				+
				\big((\ell_{g^{-1}})^* x^i\big) 
				\big( (\ell_{g^{-1}})^* y^j \big) 
				\big(\alpha(A_i, A_j) \big)^L \\
				&=
				\Big(\big((\ell_{g^{-1}})^*x^i \big)
				\Lie_{(\ell_{g^{-1}})^* A_i^L} \big( (\ell_{g^{-1}})^*y^j \big) \Big) A_j^L 
				+ 
				\big( (\ell_{g^{-1}})^* (x^i y^j) \big)
				\big( \alpha(A_i, A_j)\big)^L \\
				&=
				\big((\ell_{g^{-1}})^* x^i \big)
				\Big( (\ell_{g^{-1}})^*\big( \Lie_{A_i^L} y^j\big) \Big) A_j^L 
				+  
				\big((\ell_{g^{-1}})^* (x^i y^j) \big) \big(\alpha(A_i, A_j)\big)^L \\
				&=
				\big( (\ell_{g^{-1}})^* \big( \Lie_{x^i A_i^L} y^j \big) \big) A_j^L
				+
				\big((\ell_{g^{-1}})^* (x^i y^j) \big) \big(\alpha(A_i, A_j)\big)^L \\
				&=
				(\ell_g)_* \big( \nablaHorAlpha_{\overline{X}}  \overline{Y} \big),
			\end{split}
		\end{equation*}
		where we used
		$\Lie_{(\ell_{g^{-1}})^* A_i^L}\big( (\ell_{g^{-1}})^*y^j \big) 
		=
		(\ell_{g^{-1}})^*(\Lie_{A_i^L} y^j)$,
		see e.g.~\cite[Prop. 8.16]{lee:2013}
		and the last equality follows by
		Lemma~\ref{lemma:push_forward_left_invariant_vectorfield_function},
		Claim~\ref{item:lemma:push_forward_left_invariant_vectorfield_function_pushforward_function}.
		Thus Claim~\ref{item:lemma_nablaHor_invariance} is shown.
		
		It remains to prove
		Claim~\ref{item:lemma_nablaHor_uniqueness}.
		To this end, let $X, Y \in \liealg{m}$.
		Then we can write 
		\begin{equation*}
			X^L = x^i A_i^L 
			\quad \text{ and } \quad
			Y^L = y^j A_j^L,
		\end{equation*}
		where the functions $x^i, y^j \colon G \to \field{R}$
		are clearly constant.
		By this notation, we compute
		by exploiting that $\Lie_{X^L} y^j = 0$
		since $y^j \colon G \to \field{R}$ is constant
		for all $j \in \{1, \ldots, N\}$
		\begin{equation*}
			\nablaHorAlpha_{X^L} Y^L
			=
			\big( \Lie_{X^L} y^j \big) A_j^L 
			+ x^i y^j \big(\alpha(A_i, A_j))^L
			=
			x^i y^j \big(\alpha(A_i, A_j))^L
			=
			\big(\alpha(X, Y)\big)^L .
		\end{equation*}
		For $g = e$, the equation above yields
		\begin{equation}
			\nablaHorAlpha_{X^L} Y^L \at{e}
			=
			\big(\alpha(X, Y)\big)^L(e)
			=
			\alpha(X, Y)
		\end{equation}
		as desired.
	\end{proof}
\end{lemma}

\begin{remark}
	The map $\nablaHorAlpha \colon \Secinfty\big(\Hor(G)\big) \times \Secinfty\big(\Hor(G)\big) \to \Secinfty\big(\Hor(G)\big)$
	in Lemma~\ref{lemma:nablaHor}
	has properties that are similar to those of a covariant derivative
	on $\Hor(G) \to G$ although its first argument is only
	defined on $\Secinfty\big(\Hor(G)\big) \subsetneq \Secinfty(T G)$.
\end{remark}
Horizontal lifts are compatible with push-forwards in the following sense.

\begin{lemma}
	\label{lemma:horizontal_lifts_vectorfieds_pushforwards_compatibility}
	Let $X \in \Secinfty\big(T (G / H) \big)$ and let
	$\overline{X} \in \Secinfty\big(\Hor(G) \big)$ be its horizontal lift.
	Then
	\begin{equation}
		\overline{(\tau_g)_* X} = (\ell_g)_* \overline{X}
	\end{equation}
	holds for $g \in G$, where $\overline{(\tau_g)_* X}$ denotes
	the horizontal lift of $(\tau_g)_* X \in \Secinfty\big(T (G / H) \big)$.
	\begin{proof}
		Let $g \in G$.
		We have $\pr \circ \ell_g = \tau_g \circ \pr$
		implying $T \pr \circ T \ell_g = T \tau_g \circ T \pr$.
		Using this equality as well as
		$T\pr \circ \overline{X} = X \circ \pr$
		we compute
		\begin{equation}				\label{equation:lemma_horizontal_lifts_vectorfieds_pushforwards_compatibility_computation}
			\begin{split}
				T \pr \circ \big((\ell_g)_* \overline{X} \big)
				&=
				T \pr \circ \big( T \ell_g \circ \overline{X} 
				\circ \ell_{g^{-1}} \big) \\
				&=
				T \tau_g \circ \big( T \pr \circ \overline{X} \big) 
				\circ \ell_{g^{-1}} \\
				&= 
				T \tau_g \circ \big( X \circ \pr \big)  \circ \ell_{g^{-1}} \\
				%					&= 
				%					T \tau_g \circ  X \circ ( \pr   \circ \ell_{g^{-1}} )  \\
				&= 
				T \tau_g  \circ X \circ \tau_{g^{-1}} \circ \pr \\
				&= 
				\big((\tau_g)_* X \big) \circ \pr.
			\end{split}
		\end{equation}
		Since $(\ell_g)_* \overline{X} \in \Secinfty\big(\Hor(G) \big)$ is
		horizontal and
		$T \pr \circ \big((\ell_g)_* \overline{X} \big)
		=
		\big((\tau_g)_* X \big) \circ \pr$
		holds
		by~\eqref{equation:lemma_horizontal_lifts_vectorfieds_pushforwards_compatibility_computation}, 
		we obtain 
		$\overline{(\tau_g)_* X} = (\ell_g)_* \overline{X}$
		as desired.
	\end{proof}
\end{lemma}

\begin{lemma}
	\label{lemma:horizontal_lifts_linear_combination_of_left_invariant_vf}
	Let $G / H$ be a reductive homogeneous space with reductive split $\liealg{g} = \liealg{h} \oplus \liealg{m}$.
	Moreover, let $\{ A_1, \ldots, A_M \} \subseteq  \liealg{m}$ be
	some vectors, not necessarily forming a basis of $\liealg{m}$,
	and let $x^i \colon G \to \field{R}$ for $i \in \{1, \ldots, M\}$
	be smooth.
	Define the horizontal vector field
	$\overline{X} \in \Secinfty\big(\Hor(G)\big)$ by
	\begin{equation}
		\label{equation:lemma_horizontal_lifts_linear_combination_of_left_invariant_vf}
		\overline{X}(g) = x^i(g) A_i^L(g), \quad g \in G .
	\end{equation}
	Then $\overline{X}$ is the horizontal lift of
	$X \in \Secinfty\big(T (G / H) \big)$ given by
	\begin{equation}
		\label{equation:lemma_horizontal_lifts_horizontal_lift}
		T \pr \circ \overline{X}  = X \circ \pr
	\end{equation}
	iff
	\begin{equation}
		\label{equation:lemma_horizontal_lifts_linear_combination_of_left_invariant_vf_equivalent_characterization}
		x^i(g) A^L(g) = x^i(g h) \big( \Ad_h(A_i)\big)^L(g)
		\quad
		\iff
		\quad
		x^i(g) A_i = x^i(g h) \Ad_h(A_i) 
	\end{equation}
	holds for all $g \in G$ and $h \in H$.
	\begin{proof}
		We first assume that
		$\overline{X} = x^i A_i^L$ is the horizontal lift
		of the vector field $X \in \Secinfty\big(T(G / H) \big)$.
		Then~\eqref{equation:lemma_horizontal_lifts_horizontal_lift}
		holds.
		Using $\pr(g h) = \pr(g) $ for all $g \in G$ and $h \in H$,
		we can
		rewrite~\eqref{equation:lemma_horizontal_lifts_horizontal_lift}
		equivalently as
		\begin{equation}
			\label{equation:lemma_horizontal_lifts_linear_combination_of_left_invariant_vf_pr_related_computation}
			\begin{split}
				T_g \pr \big(x^i(g) A_i^L(g) \big) 
				&=
				(T_g \pr) \overline{X}(g) \\
				&=
				X \circ \pr(g) \\
				&=
				X \circ \pr(g h) \\
				&=
				\big(T \pr \circ \overline{X} \big) (g h)  \\
				&=
				T_{ g h} \pr \big( x^i(gh) A_i^L(gh)  \big) \\
				&=
				x^i(g h) \big( ( T_{g h} \pr \circ T_e \ell_{g h} ) A_i \big)\\
				&=
				x^i(g h) \big( T_e (\pr \circ \ell_{g h}) A_i \big) \\
				&= 
				x^i(g h) \big( (T_{\pr(e)} \tau_{g h}  \circ T_e \pr)  A_i \big) \\
				&= 
				x^i (g h) \big( T_{\pr(e)} (\tau_g \circ \tau_h) \circ T_e \pr A_i \big) \\
				&= 
				x^i(g h) \big( T_{\pr(e)} \tau_{g}
				\circ (T_{\pr(e)} \tau_h \circ  T_e \pr) A_i \big)\\
%				&=
%				x^i(g h) \big( ( T_{\pr(e)}  \tau_g 
%				\circ (T_{e} \tau_h \circ \pr) )  A_i \big) \\
				&=
				x^i(g h) \big(  T_{\pr(e)} \tau_g  \circ (T_e \pr  \circ \Ad_h ) A_i \big) \\
				&=
				x^i(g h) \big( T_e (\tau_g \circ \pr) (\Ad_h(A_i)) \big) \\
				&=
				x^i(g h) \big( T_e (\pr \circ \ell_g) (\Ad_h(A_i)) \big) \\
				&=  x^i(g h) \big( (T_g \pr \circ T_e \ell_g) \Ad_h(A_i) \big) \\
				&=
				x^i(g h) \big( T_g \pr  \big(\Ad_h(A_i)\big)^L(g) \big),
			\end{split}
		\end{equation}
		where we exploited Lemma~\ref{lemma:isotropy_representation_equivalence},
		i.e.
		$T_{\pr(e)} \tau_h  \circ T_e \pr\at{\liealg{m}}
		= 
		T_e \pr \circ \Ad_h\at{\liealg{m}}$
		for all $h \in H$.
		Since $T_g \pr \colon \Hor(G)_g \to T_{\pr(g)} (G / H)$ is 
		a linear isomorphism for each $g \in G$,
		Equation~\eqref{equation:lemma_horizontal_lifts_linear_combination_of_left_invariant_vf_pr_related_computation}
		is equivalent to
		the left-hand side
		of~\eqref{equation:lemma_horizontal_lifts_linear_combination_of_left_invariant_vf_equivalent_characterization}.
		Applying the linear isomorphism $(T_e \ell_g)^{-1} \colon \Hor(G)_g \to \liealg{m}$ to both sides of this equality
		shows the equivalence to right-hand side
		of~\eqref{equation:lemma_horizontal_lifts_linear_combination_of_left_invariant_vf_equivalent_characterization}.
		
		Conversely, assuming that the functions
		$x^i  \colon G \to \field{R}$
		in
		the definition of $\overline{X} \in \Secinfty\big( \Hor(G)\big)$
		in~\eqref{equation:lemma_horizontal_lifts_linear_combination_of_left_invariant_vf}
		fulfill~\eqref{equation:lemma_horizontal_lifts_linear_combination_of_left_invariant_vf_equivalent_characterization}
		for all $i \in \{1, \ldots, M\}$,
		we define the map
		\begin{equation*}
			X \colon G / H \to T (G / H),
			\quad
			\pr(g) = g \cdot H \mapsto (T_g \pr) \circ \overline{X}(g) ,
		\end{equation*}
		where the coset $\pr(g) = g \cdot H \in G / H$ is represented by $g \in G$.
		Then the computation
		in~\eqref{equation:lemma_horizontal_lifts_linear_combination_of_left_invariant_vf_pr_related_computation}
		shows that $X \colon G / H \to T(G / H)$ is well-defined, i.e.
		we have
		for all $g \in G$ and $h \in H$
		\begin{equation}
			X(\pr(g)) 
			=
			T_g \pr \circ \overline{X}(g) 
			=
			T_{g h} \pr \circ \overline{X}(g h) 
			= 
			X(\pr(g h)) .
		\end{equation}
		Then $X \circ \pr = T \pr \circ \overline{X}$ holds by construction.
		Since $\pr \colon G \to G / H$ is a surjective submersion and
		$T \pr \circ \overline{X}  \colon G \to T (G / H)$ is smooth, 
		the map $X \colon G / H \to T (G / H)$ is smooth
		by~\cite[Thm. 4.29]{lee:2013}.
		Clearly, for $\pr(g) \in G /H$, one has $X(\pr(g)) \in T_{\pr(g)} (G / H)$.
		Hence $X \in \Secinfty\big(T (G / H)\big)$
		is a smooth vector field on $G / H$.
		Obviously, its horizontal lift
		is given by $\overline{X}$.
	\end{proof}
\end{lemma}

\begin{lemma}
	\label{lemma:horizontal_lift_lie_derivative}
	Let $\overline{X}, \overline{Y} \in \Secinfty\big( \Hor(G) \big)$
	be the horizontal lifts of $X, Y \in \Secinfty\big(T (G / H) \big)$,
	respectively, and let $\{ A_1, \ldots, A_N \} \subseteq \liealg{m}$
	be a basis of $\liealg{m}$.
	Denote by $A_1^L, \ldots, A_N^L$ the corresponding left-invariant
	frame of $\Secinfty(\Hor(G))$. 
	Moreover, expand $\overline{X} = x^i A_i^L$ and
	$\overline{Y} = y^j A_j^L$, where $x^i, y^j \colon G \to \field{R}$
	are smooth.
	Then
	\begin{equation}
		\label{equation:proposition:covariant_derivative_Lie_derivative_new}
		(\Lie_{\overline{X}} y^j )(g) A_j^L(g)
		=
		(\Lie_{\overline{X}} y^j )(g h) \big(\Ad_h(A_j) \big)^L(g) 
	\end{equation}
	holds. In particular,
	$\big(\Lie_{\overline{X}} y^j \big) A_j^L \in \Secinfty\big(\Hor(G) \big)$ 
	is the horizontal lift of the vector field
	$X \in \Secinfty\big(T (G / H)\big)$
	given by
	$X \circ \pr 
	=
	T \pr \circ 	\big( (\Lie_{\overline{X}} y^j ) A_j^L \big)$.
	\begin{proof}
		Let $\{ A^1, \ldots A^N \} \subseteq \liealg{m}^*$ be the dual basis of
		$\{ A_1, \ldots, A_N \}$, i.e. $A^i(A_j) = \delta^i_j$ for
		all $i, j \in \{1, \ldots, N\}$ with $\delta^i_j$
		denoting Kronecker deltas.
		Since $\overline{Y} = y^j A_j^L$ is the horizontal lift of $Y$,
		one has
		\begin{equation}
			\label{equation:lemma_horizontal_lift_lie_derivative_y_j}
			y^j(g) A_j = y^j(gh) \Ad_h(A_j)
		\end{equation}
		for all $g \in G$ and $h \in H$
		by Lemma~\ref{lemma:horizontal_lifts_linear_combination_of_left_invariant_vf}.
		Let $j \in \{1, \ldots, N\}$.
		Applying $A^j \in \liealg{m}^*$
		to~\eqref{equation:lemma_horizontal_lift_lie_derivative_y_j}
		yields by $A^j(A_k) = \delta_k^j$
		\begin{equation}
			\label{equation:lemma_horizontal_lift_lie_derivative_y_j_identity_2}
			y^j(g)
%			=
%			y^k(g) \delta^j_k
			=
			A^j\big(y^k(g) A_k \big) 
			=
			A^j\big( y^k(g h) \Ad_h(A_k) \big) 
			=
			y^k(g h)  A^j \big(\Ad_h(A_k)\big) 
		\end{equation}
		for all $g \in G$ and $h \in H$.
		Next we define the curves
		$c_1 \colon \field{R} 
		\ni t \mapsto g \exp\big(t x^i(g) A_i \big) \in G$
		and
		$c_2 \colon \field{R} 
		\ni t \mapsto g h \exp\big(t x^i(g h) A_i \big) \in G$.
		Then 
		\begin{equation*}
			\dot{c_1}(0)
			=
			\tfrac{\D}{\D t} \big(
			g \exp\big(t x^i(g) A_i \big) \big) \at{t = 0}
			=
			T_e \ell_g \big( x^i(g) A_i \big)
			=
			x^i(g) A_i^L(g)
			=
			\overline{X}(g)
		\end{equation*}
		holds and analogously one obtains
		\begin{equation*}
			\dot{c_2}(0)
			=
			\tfrac{ \D}{\D t} \big( g h \exp\big(t x^i(g h) A_i\big) \big) \at{t = 0}
			=
			T_e \ell_{g h} \big( x^i(g h) A_i \big)
			=
			x^i(gh) A_i^L(g h)
			=
			\overline{X}(g h) .
		\end{equation*}
		Expressing $y^j \colon G \to \field{R}$
		by~\eqref{equation:lemma_horizontal_lift_lie_derivative_y_j_identity_2}
		and using the definition of $c_1$ and $c_2$,
		we compute for $g \in G$ and $h \in H$
		by $\Conj_h  \circ \exp = \exp \circ \Ad_h$
		\begin{equation*}
			\begin{split}
				\big(\Lie_{\overline{X}} y^j \big) A_j^L(g)
				&=
				\Big( \tfrac{\D}{\D t} \big(
				y^j\big(c_1(t) \big) \at{t = 0}\big) \Big) A_j^L(g) \\
				&=
				\Big( \tfrac{\D}{\D t} \Big(
				y^k\big( c_1(t) h  \big) A^j\big(\Ad_h(A_k) \big) \Big) \at{t = 0} \Big) A_j^L(g) \\
				&=
				\Big( \tfrac{\D}{\D t}
				y^k\Big(g \exp(t x^i(g) A_i) h  \Big) \at{t = 0} \Big)
				\Big(A^j\big(\Ad_h(A_k) \big) A_j^L(g) \Big) \\
				&=
				\Big(\tfrac{\D}{\D t}
				y^k\Big(g \exp\Big(t x^i(gh) \Ad_h(A_i) \Big) h \Big) \at{t = 0} \Big) 
				\big(\Ad_h(A_k)\big)^L(g) \\
				&=
				\Big(\tfrac{\D}{\D t}
				y^k\Big(g \exp\Big(t \Ad_h\big( x^i(gh) A_i\big) \Big) h \Big) \at{t = 0} \Big) \big(\Ad_h(A_j)\big)^L(g) \\
				&=
				\Big(\tfrac{\D}{\D t}
				y^k \Big( g \Conj_h\Big( \exp\big(t x^i(g h) A_k \big)\Big) h \Big) \at{t = 0} \Big)  \big(\Ad_h(A_j)\big)^L(g) \\
				&=
				\Big( \tfrac{\D}{\D t}
				y^k\Big( g h \exp\big(t x^i(g h) A_k \big) \Big) \at{t = 0} \Big)
				\big(\Ad_h(A_j)\big)^L(g) \\
				&= 
				\Big( \tfrac{\D}{\D t} y^k\big(c_2(t) \big) \at{t = 0} \Big)
				\big(\Ad_h(A_j)\big)^L(g) \\
				&= 
				\big( \Lie_{\overline{X}} y^k \big)(g h) 
				\big(\Ad_h(A_j)\big)^L(g) 
			\end{split}
		\end{equation*}	
		showing~\eqref{equation:proposition:covariant_derivative_Lie_derivative_new}.
		Thus $\big(\Lie_{\overline{X}} y^j\big) A_j^L$ is the horizontal
		lift of the vector field $X$ on $G / H$ given by
		$X \circ \pr 
		=
		T \pr \circ \big( \Lie_{\overline{X}} y^j \big) A_j^L$
		according to
		Lemma~\ref{lemma:horizontal_lifts_linear_combination_of_left_invariant_vf}.
	\end{proof}
\end{lemma}
After this preparation, we are in the position to prove the existence of 
an invariant covariant derivative on $G / H$ associated with an
$\Ad(H)$-invariant bilinear map by expressing it in terms of 
horizontally lifted vector fields.

\begin{theorem}	
	\label{theorem:covariant_derivative}
	Let $X, Y \in \Secinfty\big(T (G / H)\big)$ and let
	$\overline{X}, \overline{Y} \in \Secinfty\big(\Hor(G)\big)$ denote their
	horizontal lifts. Let $\{ A_1, \ldots, A_N \} \subseteq \liealg{m}$ be a basis
	and let $A_1^L, \ldots A_N^L$ denote the corresponding
	left-invariant vector fields.
	Moreover, expand
	$\overline{X} = x^i A_i^L$ and $\overline{Y} = y^j A_j^L$
	with smooth functions $x^i, y^j \colon G \to \field{R}$
	for $i, j \in \{1, \ldots, N\}$.
	Let $\alpha \colon \liealg{m} \times \liealg{m} \to \liealg{m}$
	be an $\Ad(H)$-invariant bilinear map.
	Then 
	\begin{equation}
		(\nablaAlpha_X Y ) \circ \pr 
		=
		T \pr \big( (\Lie_{\overline{X}} y^j) A_j^L
		+ x^i y^j \big( \alpha(A_i, A_j) \big)^L \big)
	\end{equation}
	defines an invariant covariant derivative
	$\nablaAlpha \colon \Secinfty\big(T (G / H) \big) \times \Secinfty\big(T (G / H)\big) \to \Secinfty\big(T (G / H)\big)$
	and
	\begin{equation}
		\overline{\nablaAlpha_X Y} 
		=
		\big(\Lie_{\overline{X}} y^j \big)A_j^L 
		+  x^i y^j \big( \alpha(A_i, A_j)\big)^L 
	\end{equation}
	holds,
	where $\overline{\nablaAlpha_X Y}$
	denotes the horizontal lift of $\nablaAlpha_X Y$.
	Moreover, for all $X, Y \in \liealg{m}$
	\begin{equation}
		\label{equation:theorem_covariant_derivative_fundamental_vf_at_pr_e}
		\nablaAlpha_{X_{G / H}} Y_{G / H} \at{\pr(e)}
		= 
		T_e \pr \big( - [X, Y]_{\liealg{m}} + \alpha(X, Y) \big)
	\end{equation}
	is fulfilled. In addition, $\nablaAlpha$ is the unique invariant
	covariant derivative on $G / H$
	satisfying~\eqref{equation:theorem_covariant_derivative_fundamental_vf_at_pr_e}.
	\begin{proof}
		We define the covariant derivative $\nablaAlpha$ on $G / H$ by
		\begin{equation}
			\label{equation:proposition_covariant_derivative_definition}
			\big( \nablaAlpha_X Y \big) \circ \pr 
			=
			T \pr \circ \Big( \nablaHorAlpha_{\overline{X}} \overline{Y} \Big) ,
		\end{equation}
		where $\nablaHorAlpha$ is given by
		Lemma~\ref{lemma:nablaHor}.
		We first show that this definition yields a well-defined
		expression, i.e.
		$\big(\nablaAlpha_X Y) \circ \pr(g) 
		=
		\big(\nablaAlpha_X Y\big) \circ \pr(g h)$
		holds for $g \in G$ and $h \in H$.
		To this end, we calculate
		by exploiting
		Lemma~\ref{lemma:horizontal_lifts_linear_combination_of_left_invariant_vf}
		and
		Lemma~\ref{lemma:horizontal_lift_lie_derivative}
		as well as the $\Ad(H)$-invariance of
		$\alpha \colon  \liealg{m} \times \liealg{m} \to \liealg{m}$
		\begin{equation*}
			\begin{split}
				\nablaHorAlpha_{\overline{X}} \overline{Y} \at{g}
				&=
				\big( \Lie_{\overline{X}} y^j \big)(g) A_j^L(g)
				+
				x^i(g) y^j(g) \big(\alpha (A_i, A_j)\big)^L(g)
				\\
				&=
				\big(\Lie_{\overline{X}} y^j \big)(g) A_j^L(g) 
				+
				\big( \alpha\big(x^i (g ) A_i, y^j(g) A_j \big) \big)^L(g ) \\
				&=
				\big( \Lie_{\overline{X}} y^j \big)(gh) \big(\Ad_h(A_j)\big)^L(g)
				+
				\big(\alpha\big( x^i(g h) \Ad_h(A_i), y^j(g h) \Ad_h(A_j)\big)^L(g) \\
				&=
				\big( \Lie_{\overline{X}} y^j \big)(gh) \big(\Ad_h(A_j)\big)^L(g)
				+
				x^i(g h)  y^j(g h)
				\big(\alpha\big( \Ad_h(A_i), \Ad_h(A_j)\big)^L(g) \\
				&=
				\big( \Lie_{\overline{X}} y^j \big)(gh) \big(\Ad_h(A_j)\big)^L(g)
				+
				x^i(g h)  y^j(g h)
				\big(\Ad_h\big(\alpha(A_i, A_j)\big) \big)^L(g) .
			\end{split}
		\end{equation*}
		Hence~\eqref{equation:proposition_covariant_derivative_definition}
		yields a well-defined vector field on $G / H$ by
		Lemma~\ref{lemma:horizontal_lifts_linear_combination_of_left_invariant_vf}.
		
		Next we show that $\nablaAlpha$ yields
		a covariant derivative on $G / H$.
		Let $f \colon G / H \to \field {R}$ be smooth.
		By $\overline{f X} = \pr^*(f) \overline{X}$
		and the properties of $\nablaHorAlpha$
		from~\eqref{equation:lemma_nablaHor_definition_covariant_derivative_like_properties}
		in Lemma~\ref{lemma:nablaHor},
		we obtain
		\begin{equation*}
			\begin{split}
				\nablaAlpha_X (f Y) \circ \pr
				&=
				T \pr \circ \big( \nablaHor_{\overline{X}} (\pr^*(f) \overline{Y} )\big) \\
				&=
				T \pr \circ \big(  (\Lie_{\overline{X}} (\pr^* f)) \overline{Y}
				+
				\pr^*(f) \nablaHor_{\overline{X}} \overline{Y}  \big) \\
				&=
				T \pr \circ \big( \pr^* (\Lie_X f) \overline{Y}\big)
				+ T  \pr \circ \big( \pr^*(f) \nablaHor_{\overline{X}} \overline{Y}\big) \\
				&=
				\big( (\Lie_X f ) Y + f \nablaAlpha_X Y \big) \circ \pr
			\end{split}
		\end{equation*}
		due to $\Lie_{\overline{X}} (\pr^* f) = \pr^* (\Lie_X f)$
		by~\cite[Prop. 8.16]{lee:2013}
		since $X$ and $\overline{X}$ are $\pr$-related.
		Moreover, we have
		\begin{equation*}
			\begin{split}
				\nablaAlpha_{f X} Y \circ \pr 
				&=
				T \pr \circ \big( \nablaHorAlpha_{\pr^*(f) \overline{X}} \overline{Y} \big) \\
				&=
				T \pr \circ \big( \pr^*(f) \nablaHorAlpha_{\overline{X}} \overline{Y} \big)  \\
				&= 
				(\pr^* f) \big( T \pr \big( \nablaHorAlpha_{\overline{X}} \overline{Y} \big) \big) \\
				&=
				\big(f \nablaAlpha_X Y \big) \circ \pr 
			\end{split}
		\end{equation*}
		by Lemma~\ref{lemma:nablaHor}.
		Hence $\nablaAlpha$ is indeed a covariant derivative. 
		In addition, $\nablaAlpha$ is invariant.
		Indeed, by Lemma~\ref{lemma:nablaHor},
		Claim~\ref{item:lemma_nablaHor_invariance}
		and
		Lemma~\ref{lemma:horizontal_lifts_vectorfieds_pushforwards_compatibility},
		one has
		\begin{equation*}
			\begin{split}
				\big(\nablaAlpha_{(\tau_g)_* X} (\tau_g)_* Y \big) \circ \pr 
				&= 
				T \pr \circ \big(
				\nablaHorAlpha_{ (\ell_g)_* \overline{X}} (\ell_g)_* \overline{Y} \big) \\
				&=
				T \pr \circ \big( (\ell_g)_* \nablaHorAlpha_{\overline{X}} \overline{Y} \big) \\
				&=
				\big((\tau_g)_* \nablaAlpha_X Y \big) \circ \pr .
			\end{split}
		\end{equation*}
		Next let $X, Y \in \liealg{m}$
		and let $\{A^1, \ldots, A^N \} \subseteq \liealg{m}^*$ be
		the dual basis of $\{A_1, \ldots, A_N\}$.
		By Lemma~\ref{lemma:reductive_homogeneous_space_fundamental_vfs},
		Claim~\ref{item:lemma_reductive_homogeneous_space_fundamental_vfs_horizontal_lifts}, 
		we have $\overline{Y_{G / H}} = y^j A_j^L$
		with
		$y^j \colon G \ni g \mapsto y^j(g) 
		= A^j(\Ad_{g^{-1}}(Y)) \in \field{R}$
		for $j \in \{1, \ldots, N\}$.
		Thus we obtain by
		Lemma~\ref{lemma:reductive_homogeneous_space_fundamental_vfs},
		Claim~\ref{lemma:reductive_homogeneous_space_fundamental_vfs_lie_derivative_horizontal_lift_at_identity}
		\begin{equation}
			\label{equation:horizontal_lift_inv_cov_fundamental_vector_field}
			\begin{split}
%				&\overline{\nablaAlpha_{X_{G / H}} Y_{G / H}}(e) \\
%				&=
				&\nablaHorAlpha_{\overline{X_{G / H}}} \overline{Y_{G / H}} \at[\Big]{e} \\
				&=
				- [X, Y]_{\liealg{m}}
				+ 
				A^i(\Ad_{e}(X)) A^j(\Ad_{e}(Y)) \alpha(A_i, A_j)^L(e) \\
%				&=
%				-[X, Y]_{\liealg{m}} + A^i(X) A^j(Y) \alpha(A_i, A_j) \\
				&=
				- [X, Y]_{\liealg{m}} + \alpha(X, Y) ,
			\end{split}
		\end{equation}	
		where we used that
		$\overline{X_{G / H}}(e) =	A^i\big(\Ad_{e^{-1}}(X)\big) A_i^L(e) = X$
		is fulfilled for all $X \in \liealg{m}$ by
		Lemma~\ref{lemma:reductive_homogeneous_space_fundamental_vfs},
		Claim~\ref{item:lemma_reductive_homogeneous_space_fundamental_vfs_horizontal_lifts}.
		Equation~\eqref{equation:horizontal_lift_inv_cov_fundamental_vector_field} 
		is equivalent 
		to~\eqref{equation:theorem_covariant_derivative_fundamental_vf_at_pr_e}
		because of
		\begin{equation}
			\nablaAlpha_{X_{G / H}} Y_{G / H} \at{\pr(e)}
			=
			T_e \pr \Big( 
			\overline{\nablaAlpha_{X_{G / H}} Y_{G / H}}(e)  \Big)
			=
			T_e \pr \Big( \nablaHorAlpha_{\overline{X_{G / H}}}
			\overline{Y_{G / H}} \at[\Big]{e} \Big) .
		\end{equation}
		Moreover, $\nablaAlpha$ is uniquely determined
		by~\eqref{equation:theorem_covariant_derivative_fundamental_vf_at_pr_e} according to
		Lemma~\ref{lemma:G_invariant_covariant_derivative_Ad_H_invariant_bilinear_map},
		Claim~\ref{item:lemma_G_invariant_covariant_derivative_Ad_H_invariant_bilinear_map_uniqueness}.
		This yields the desired result.
	\end{proof}
\end{theorem}
The next definition makes sense due to
Theorem~\ref{theorem:covariant_derivative}.

\begin{definition}
	\label{definition:invariant_covariant_derivative}
	Let $\alpha \colon \liealg{m} \times \liealg{m} \to \liealg{m}$ be an 
	$\Ad(H)$-invariant bilinear map.
	Then the invariant covariant derivative
	$\nablaAlpha \colon \Secinfty\big(T (G / H) \big) \times \Secinfty\big( T(G / H) \big) \to \Secinfty \big( T (G / H) \big)$
	which is uniquely determined by
	\begin{equation}
		\label{equation:definition_invariant_covariant_derivative}
		\nablaAlpha_{X_{G / H}} Y_{G / H} \at{\pr(e)}
		= 
		T_e \pr
		\big( - [X, Y]_{\liealg{m}} + \alpha(X, Y) \big),
		\quad
		X, Y \in \liealg{m},
	\end{equation}
	is called the invariant covariant derivative associated
	with $\alpha$ or corresponding to $\alpha$.
\end{definition}

\begin{remark}
	The right hand side of~\eqref{equation:definition_invariant_covariant_derivative}
	in Definition~\ref{definition:invariant_covariant_derivative}
	is chosen such that the invariant covariant derivative
	$\nablaAlpha$ corresponds to the invariant affine connection
	from~\cite[Thm. 8.1]{nomizu:1954} associated with
	the $\Ad(H)$-invariant bilinear map $\alpha$,
	see
	Proposition~\ref{proposition:one-to-one-correspondence-invariant-affine-connection-covariant-derivatives}
	below.
\end{remark}

%\begin{remark}
%	\label{remark:one-to-one-correspondece-billinear-map-invariant-cov}
	As already mentioned above, the one-to-one correspondence
	between invariant
	affine connections on $G / H$ and $\Ad(H)$-invariant bilinear
	maps $\liealg{m} \times \liealg{m} \to \liealg{m}$ 
	is proven in~\cite[Thm. 8.1]{nomizu:1954}.
	Clearly, an invariant covariant derivative on $G / H$
	yields an invariant affine connection on $G / H$
	and vice versa by
	Remark~\ref{remark:affine_connection_covariant_derivatives}.
	In addition, Theorem~\ref{theorem:covariant_derivative}
	provides another proof for the existence and uniqueness 
	of an invariant covariant derivative on $G / H$
	corresponding to an $\Ad(H)$-invariant bilinear map
	$\alpha \colon \liealg{m} \times \liealg{m} \to \liealg{m}$
	via~\eqref{equation:definition_invariant_covariant_derivative}
	from
	Definition~\ref{definition:invariant_covariant_derivative_on_G_H_reductive}.
	The next proposition shows that $\nablaAlpha$
	associated to the $\Ad(H)$-invariant bilinear map
	$\alpha \colon \liealg{m} \times \liealg{m} \to \liealg{m}$
	corresponds indeed the invariant
	affine connection associated with $\alpha$ from~\cite[Thm. 8.1]{nomizu:1954}.
%\end{remark}

\begin{proposition}
	\label{proposition:one-to-one-correspondence-invariant-affine-connection-covariant-derivatives}
	Let $G / H$ be a reductive homogeneous space with 
	fixed reductive decomposition
	$\liealg{g} = \liealg{h} \oplus \liealg{m}$ and let
	$\alpha \colon \liealg{m} \times \liealg{m} \to \liealg{m}$ 
	be an $\Ad(H)$-invariant bilinear map.
	Moreover, let
	$t^{\alpha} \colon \Secinfty\big( T (G / H)\big)
	\to \End_{\Cinfty(G / H)}\big( \Secinfty\big( T (G / H) \big) \big)$
	denote the
	invariant affine connection corresponding to $\alpha$ 
	from~\cite[Thm. 8.1]{nomizu:1954}.
	Then $\nablaAlpha_X Y = t^{\alpha}(Y)(X)$ holds for all
	$X, Y \in \Secinfty\big(T (G / H) \big)$, i.e. 
	$t^{\alpha}$ is the affine connection corresponding to $\nablaAlpha$
	by Remark~\ref{remark:affine_connection_covariant_derivatives}.
	\begin{proof}
		Obviously, an invariant affine connection corresponds
		to an invariant covariant derivative
		and vice versa by
		Remark~\ref{remark:affine_connection_covariant_derivatives}.
		We now briefly recall some parts of the construction of
		the invariant affine connections from~\cite[Sec. 7-8]{nomizu:1954},
		where we adapt some notations.
		Let $N = \dim(\liealg{m})$ and $n = \dim(\liealg{g})$.
		Let $(V, x)$ be a chart of $G$, where $V \subseteq G$ is an open
		neighbourhood of $e \in G$ such that $V$ is diffeomorphic to
		$M \times K$, where $M$ and $K$ are the submanifolds of $V$
		defined by
		\begin{equation*}
			\begin{split}
				M
				&= 
				\big\{ g \in V \mid x^{N + 1}(g) = \cdots = x^n(g) = 0  \big\} , \\
				K
				&=
				\big\{g \in V  \mid x^1(g) = \cdots = x^N(g) = 0 \big\} ,
			\end{split}
		\end{equation*}
		where $M$ is denoted by $N$ in~\cite[Sec. 7]{nomizu:1954}.
		Moreover, it is assumed that $V$ is chosen such that
		the restriction of the canonical projection
		$\pr\at{M} \colon M \to G / H$ is a diffeomorphism onto its image
		denoted by $M^* = \pr(M)$.
		It is pointed out in~\cite[Sec. 7]{nomizu:1954} that the
		existence of such a chart is well-known
		referring to~\cite[Chap. IV, \S V]{chevalley:1946}.
		In addition to the assumptions
		from~\cite[Sec. 7]{nomizu:1954}, we assume that
		$T_e M = \liealg{m}$ holds.
		Clearly,
		a chart $(V, x)$ of $G$ centered at $e \in G$
		with the properties listed above can be constructed
		by exploiting that the map
		\begin{equation}
			\liealg{g} \to G,
			\quad
			X \mapsto \exp(X_{\liealg{m}}) \exp(X_{\liealg{h}})
		\end{equation}
		restricted to a suitable open neighbourhood of $0 \in \liealg{g}$
		is a diffeomorphism onto its image
		which is an open neighbourhood of $e \in G$, see
		e.g.~\cite[p. 76]{knapp:2002}.
		Obviously, $M^*$ is an open submanifold of $G / H$.
		Following~\cite[Eq. (7.1)]{nomizu:1954}, 
		we now define for $X \in \liealg{m}$ the vector field
		$X^* \in \Secinfty\big(T M^*\big)$ by 
		\begin{equation}
			\label{equation:one-to-one-correspondence-affine-connection-defi-vector-field-M-star-nomizu}
			X^*(\pr(g))
			=
			X^*(\tau_g(\pr(e))) 
			=
			\big( T_{\pr(e)} \tau_g \big) (T_e \pr X) , 
			\quad
			\pr(g) \in M^*,
			\quad
			g \in M ,
		\end{equation}
		where we exploit that $\pr\at{M} \colon M \to M^*$
		is a diffeomorphism.
		We now relate $\nablaAlpha$ to $t^{\alpha}$ which is
		uniquely determined by
		\begin{equation}
			\label{equation:one-to-one-correspondence-affine-connection-correspondence-affine-connection-nomizu}
			t^{\alpha}(Y^*) (X^*)\at{\pr(e)} 
			=
			T_e \pr \big( \alpha(X, Y) \big), 
			\quad
			X, Y \in \liealg{m}
		\end{equation}
		according to~\cite[Thm. 8.1]{nomizu:1954},
		see in particular~\cite[Eq. (8.1)]{nomizu:1954}.
		To this end, we
		rewrite~\eqref{equation:one-to-one-correspondence-affine-connection-defi-vector-field-M-star-nomizu}
		as
		\begin{equation}
			\label{equation:one-to-one-correspondence-affine-connection-horizontal-lift-restriction}
			\begin{split}
				X^*(\pr(g))
				&=
				\big( T_{\pr(e)} \tau_g \big) \circ (T_e \pr X) \\
				&=
				T_e (\tau_g \circ \pr) X \\
				&=
				T_e (\pr \circ \ell_ g) X \\
				&= 
				T_g \pr \circ T_e \ell_g X \\
				&= 
				T_g \pr \circ X^L(g)
			\end{split}
		\end{equation}
		for all $g \in M$,
		where we used $\tau_g \circ \pr = \pr \circ \ell_g$
		and $\tau_g(\pr(e)) = \pr(g)$.
		Thus the horizontal lift
		$\overline{X^*} \in \Secinfty\big(T \pr^{-1}(M^*)\big)$
		of $X^*$ restricted to
		$M \subseteq \pr^{-1}(M^*) \subseteq G$
		fulfills $\overline{X^*} \at{M} = X^L \at{M}$
		due
		to~\eqref{equation:one-to-one-correspondence-affine-connection-horizontal-lift-restriction}
		since $X^L$ is horizontal.
		Next let $\{A_1, \ldots, A_N\} \subseteq \liealg{m}$ be a basis
		of $\liealg{m}$ and
		expand $\overline{X^*} = x^i A_i^L$ with uniquely
		determined smooth functions $x^i \colon \pr^{-1}(M^*) \to \field{R}$.
		Analogously, one defines for $Y \in \liealg{m}$ the vector field
		$Y^*$ on $M^*$ whose horizontal lift
		$\overline{Y^*} \in \Secinfty\big(T \pr^{-1}(M^*) \big)$
		is expanded as
		$\overline{Y^*} = y^j A_j^L$.
		Clearly, the unique smooth functions
		$y^j \colon \pr^{-1}(M^*) \to \field{R}$
		restricted to $M$,
		i.e. $y^j \at{M} \colon M \to \field{R}$, are constant
		for all $j \in \{1, \ldots, N\}$.
		We now compute $\nablaAlpha_{X^*} Y^*\at{\pr(e)}$ 
		which makes sense since $Y^*$ is defined on the open neighbourhood
		$M^*$ of $\pr(e) \in G / H$.
		Moreover, by the assumption $T_e M = \liealg{m}$, there exists
		a smooth curve 
		$c \colon (- \epsilon, \epsilon ) \to M$
		for some $\epsilon > 0$ 
		with $c(0) = e$ and $\dot{c}(0) = X \in \liealg{m}$.
		Since the functions
		$y^j \at{M} \colon M \to \field{R}$ are constant
		for all $j \in \{1, \ldots, N\}$,
		Theorem~\ref{theorem:covariant_derivative} yields
		for $X, Y \in \liealg{m}$
		\begin{equation}
			\label{equation:one-to-one-correspondence-affine-connection-covariant-derivative-at-special-vf}
			\begin{split}
				\nablaAlpha_{X^*} Y^*\at{\pr(e)}
				&=
				T_e \pr \Big(
				\big( \Lie_{\overline{X^*}} y^j \big) A_j^L(e) + x^i(e) y^j(e) \alpha(A_i, A_j)^L(e) \Big) \\
				&=
				T_e \pr \Big( \big( \tfrac{\D}{\D t} y^j(c(t)) \at{t = 0} \big) A_i + \alpha(X, Y) \Big) \\
				&=
				T_e \pr \big( \alpha(X, Y) \big) \\
				&=
				t^{\alpha}(Y^*)(X^*) \at{\pr(e)},
			\end{split}
		\end{equation}
		where we
		used~\eqref{equation:one-to-one-correspondence-affine-connection-correspondence-affine-connection-nomizu}
		in the last equality.
		Moreover, $\nablaAlpha$ is the unique invariant covariant derivative
		on $G / H$
		satisfying~\eqref{equation:one-to-one-correspondence-affine-connection-covariant-derivative-at-special-vf}.
		Indeed, let $\nabla^{\beta}$ be the invariant covariant
		derivative associated with the
		$\Ad(H)$-invariant bilinear map
		$\beta \colon \liealg{m} \times \liealg{m} \to \liealg{m}$
		fulfilling
		$\nabla^{\beta}_{X^*} Y^* \at{\pr(e)}
		=
		t^{\alpha}(Y^*)(X^*)\at{\pr(e)}$
		for all $X, Y \in \liealg{m}$.
		Then
		\begin{equation*}
			\nabla^{\beta}_{X^*} Y^* \at{\pr(e)} 
			=
			T_e \pr \big( \beta(X, Y) \big)
			=
			t^{\alpha}(Y^*)(X^*) \at{\pr(e)}
			=
			T_e \pr \big( \alpha(X, Y) \big)
			= \nablaAlpha_{X^*} Y^*\at{\pr(e)}
		\end{equation*}
		yields $\beta = \alpha$
		implying $\nablaAlpha = \nabla^{\beta}$.
		In addition, $t^{\alpha}$ is uniquely determined
		by~\eqref{equation:one-to-one-correspondence-affine-connection-correspondence-affine-connection-nomizu}.
		Hence $\nablaAlpha$ and $t^{\alpha}$
		are both uniquely determined
		by~\eqref{equation:one-to-one-correspondence-affine-connection-covariant-derivative-at-special-vf}.
		Thus~\eqref{equation:one-to-one-correspondence-affine-connection-covariant-derivative-at-special-vf}
		implies 
		$\nablaAlpha_X Y = t^{\alpha}(Y)(X)$ for all $X, Y \in \Secinfty\big( T (G / H)\big)$ as desired.
	\end{proof}
\end{proposition}

\subsection{Torsion and Curvature}
\label{subsec:torsion_curvature}

Next we consider the torsion of an invariant covariant derivative.
This is the next lemma whose result coincides
with~\cite[Eq. (9.2)]{nomizu:1954}.

\begin{lemma}
	\label{lemma:invariant_covariant_derivative_torsion}
	Let $\nabla^{\alpha}$ be the invariant covariant derivative
	on $G / H$
	associated to the $\Ad(H)$-invariant
	bilinear map $\alpha \colon \liealg{m} \times \liealg{m} \to \liealg{m}$.
	The torsion of $\nablaAlpha$
	is the $G$-invariant tensor field $\Tor^{\alpha} \in \Secinfty\big( \Anti^2 (T^* (G / H)) \tensor T (G / H) \big)$
	defined by
	\begin{equation}
		\Tor^{\alpha}\big(X_{G / H}, Y_{G / H} \big)\at{\pr(e)} 
		=
		T_e \pr \big( \alpha(X, Y) - \alpha(Y, X) - [X, Y]_{\liealg{m}} \big) 
	\end{equation}
	for all $X, Y \in \liealg{m}$
	\begin{proof}
		We first note that
		$(\tau_g)_* [X_{G / H}, Y_{G / H}] 
		=
		[ (\tau_g)_* X_{G / H}, (\tau_g)_* Y_{G / H}]$
		holds all for $g \in G$,
		see e.g.~\cite[Cor. 8.31]{lee:2013}.
		This identity and the
		invariance of $\nablaAlpha$
		yields  that $\Tor^{\alpha}$ is $G$-invariant.
		Thus $\Tor^{\alpha}$ corresponds to an $\Ad(H)$-invariant
		bilinear map
		$\liealg{m} \times \liealg{m} \to \liealg{m}$
		by
		Lemma~\ref{lemma:one-to-one-correspondence_invariant_tensor_fields-Ad_H_invariant_bilinear_map}.
		In order to determine this bilinear map,
		writing $\pr(e) = o$, we compute
		\begin{equation*}
			\begin{split}
				\Tor^{\alpha}\big(X_{G / H}, Y_{G / H} \big)\at{o} 
				&=
				\nablaAlpha_{X_{G / H}} Y_{G / H} \at{o}
				- 
				\nablaAlpha_{Y_{G / H}} X_{G / H} \at{o}
				- [X_{G / H}, Y_{G / H}] \at{o} \\
				&=
				\nablaAlpha_{X_{G / H}} Y_{G / H} \at{o}
				- 
				\nablaAlpha_{Y_{G / H}} X_{G / H} \at{o}
				+ [X, Y]_{G / H}(o) \\
				&=
				T_e \pr \big( - [X, Y]_{\liealg{m}} + \alpha(X, Y) -
				\big(\alpha(Y, X) - [Y, X]_{\liealg{m}} \big) + [X, Y]_{\liealg{m}} \big) \\
%				&=
%				T_e \pr \Big(  \alpha(X, Y) -
%				\big( [X, Y]_{\liealg{m}} + \alpha(Y, X) \big) \Big) \\
				&=
				T_e \pr \big( \alpha(X, Y) - \alpha(Y, X) - [X, Y]_{\liealg{m}} \big)
			\end{split}			
		\end{equation*}
		for all $X, Y \in \liealg{m}$,
		where we exploited that
		$\liealg{g} \ni X \mapsto X_{G / H} \in \Secinfty\big(T (G / H) \big)$
		is an anti-morphism of Lie algebras,
		see e.g.~\cite[Sec. 6.2]{michor:2008}.
	\end{proof}
\end{lemma}

Moreover, one can compute the curvature
of $\nablaAlpha$ given by
\begin{equation}
	\label{equation:curvatuve_definition}
	R^{\alpha}(X, Y) Z 
	=
	\nablaAlpha_X \nablaAlpha_Y Z
	- \nablaAlpha_Y \nablaAlpha_X Z
	- \nablaAlpha_{[X, Y]} Z,
	\quad
	X, Y, Z \in \Secinfty\big(T (G / H) \big)
\end{equation}
by using the
expression for $\nablaAlpha$ from Theorem~\ref{theorem:covariant_derivative}.
This is the next proposition which yields an alternative
derivation for the
curvature obtained in~\cite[Eq. (9.6)]{nomizu:1954}.

\begin{proposition}
	\label{proposition:curvature_invariant_covariant_derivative}
	Let $\nabla^{\alpha}$ be the invariant covariant derivative
	on $G / H$
	associated to the $\Ad(H)$-invariant
	bilinear map $\alpha \colon \liealg{m} \times \liealg{m} \to \liealg{m}$.
	The curvature of $\nablaAlpha$
	is the $G$-invariant tensor field
	$R^{\alpha} \in \Secinfty\big( (\Anti^2 (T^* (G / H))) \tensor T^* (G  / H) \tensor T (G / H) \big) \big)$
	given by
	\begin{equation}
		\begin{split}
			&R^{\alpha}(X_{G / H}, Y_{G / H}) Z_{G / H} \at{\pr(e)} \\
			&= 
			T_e \pr \Big( 	\alpha\big(X,  \alpha(Y, Z) \big)
			-
			[[X, Y]_{\liealg{h}}, Z]
			- \alpha\big( [ X, Y]_{\liealg{m}}, Z \big) 
			- \alpha\big(Y,  \alpha(X , Z) \big) \Big)
		\end{split}
	\end{equation}
	for all $X, Y, Z \in \liealg{m}$.
	\begin{proof}
		Obviously, the curvature $R^{\alpha}$ fulfills
		\begin{equation*}
			(\tau_g)_* \big(R^{\alpha}(X, Y) Z \big)
			=
			R^{\alpha} \big ((\tau_g)_* X, (\tau_g)_* Y \big) (\tau_g)_* Z 
		\end{equation*}
		for all vector fields $X, Y, Z \in \Secinfty\big( T( G / H) \big)$
		by the invariance of $\nablaAlpha$.
		Hence $R^{\alpha}$ is uniquely determined by an $\Ad(H)$-invariant 
		$3$-linear map $\liealg{m}^3 \to \liealg{m}$
		according to
		Lemma~\ref{lemma:one-to-one-correspondence_invariant_tensor_fields-Ad_H_invariant_bilinear_map}.
		We now determine this $3$-linear map.
		To this end,
		let $X, Y, Z \in \liealg{m}$ and let
		$X_{G / H}, Y_{G / H}, Z_{G / H} \in \Secinfty\big(T (G / H)\big)$
		be the associated fundamental vector fields.
		In order to compute the curvature
		$R^{\alpha}(X_{G / H}, Y_{G / H}) Z_{G / H}$ defined
		by~\eqref{equation:curvatuve_definition}
		evaluated at the point $\pr(e) = e \cdot H \in G / H$, we
		need some computations as preparation.
		Let $\{ A_1, \ldots, A_N \} \subseteq \liealg{m}$ be a basis
		of $\liealg{m}$
		and denote by 
		$\{A^1, \ldots, A^N\} \subseteq \liealg{m}^*$ its dual basis.
		By Lemma~\ref{lemma:reductive_homogeneous_space_fundamental_vfs},
		 Claim~\ref{item:lemma_reductive_homogeneous_space_fundamental_vfs_horizontal_lifts}, we have
		$\overline{X_{G / H}} = x^i A_i^L$
		and  $\overline{Y_{G / H}} = y^j A_j^L$ 
		as well as $\overline{Z_{G / H}} = z^k A_k^L$,
		where the functions 
		$x^i, y^j, z^k \colon G \to \field{R}$ are defined by
		\begin{equation*}
			x^i(g) = A^i\big(\Ad_{g^{-1}}(X) \big),
			\quad
			y^j(g) = A^j\big(\Ad_{g^{-1}}(Y) \big),
			\quad
			z^k(g) = A^k\big(\Ad_{g^{-1}}(Z) \big)
		\end{equation*}
		for all $g \in G$ and $i, j, k \in \{1, \ldots, N\}$.
		Using this notation, we obtain by Theorem~\ref{theorem:covariant_derivative}
		\begin{equation*}
			\overline{\nablaAlpha_{X_{G / H}} Z_{G / H}} 
			=
			\big( \Lie_{\overline{Y_{G / H}}} z^k  \big) A_k^L 
			+ 
			y^j z^k \big( \alpha(A_j, A_k)\big)^L
			= 
			a^{\ell} A_{\ell}^L ,
		\end{equation*}
		where the functions $a^{\ell} \colon G \to \field{R}$
		for $\ell \in \{1, \ldots, N\}$ are given by
		\begin{equation}
			\label{equation:propositio_curvature_invariant_covariant_derivative_formula_a_ell}
			a^{\ell}
			=
			A^{\ell} \big( 
			\big( \Lie_{\overline{Y_{G / H}}} z^k  \big) A_k 
			+
			y^j z^k \big( \alpha(A_j, A_k)\big)  \big)
			=
			\Lie_{\overline{Y_{G / H}}} z^{\ell} + 	y^j z^k A^{\ell} \big( \alpha (A_j, A_k) \big) .
		\end{equation}
		In particular, evaluating $a^{\ell} \colon G \to \field{R}$
		at $g = e$ yields by
		 Lemma~\ref{lemma:reductive_homogeneous_space_fundamental_vfs},
		 Claim~\ref{lemma:reductive_homogeneous_space_fundamental_vfs_lie_derivative_horizontal_lift_at_identity}
		\begin{equation}
			\label{equation:propositio_curvature_invariant_covariant_derivative_formula_a_ell_at_e}
			\begin{split}
				a^{\ell}(e) A_{\ell}
				&=
				\big( \Lie_{Y_{G / H}} z^{\ell} \big)(e) A_{\ell}^L(e)
				+
				y^j(e) z^k(e) A^{\ell} \big( \alpha (A_j, A_k) \big) A_{\ell} \\
				&=
				- [Y, Z]_{\liealg{m}} + \alpha(Y, Z) .
			\end{split}
		\end{equation}
		Moreover, we obtain by
		Theorem~\ref{theorem:covariant_derivative}
		and~\eqref{equation:propositio_curvature_invariant_covariant_derivative_formula_a_ell_at_e}
		\begin{equation}
			\label{equation:propositio_curvature_invariant_covariant_derivative_nabla_X_Y_Z_at_e}
			\begin{split}
				\overline{\nablaAlpha_{X_{G / H}} \nablaAlpha_{Y_{G / H}} Z_{G / H}} \at{e}
				&=
				\big( \Lie_{\overline{X_{G / H}}} a^{\ell}  \big)(e) A_{\ell}
				+
				x^i(e) a^{\ell}(e) \alpha(A_i, A_{\ell}) \\
%				&=
%				\big( \Lie_{\overline{X_{G / H}}} a^{\ell}  \big)(e) A_{\ell}
%				+
%				\alpha\big(X, - [Y, Z]_{\liealg{m}} + \alpha(Y, Z) \big) \\
				&=
				\big( \Lie_{\overline{X_{G / H}}} a^{\ell}  \big)(e) A_{\ell}
				+
				\alpha\big(X,  \alpha(Y, Z) \big)
				- 
				\alpha\big(X, [Y, Z]_{\liealg{m}} \big) .
			\end{split}
		\end{equation}
		In order to obtain a more explicit expression
		for~\eqref{equation:propositio_curvature_invariant_covariant_derivative_nabla_X_Y_Z_at_e},
		we consider the first summand on the right-hand side.
		Recalling that $a^{\ell}$ is given
		by~\eqref{equation:propositio_curvature_invariant_covariant_derivative_formula_a_ell}
		one obtains
		by the Leibniz rule
		\begin{equation}
			\label{equation:propositio_curvature_invariant_covariant_derivative_Lie_a_ell}
			\begin{split}
				&\big( \Lie_{\overline{X_{G / H}}} a^{\ell}  \big) A_{\ell}  \\
%				&=
%				 \Lie_{\overline{X_{G / H}}} \big( \Lie_{\overline{Y_{G / H}}} z^{\ell} + 	y^j z^k A^{\ell} \big( \alpha (A_j, A_k) \big) \big) A_{\ell} \\
				 &=
				 \Lie_{\overline{X_{G / H}}} \big( \Lie_{\overline{Y_{G / H}}} z^{\ell} \big) A_{\ell}
				 +
				 \big(  (\Lie_{\overline{X_{G / H}}} y^j ) z^k
				 + y^j  (\Lie_{\overline{X_{G / H}}} z^k ) \big) 
				 A^{\ell} \big( \alpha(A_j, A_k) \big) A_{\ell} .
			\end{split}
		\end{equation}
		We now take a closer look
		at~\eqref{equation:propositio_curvature_invariant_covariant_derivative_Lie_a_ell}
		evaluated at $g = e$.
		We obtain for second summand of its right-hand side
		by Lemma~\ref{lemma:reductive_homogeneous_space_fundamental_vfs},
		Claim~\ref{lemma:reductive_homogeneous_space_fundamental_vfs_lie_derivative_horizontal_lift_at_identity}
		and
		$A^{\ell}\big(\alpha (A_j, A_k) \big) A_{\ell} 
		=
		\alpha(A_j, A_k)$
		\begin{equation}
			\label{equation:propositio_curvature_invariant_covariant_derivative_Lie_X_z_ell_1}
			\begin{split}
				&\big(  (\Lie_{\overline{X_{G / H}}} y^j ) z^k
				+ y^j  (\Lie_{\overline{X_{G / H}}} z^k ) \big)(e) \alpha(A_j, A_k)  \\
				&=
				\alpha\big( \big(  (\Lie_{\overline{X_{G / H}}} y^j )(e) A_j, Z \big)
				+
				\alpha\big( Y,  (\Lie_{\overline{X_{G / H}}} z^k )(e) A_k \big) \\
				&=
				- \alpha\big( [ X, Y]_{\liealg{m}}, Z \big) 
				- \alpha\big(Y, [X, Z]_{\liealg{m}} \big) .
			\end{split}
		\end{equation}
		Next we consider the first summand of the right-hand side
		of~\eqref{equation:propositio_curvature_invariant_covariant_derivative_Lie_a_ell}.
		As preparation, we note that for fixed $g \in G$, the curve
		$\gamma_Y \colon \field{R} \to G$ defined by 
		\begin{equation*}
			\gamma_Y(t) = g \exp\big( t  A^j\big(\Ad_{g^{-1}} (X) \big) A_j \big),
			\quad t \in \field{R}
		\end{equation*}
		fulfills $\gamma_Y(0) = g$ and 
		$\dot{\gamma}_Y(0) = T_e \ell_g A^j(\Ad_{g^{-1}} (Y)) A_j
		= \overline{Y_{G / H}}(g)$,
		where the last equality follows by Lemma~\ref{lemma:reductive_homogeneous_space_fundamental_vfs},
		Claim~\ref{item:lemma_reductive_homogeneous_space_fundamental_vfs_horizontal_lifts}.
		Thus we obtain
		\begin{equation}
			\label{equation:propositio_curvature_invariant_covariant_derivative_Lie_Z_ell}
			\begin{split}
				\Lie_{\overline{Y_{G / H}}} z^{\ell} (g)
				&=
				\tfrac{\D}{\D t} z^{\ell} \big(\gamma_Y(t) \big) \at{t = 0} \\
				&=
				\tfrac{\D}{\D t}
				A^{\ell}
				\Big(\Ad_{\big(g \exp\big(t A^j(\Ad_{g^{-1}} (Y) A_j) \big) \big)^{-1}}(Z) \Big)
				\at[\Big]{t = 0} \\
				&=
				\tfrac{\D}{\D t} A^{\ell} \Big( \Ad_{\exp\big(-A^j(\Ad_{g^{-1}} (Y) A_j) \big)} \big( \Ad_{g^{-1}} (Z)\big) \Big) \at[\Big]{t = 0} \\
				&=
				- A^{\ell} \big( \big[ A^j\big(\Ad_{g^{-1}} (Y) \big) A_j , \Ad_{g^{-1}}(Z) \big] \big) .
			\end{split}
		\end{equation}
		Since the curve
		$\gamma_X \colon \field{R} \ni t \mapsto \exp(t X) \in G$ fulfills
		$\gamma_X(0) = e$ and $\dot{\gamma}(0) = X = \overline{X_{G / H}}(e)$,
		Equation~\eqref{equation:propositio_curvature_invariant_covariant_derivative_Lie_Z_ell}
		yields
		\begin{equation}
			\label{equation:propositio_curvature_invariant_covariant_derivative_Lie_X_Lie_Y_z_ell}
			\begin{split}
				\big( \Lie_{\overline{X_{G / H}}} \big( \Lie_{\overline{Y_{G / H}}} z^{\ell} \big) \big) (e) A_{\ell}
				&=
				\tfrac{\D}{\D t} \big(
				\Lie_{\overline{Y_{G / H}}} z^{\ell} \big(\exp(t X) \big)
				\big) \at{t = 0} \\
				&=
				- \tfrac{\D}{\D t} A^{\ell} \big(  \big[ A^j(\Ad_{\exp(- t X)} (Y) ) A_j , \Ad_{\exp(- t X)}(Z) \big] \big) A_{\ell}
				\at{t = 0} \\
				&=
				[[X, Y]_{\liealg{m}}, Z]_{\liealg{m}}
				+
				[Y, [X, Z]]_{\liealg{m}} .
			\end{split}
		\end{equation}
		Plugging~\eqref{equation:propositio_curvature_invariant_covariant_derivative_Lie_X_z_ell_1}
		and~\eqref{equation:propositio_curvature_invariant_covariant_derivative_Lie_X_Lie_Y_z_ell}
		into~\eqref{equation:propositio_curvature_invariant_covariant_derivative_Lie_a_ell}
		yields by~\eqref{equation:propositio_curvature_invariant_covariant_derivative_nabla_X_Y_Z_at_e}
		\begin{equation}
			\label{equation:propositio_curvature_invariant_covariant_derivative_nabla_X_Y_Z_at_e_final}
			\begin{split}
				&\overline{\nablaAlpha_{X_{G / H}} \nablaAlpha_{Y_{G / H}} Z_{G / H}} \at{e} \\
%				&=
%				\big( \Lie_{\overline{X_{G / H}}} a^{\ell}  \big)(e) A_{\ell}
%				+
%				\alpha\big(X,  \alpha(Y, Z) \big)
%				- 
%				\alpha\big(X, [Y, Z]_{\liealg{m}} \big) \\
				&=
				[[X, Y]_{\liealg{m}}, Z]_{\liealg{m}}
				+
				[Y, [X, Z]]_{\liealg{m}}
				- \alpha\big( [ X, Y]_{\liealg{m}}, Z \big) 
				- \alpha\big(Y, [X, Z]_{\liealg{m}} \big) \\
				&\quad+
				\alpha\big(X,  \alpha(Y, Z) \big)
				- 
				\alpha\big(X, [Y, Z]_{\liealg{m}} \big) .
			\end{split}
		\end{equation}
		By exchanging $X$ with $Y$
		in~\eqref{equation:propositio_curvature_invariant_covariant_derivative_nabla_X_Y_Z_at_e_final},
		one obtains
		\begin{equation}
			\label{equation:propositio_curvature_invariant_covariant_derivative_nabla_Y_X_Z_at_e_final}
			\begin{split}
				&\overline{\nablaAlpha_{Y_{G / H}} \nablaAlpha_{X_{G / H}} Z_{G / H}} \at{e} \\
				&=
				[[Y, X]_{\liealg{m}}, Z]_{\liealg{m}}
				+
				[X, [Y, Z]]_{\liealg{m}}
				- \alpha\big( [ Y, X]_{\liealg{m}}, Z \big) 
				- \alpha\big(X, [Y, Z]_{\liealg{m}} \big) \\
				&\quad+
				\alpha\big(Y,  \alpha(X , Z) \big)
				- 
				\alpha\big(Y, [X, Z]_{\liealg{m}} \big) .
			\end{split}
		\end{equation}
		Moreover, we obtain by Theorem~\ref{theorem:covariant_derivative}
		\begin{equation}
			\label{equation:propositio_curvature_invariant_covariant_derivative_nabla_Lie_bracket_X_Y_Z_at_e_final}
			\begin{split}
				\overline{\nablaAlpha_{[X_{G / H}, Y_{G / H}]} Z_{G / H}} \at{e}
				&= 
				\overline{\nablaAlpha_{- [X, Y]_{G / H}} Z_{G / H}} \\
				&=
				- \big( \big( \Lie_{\overline{[X, Y]_{G / H}}} z^k \big)(e) A_k
				+
				A^i([X, Y]) z^k \alpha(A_i, A_k) \big) \\	
				&=
				[ [X, Y]_{\liealg{m}} , Z]_{\liealg{m}} - \alpha([X, Y]_{\liealg{m}}, Z) ,
			\end{split}
		\end{equation}
		where we exploited that
		$\liealg{g} \ni X \mapsto X_{G / H} \in \Secinfty( T (G / H))$
		is an anti-morphism of Lie algebras,
			see e.g.~\cite[Sec. 6.2]{michor:2008}.
		Combining~\eqref{equation:propositio_curvature_invariant_covariant_derivative_nabla_X_Y_Z_at_e_final}
		with~\eqref{equation:propositio_curvature_invariant_covariant_derivative_nabla_Y_X_Z_at_e_final}
		and~\eqref{equation:propositio_curvature_invariant_covariant_derivative_nabla_Lie_bracket_X_Y_Z_at_e_final}
		yields the following expression for the curvature
		\begin{equation*}
			\begin{split}
				&\overline{ R^{\alpha}(X_{G / H}, Y_{G / H}) Z_{G / H}}(e) \\
%				&=
%				\nablaAlpha_{X_{G / H}} \nablaAlpha_{Y_{G / H}} Z_{G / H}
%				-
%				\nablaAlpha_{Y_{G / H}} \nablaAlpha_{X_{G / H}} Z_{G / H}
%				-\nablaAlpha_{ [X_{G / H}, Y_{G / H}]} Z_{G / H}					\\
				&=
				\Big(
				[[X, Y]_{\liealg{m}}, Z]_{\liealg{m}}
				+
				[Y, [X, Z]]_{\liealg{m}}
				- \alpha\big( [ X, Y]_{\liealg{m}}, Z \big) 
				- \alpha\big(Y, [X, Z]_{\liealg{m}} \big) \\
				&\quad+
				\alpha\big(X,  \alpha(Y, Z) \big)
				- 
				\alpha\big(X, [Y, Z]_{\liealg{m}} \big)
				\Big) \\
				&\quad -
				\Big(
				[[Y, X]_{\liealg{m}}, Z]_{\liealg{m}}
				+
				[X, [Y, Z]]_{\liealg{m}}
				- \alpha\big( [ Y, X]_{\liealg{m}}, Z \big) 
				- \alpha\big(X, [Y, Z]_{\liealg{m}} \big) \\
				&\quad+
				\alpha\big(Y,  \alpha(X , Z) \big)
				- 
				\alpha\big(Y, [X, Z]_{\liealg{m}} \big)
				\Big) \\
				&\quad -
				\Big( 
				[ [X, Y]_{\liealg{m}} , Z]_{\liealg{m}} - \alpha([X, Y]_{\liealg{m}}, Z) \Big) \\
%				&=
%				[[X, Y]_{\liealg{m}}, Z]_{\liealg{m}}
%				+
%				[Y, [X, Z]]_{\liealg{m}}
%				- \alpha\big( [ X, Y]_{\liealg{m}}, Z \big) 
%				- \alpha\big(Y, [X, Z]_{\liealg{m}} \big) \\
%				&\quad+
%				\alpha\big(X,  \alpha(Y, Z) \big)
%				- 
%				\alpha\big(X, [Y, Z]_{\liealg{m}} \big)
%				\\
%				&\quad
%				- [[Y, X]_{\liealg{m}}, Z]_{\liealg{m}}
%				-
%				[X, [Y, Z]]_{\liealg{m}}
%				+ \alpha\big( [ Y, X]_{\liealg{m}}, Z \big) 
%				+ \alpha\big(X, [Y, Z]_{\liealg{m}} \big) \\
%				&\quad-
%				\alpha\big(Y,  \alpha(X , Z) \big)
%				+ 
%				\alpha\big(Y, [X, Z]_{\liealg{m}} \big)
%				 \\
%				&\quad -
%				[ [X, Y]_{\liealg{m}} , Z]_{\liealg{m}} 
%				+ \alpha([X, Y]_{\liealg{m}}, Z)  \\
				&=
				[Y, [X, Z]]_{\liealg{m}}
				- \alpha\big( [ X, Y]_{\liealg{m}}, Z \big) 
				+	\alpha\big(X,  \alpha(Y, Z) \big) \\
				&\quad
				- [[Y, X]_{\liealg{m}}, Z]_{\liealg{m}}
				-
				[X, [Y, Z]]_{\liealg{m}}
				- \alpha\big(Y,  \alpha(X , Z) \big) \\
				&=
				-
				[[X, Y]_{\liealg{h}}, Z]
				- \alpha\big( [ X, Y]_{\liealg{m}}, Z \big) 
				+	\alpha\big(X,  \alpha(Y, Z) \big)
				- \alpha\big(Y,  \alpha(X , Z) \big),
			\end{split}
		\end{equation*}
		where the last holds due to
		\begin{equation*}
			\begin{split}
%				&
				[Y, [X, Z]]_{\liealg{m}} - [[Y, X]_{\liealg{m}}, Z]_{\liealg{m}}
				-
				[X, [Y, Z]]_{\liealg{m}}
%				 \\
				&=
				\big( [Y, [X, Z] ] - [X, [Y, Z]] \big)_{\liealg{m}} -[[Y, X]_{\liealg{m}}, Z]_{\liealg{m}} \\
				&=
				- [ [X, Y]_{\liealg{h}}, Z]
			\end{split}
		\end{equation*}
		by the Jacobi identity and $[\liealg{h}, \liealg{m} ] \subseteq \liealg{m}$.
		This yields the desired result.
	\end{proof}
\end{proposition}

\subsection{Invariant Metric Covariant Derivatives}

In this short subsection, we assume that $G / H$ carries an invariant
pseudo-Riemannian metric
defined by an $\Ad(H)$-invariant scalar product 
$\langle \cdot, \cdot \rangle \colon \liealg{m} \times \liealg{m} \to \field{R}$.
We characterize all $\Ad(H)$-invariant
bilinear maps
$\alpha \colon \liealg{m} \times \liealg{m} \to \liealg{m}$
such that $\nablaAlpha$
is an invariant metric covariant derivative with respect to the
invariant pseudo-Riemannian metric corresponding to
$\langle \cdot, \cdot \rangle$.
To this end, we first recall that a covariant derivative $\nabla$
on a manifold $M$
is called compatible with
the pseudo-Riemannian metric
$g
\in \Secinfty\big(\Sym^2( T^* M)\big)$,
or metric for short,
if
\begin{equation}
	\Lie_Z \big( g(X, Y) \big)
	=
	g\big(\nabla_Z X, Y \big) + g\big(X, \nabla_Z Y \big),
	\quad
	X, Y, Z \in \Secinfty(T M)
\end{equation}
holds, see e.g.~\cite[Sec. 22.5]{michor:2008}.
\begin{notation}
	In this subsection,
	we denote by $g$ and $\overline{g}$
	a pseudo-Riemannian metric on $G$ and a
	fiber metric on $\Hor(G)$, respectively, while in the previous
	sections as well as in the sequel, we usually
	denote by $g$ an element in a Lie group $G$.
\end{notation}

\begin{proposition}
	\label{proposition:invariant_covariant_derivative_compatible_with_structure_skew_adjoint}
	Let $\alpha \colon \liealg{m} \times \liealg{m} \to \liealg{m}$ be
	an $\Ad(H)$-invariant bilinear map defining the invariant
	covariant derivative $\nablaAlpha$ on $G / H$.
	Then $\nablaAlpha$ is metric with respect to the
	invariant pseudo-Riemannian
	metric on $G / H$ defined by
	the $\Ad(H)$-invariant scalar product 
	$\langle \cdot, \cdot \rangle \colon \liealg{m} \times \liealg{m} \to \field{R}$ iff
	for each $X \in \liealg{m}$
	the linear map
	\begin{equation}
		\alpha(X, \cdot) \colon \liealg{m} \to \liealg{m},
		\quad
		Y \mapsto \alpha(X, Y)
	\end{equation}
	is skew-adjoint
	with respect to $\langle \cdot, \cdot \rangle$, i.e. 
	\begin{equation}
		\label{equation:lemma_invariant_covariant_derivative_compatible_with_structure_skew_adjoint_skew_adjoint}
		\big\langle \alpha(X, Y) , Z \big\rangle 
		=
		- \big\langle Y, \alpha(X, Z) \big\rangle
	\end{equation} 
	holds for all $X, Y, Z \in \liealg{m}$.
	\begin{proof}
		We denote the invariant pseudo-Riemannian
		metric on $G /H$ corresponding to
		$\langle \cdot, \cdot \rangle$ by
		$g \in \Secinfty\big(\Sym^2 (T^* (G / H) ) \big)$.
		Let $X, Y, Z \in \Secinfty(T (G / H))$ be vector fields with
		horizontal lifts
		$\overline{X}, \overline{Y}, \overline{Z} \in \Secinfty(\Hor(G))$.
		We expand these vector fields in a left-invariant frame, i.e.
		\begin{equation*}
			\overline{X} = x^i A_i^L ,
			\quad
			\overline{Y} = y^j A_j^L
			\quad \text{ and } \quad
			\overline{Z} = z^k A_k^L ,
		\end{equation*}
		where $\{ A_1, \ldots, A_N \} \subseteq \liealg{m}$ is a basis
		of $\liealg{m}$ and
		$x^i, y^j , z^k \colon G \to \field{R}$
		are uniquely determined smooth functions
		for $i, j, k \in \{1, \ldots, N\}$.
		Moreover, we endow $\Hor(G) \to G$ with the fiber metric
		$\overline{g} \in \Secinfty\big(\Sym^2 \Hor(G)^* \big)$
		defined by left translating
		the scalar product $\langle \cdot, \cdot \rangle$.
		Then
		\begin{equation}
			\label{equation:proposition_invariant_covariant_derivative_compatible_with_structure_skew_adjoint_functions_defined_by_metric}
			g(X, Y) \circ \pr
			=
			\pr^* \big(g(X, Y)\big)
			=
			\overline{g}(\overline{X}, \overline{Y})
			\colon
			G \to \field{R}
		\end{equation}
		holds by the definition of
		$g \in \Secinfty\big(\Sym^2 T^* (G / H) \big)$.
		Since $Z$ and $\overline{Z}$ are $\pr$-related, we obtain
		by~\cite[Prop. 8.16]{lee:2013}
		and~\eqref{equation:proposition_invariant_covariant_derivative_compatible_with_structure_skew_adjoint_functions_defined_by_metric}
		\begin{equation}
			\label{equation:lemma_invariant_covariant_derivative_compatible_with_structure_skew_adjoint_computation_1}
			\begin{split}
				\pr^* \big( \Lie_Z \big(g(X, Y) \big) \big)
				&= 
				\Lie_{\overline{Z}} \big(\pr^* \big( g(X, Y) \big) \big) \\
				&=
				\Lie_{\overline{Z}} \big( \overline{g}(\overline{X}, \overline{Y}) \big) \\
				&=
				\Lie_{\overline{Z}} \big( \overline{g}\big(x^i A_i^L, y^j A_j^L \big) \big) \\
				&=
				\Lie_{\overline{Z}} \big( x^i y^j \langle A_i , A_j \rangle \big) \\
				&=
				\big(\Lie_{\overline{Z}} x^i \big) y^j \langle A_i, A_j \rangle
				+
				x^i \big(\Lie_{\overline{Z}} y^j \big)  \langle A_i, A_j \rangle  ,
			\end{split}
		\end{equation}
		where we exploited that
		$\overline{g}\big(A_i^L, A_j^L \big) 
		=
		\langle A_i, A_j \rangle$
		holds by the definition of
		$\overline{g} \in \Secinfty\big(  \Sym^2 \Hor(G)^*  \big)$.
		Moreover, with $\nablaHorAlpha$ from Lemma~\ref{lemma:nablaHor},
		we compute
		\begin{equation}
			\label{equation:lemma_invariant_covariant_derivative_compatible_with_structure_skew_adjoint_computation_2}
			\begin{split}
				&\overline{g}\big(\nablaHorAlpha_{\overline{Z}} \overline{X}, \overline{Y} \big)
				+
				\overline{g}\big(\overline{X}, \nablaHorAlpha_{\overline{Z}} \overline{Y} \big) \\
				&=
				\overline{g}\big( \big(\Lie_{\overline{Z}} x^i \big) A_i^L + z^k x^i (\alpha(A_k, A_i))^L, y^j A_j^L) \\
				&\quad + 
				\overline{g}\big(x^i A_i^L, \big(\Lie_{\overline{Z}} y^j \big) A_j^L + z^k y^j (\alpha(A_k, A_j))^L \big) \\
				&=
				\big(\Lie_{\overline{Z}} x^i \big) y^j \langle A_i, A_j \rangle
				+ 
				x^i \big(\Lie_{\overline{Z}} y^j \big) \langle A_i, A_j \rangle \\
				&\quad +
				z^k x^i y^j \big(\big\langle \alpha(A_k, A_i), A_j \big\rangle + \big\langle A_i, \alpha(A_k, A_j) \big\rangle \big) .
			\end{split}
		\end{equation}
		By
		comparing~\eqref{equation:lemma_invariant_covariant_derivative_compatible_with_structure_skew_adjoint_computation_1}
		with~\eqref{equation:lemma_invariant_covariant_derivative_compatible_with_structure_skew_adjoint_computation_2},
		we obtain by
		$(\nablaAlpha_X Y) \circ \pr = T \pr \circ \big(\nablaHorAlpha_{\overline{X}} \overline{Y} \big)$
		\begin{equation*}
			\begin{split}
				\pr^*\big(\Lie_{Z} \big(g(X, Y) \big) \big)
				&=
				\big( \Lie_{\overline{Z}} x^i \big) y^j \langle A_i, A_j \rangle
				+ x^i \big(\Lie_{\overline{Z}} y^j \big)  \langle A_i, A_j \rangle \\
				&=
				\overline{g}\big(\nablaHorAlpha_{\overline{Z}} \overline{X}, \overline{Y} \big) 
				+
				\overline{g} \big(\overline{X}, \nablaHorAlpha_{\overline{Z}} \overline{Y} \big) \\
				&=
				\overline{g}\big(\overline{\nablaAlpha_Z X}, \overline{Y} \big)
				+ \overline{g} \big(\overline{X}, \overline{\nablaAlpha_Z Y} \big) \\
				&=
				\pr^* \big( g \big(\nablaAlpha_Z X, Y \big) + g\big(X, \nablaAlpha_Z, Y \big) \big) ,
			\end{split}
		\end{equation*}
		where the second equality holds iff  
		\begin{equation}
			\label{equation:lemma_invariant_covariant_derivative_compatible_with_structure_skew_adjoint_condition_on_alpha_basis}
			\big\langle \alpha(A_k, A_i), A_j \big\rangle
			+
			\big\langle A_i, \alpha(A_k, A_j) \big\rangle
			= 0  
		\end{equation}
		is satisfied for all $i, j, k \in \{1, \ldots, N\}$.
		Since $\{ A_1, \ldots, A_N \} \subseteq \liealg{m}$ is a
		basis of $\liealg{m}$,
		Equation~\eqref{equation:lemma_invariant_covariant_derivative_compatible_with_structure_skew_adjoint_condition_on_alpha_basis}
		is equivalent
		to~\eqref{equation:lemma_invariant_covariant_derivative_compatible_with_structure_skew_adjoint_skew_adjoint}.
		This yields the desired result
		since the pull-back by the surjective map $\pr \colon G \to G / H$ yields clearly an injective map $\pr^* \colon \Cinfty(G / H) \to \Cinfty(G)$.
	\end{proof}
\end{proposition}

We now recall an expression for the Levi-Civita covariant derivative on 
a reductive homogeneous space $G / H$ equipped with an invariant
pseudo-Riemannian metric
corresponding to the $\Ad(H)$-invariant scalar product
$\langle \cdot, \cdot \rangle \colon \liealg{m} \times \liealg{m} \to \field{R}$.
This is the next proposition which is taken from~\cite[Sec. 23.6]{gallier.quaintance:2020},
where it is stated for the Riemannian case.
However, since its proof only relies on the non-degeneracy of the
invariant pseudo-Riemannian metric and its associated $\Ad(H)$-invariant 
scalar product, it can be generalized to the pseudo-Riemannian setting.

\begin{proposition}
	\label{proposition:Levi_Civita_covariant_derivative}
	Let $G / H$ be a reductive homogeneous space with reductive decomposition $\liealg{g} = \liealg{h} \oplus \liealg{m}$.
	Moreover, let $\langle \cdot, \cdot \rangle \colon \liealg{m} \times \liealg{m} \to \field{R}$ be an $\Ad(H)$-invariant scalar product
	corresponding to an invariant pseudo-Riemannian metric on $G / H$.
	Then the Levi-Civita covariant derivative defined by this metric 
	fulfills for all $X, Y \in \liealg{m}$
	\begin{equation}
		\nablaLC_{X_{G / H}} Y_{G / H} \at{\pr(e)}
		=
		T_e \pr \big( - \tfrac{1}{2} [X, Y]_{\liealg{m}} + U(X, Y) \big),
	\end{equation}
	where $U \colon \liealg{m} \times \liealg{m} \to \liealg{m}$ is uniquely determined by
	\begin{equation}
		\label{equation:proposition_Levi_Civita_covariant_derivative_definition_U}
		2 \big\langle U(X, Y), Z \big\rangle 
		=
		\big\langle [Z, X]_{\liealg{m}}, Y \big\rangle
		+ 
		\big\langle X, [Z, Y]_{\liealg{m}} \big\rangle
	\end{equation}
	for all $Z \in \liealg{m}$.
\end{proposition}
Proposition~\ref{proposition:Levi_Civita_covariant_derivative}
can be simplified for
naturally reductive homogeneous spaces.
This is the next corollary which can be seen as a reformulation
of~\cite[Prop. 23.25]{gallier.quaintance:2020}
adapted to the pseudo-Riemannian setting.

\begin{corollary}
	\label{corollary:Levi_Civita_covariant_derivative_naturally_reductive}
	Let $G / H$ be a naturally reductive homogeneous space.
	Then
	\begin{equation}
		\nablaLC_{X_{G / H}} Y_{G / H} \at{\pr(e)}
		= T_e \pr \big( - \tfrac{1}{2} [X, Y]_{\liealg{m}} \big)
	\end{equation}
	holds for all $X, Y \in \liealg{m}$.
	\begin{proof}
		This is proven in~\cite[Prop. 23.25]{gallier.quaintance:2020}.
		Nevertheless, we include the proof here, as well.
		Since $G / H$ is naturally reductive, we have
		$\langle [X, Y]_{\liealg{m}}, Z \rangle = \langle X, [Y, Z]_{\liealg{m}} \rangle$ for all $X, Y , Z \in \liealg{m}$
		implying $\langle [Y, X]_{\liealg{m}}, Z \rangle
		+ \langle X, [Y, Z]_{\liealg{m}} \rangle = 0$.
		Using the definition of $U \colon \liealg{m} \times \liealg{m}  \to \liealg{m}$
		in \eqref{equation:proposition_Levi_Civita_covariant_derivative_definition_U}
		of Proposition~\ref{proposition:Levi_Civita_covariant_derivative}
		yields for $X, Y, Z \in \liealg{m}$
		\begin{equation*}
			\langle U(X, Z), Y \rangle
			= \langle [Y, X]_{\liealg{m}} , Z \rangle
			+
			\langle X, [Y, Z]_{\liealg{m}} \rangle 
			= 0 
		\end{equation*}
		implying $U(X, Y) = 0$ for all $X, Y \in \liealg{m}$.
		Thus
		Proposition~\ref{proposition:Levi_Civita_covariant_derivative}
		yields the desired result.
	\end{proof}
\end{corollary}
Next we relate the Levi-Civita covariant derivative on $G / H$,
equipped with an invariant pseudo-Riemannian metric,
to an invariant covariant derivative on $G / H$.
This is the next remark which coincides with~\cite[Sec. 13]{nomizu:1954}.

\begin{remark}
	\label{remark:invariant_covariant_derivatives_levi_civita}
	Let $G / H$ be a reductive homogeneous space equipped with an invariant pseudo-Riemannian metric corresponding the $\Ad(H)$-invariant scalar product
	$\langle \cdot, \cdot \rangle \colon \liealg{m} \times \liealg{m} \to \field{R}$.
	Then the action $\tau \colon G \times G / H \to G / H$
	is isometric, see e.g.~\cite[Chap. 11, Prop. 22]{oneill:1983}.
	Therefore $\nablaLC$ is an invariant covariant derivative on $G / H$
	by~\cite[Prop. 5.13]{lee:2018}.
	Thus
	Lemma~\ref{lemma:G_invariant_covariant_derivative_Ad_H_invariant_bilinear_map},
	Claim~\ref{item:lemma_G_invariant_covariant_derivative_Ad_H_invariant_bilinear_map_uniqueness}
	implies by
	Proposition~\ref{proposition:Levi_Civita_covariant_derivative}
	that 
	$\nablaLC = \nablaAlpha$ holds,
	where $\alpha \colon \liealg{m} \times \liealg{m} \to \liealg{m}$
	is defined by
	\begin{equation}
		\alpha(X, Y) = \tfrac{1}{2} [X, Y]_{\liealg{m}} + U(X, Y)
	\end{equation}
	for all $X, Y \in \liealg{m}$ in accordance
	with~\cite[Thm. 13.1]{nomizu:1954}.
	In particular,
	one has $\nablaLC = \nablaAlpha$ for
	$\alpha(X, Y) = \tfrac{1}{2} [X, Y]_{\liealg{m}}$ 
	by
	Corollary~\ref{corollary:Levi_Civita_covariant_derivative_naturally_reductive} 
	if $G / H$ is
	naturally reductive.
	This coincides with~\cite[Eq. (13.1)]{nomizu:1954}.
\end{remark}

\subsection{Parallel Vector Fields along Curves}
\label{subsec:parallel_vector_fields}

Having an expression for $\nablaAlpha$ on a reductive homogeneous space
$G / H$ in terms of horizontally lifted vector fields on $G$
allows for determining the associated
covariant derivative of vector fields along a given curve
on $G / H$
in terms of horizontal lifts, as well.
In this subsection, an ODE for a specific curve in $\liealg{m}$
is determined which is fulfilled iff
the corresponding vector field along the given curve is parallel.
Let
\begin{equation}
	\gamma \colon I \to G / H
\end{equation}
be a curve and let
\begin{equation}
	\widehat{Z} \colon I \to T (G / H)
\end{equation} 
be a vector field along ${\gamma}$, i.e
\begin{equation}
	\widehat{Z}(t) \in T_{{\gamma}(t)}(G / H),
	\quad 
	t \in I.
\end{equation}
Moreover, let
\begin{equation}
	g \colon I \to G
\end{equation}
denote a horizontal lift of $\gamma$ with respect to the principal
connection $\mathcal{P} \in \Secinfty\big(\End(T G) \big)$ from
Proposition~\ref{proposition:principal_connection_reductive_homogeneous_space}.
It is well-known that $g$ is unique up to the initial condition
$g(t_0) = g_0 \in G_{\gamma(0)}$.
Furthermore, the curve $g$ is defined on the whole interval $I$ 
since principal connections are complete, see 
e.g.~\cite[Thm. 19.6]{michor:2008}.
Let $\overline{Z} \colon I \to \Hor(G)$ be the horizontal lift of
$\widehat{Z}$ along $g$, i.e.
\begin{equation}
	\overline{Z}(t) 
	=
	\big( T_{g(t)} \pr\at{\Hor(G)_{g(t)}} \big)^{-1} \widehat{Z}(t),
	\quad t \in I .
\end{equation}
Next we define the curves in $\liealg{m}$ associated with
$g$ and $\overline{Z}$, namely
\begin{equation}
	x \colon I \to \liealg{m}, 
	\quad t \mapsto x(t)
	=
	\big( T_{e} \ell_{g(t)} \big)^{-1} \dot{g}(t)
\end{equation}
and 
\begin{equation}
	z \colon I \to \liealg{m}, 
	\quad t \mapsto z(t)
	=
	\big( T_{e} \ell_{g(t)} \big)^{-1} \overline{Z}(t) .
\end{equation}
We now consider the covariant derivative of $\widehat{Z}$ along ${\gamma}$.
This is next proposition
which can be seen as a generalization
of~\cite[Lem. 1]{jurdjevic.markina.leite:2023},
where we use the notation which has been introduced above.

\begin{proposition}
	\label{proposition:covariant_derivative_along_curve}
	Let $G / H$ be a reductive homogeneous space and
	let ${\gamma} \colon I \to G /H$ be smooth. Let
	$g \colon I \to G$ be a horizontal lift of $\gamma$.
	Moreover, let $\widehat{Z} \colon I \to T (G / H)$ be a vector
	field along $\gamma$ with horizontal lift
	$\overline{Z} \colon I \to \Hor(G)$ along $g \colon I \to G$.
	Let $\{ A_1, \ldots, A_N \} \subseteq \liealg{m}$ be a basis and write
	\begin{equation}
		\dot{g}(t) = x^i(t) A_i^L (g(t))
		\quad \text{ and } \quad
		\overline{Z}(t) = z^j(t) A_j^L(g(t)) 
	\end{equation}
	for some uniquely determined smooth
	functions $x^i , z^j \colon I \to \field{R}$.
	Let $\alpha \colon \liealg{m} \times \liealg{m} \to \liealg{m}$
	be an $\Ad(H)$-invariant bilinear map
	and let $\nablaAlpha$ be the corresponding
	invariant covariant derivative on $G / H$. 
	Then the associated covariant derivative
	of $\widehat{Z}$ along ${\gamma}$
	lifted to a horizontal vector field along $g \colon I \to G$ 
	is given by
	\begin{equation}
		\begin{split}
			\overline{\nablaAlpha_{\dot{\gamma}(t)} \widehat{Z}} \at[\Big]{t}
			&=
			\big(\tfrac{\D}{\D t} z^j(t) \big) A_j^L(g(t)) 
			+ x^i(t) z^j(t) \big( \alpha(A_i, A_j) \big)^L(g(t)) \\
			&=
			\big( \dot{z}(t) \big)^L(g(t))
			+  \big( \alpha(x(t), z(t) \big)^L(g(t)) 
		\end{split}
	\end{equation}
	for all $t \in I$, where
	$z \colon I \ni t \mapsto z^i(t) A_i
	= (T_e \ell_{g(t)})^{-1} \widehat{Z}(t) \in \liealg{m}$
	and
	$x \colon I \ni t \mapsto x^i(t) A_i 
	= (T_e \ell_{g(t)})^{-1} \dot{g}(t) \in \liealg{m}$.
	\begin{proof}
		The proof is essentially given by applying
		Theorem~\ref{theorem:covariant_derivative}.
		To this end, we define the vector field
		$X \colon I \to T(G / H)$
		along $\gamma \colon I \to G / H$ by
		\begin{equation*}
			X(t) = \dot{{\gamma}}(t),
			\quad
			t \in I
		\end{equation*}
		and we denote by $\overline{X} \colon I \to T G$ the horizontal
		lift of $X$ along $g \colon I \to G$.
		Moreover, for fixed $t_0 \in I$, we extend $X$ and $\widehat{Z}$
		to vector fields defined on an open neighbourhood
		$O \subseteq G / H$ of $\gamma(t_0)$.
		These vector fields are denoted by 
		\begin{equation*}
			\widetilde{X} \in \Secinfty\big( T (G / H) \at{O} \big) 
			\quad \text{ and } \quad
			\widetilde{Z} \in \Secinfty\big( T (G / H) \at{O} \big), 
		\end{equation*}
		respectively. In particular, 
		\begin{equation*}
			\widetilde{X}(\gamma(t)) = X(t) = \dot{{\gamma}}(t)
			\quad \text{ and } \quad 
			\widehat{Z}(t) = \widetilde{Z}({\gamma}(t))
		\end{equation*}
		is fulfilled for all $t$ in a suitable open neighbourhood of $t_0$ in $I$.
		Moreover, their horizontal lifts
		$\overline{\widetilde{X}}, \overline{\widetilde{Y}} 
		\in \Secinfty\big( \Hor(G)\at{\pr^{-1}(O)} \big)$ 
		fulfill
		\begin{equation*}
			\overline{X}(t) = \overline{\widetilde{X}}(g(t)) = \dot{g}(t)
			\quad \text{ and } \quad
			\overline{Z}(t) = \overline{\widetilde{Z}}(g(t)) .
		\end{equation*}
		These horizontal lifts can be expanded in the global frame $A_1^L \ldots, A_N^L$ of $\Hor(G)$. 
		We write for $t \in I$ in a suitable open neighbourhood of $t_0$
		\begin{equation*}
			\overline{{X}}(t) = x^i(t) A_i^L(g(t)) = \big( x(t) \big)^L(g(t))
			\quad \text{ and } \quad
			\overline{{Z}}(t) = z^j(t) A_j^L(g(t)) = \big( z(t) \big)^L(g(t)) .
		\end{equation*}		
		Similarly, we expand
		\begin{equation*}
			\overline{\widetilde{X}} = \widetilde{x}^i A_i^L\at{\pr^{-1}(O)} 
			\quad \text{ and } \quad
			\overline{\widetilde{Z}}(t) = \widetilde{z}^j A_j \at{\pr^{-1}(O)} ,
		\end{equation*} 
		where
		$\widetilde{x}^i, \widetilde{z}^j \colon \pr^{-1}(O) \subseteq G \to \field{R}$
		are uniquely determined smooth functions for $i, j \in \{1, \ldots, N\}$.
		By construction
		\begin{equation*}
			x^i(t) = \widetilde{x}^i(g(t)) 
			\quad \text{ and } \quad
			z^j(t) = \widetilde{z}^j(g(t))
		\end{equation*}
		holds for all $t$ in a suitable open neighbourhood of $t_0$ in $I$.
		We now use~\cite[Thm. 4.24]{lee:2018}
		as well as 
		Theorem~\ref{theorem:covariant_derivative}
		to compute the horizontal lift of the covariant derivative of
		$\widehat{Z}$ along ${\gamma}$.
		We obtain for $t \in I$ in a suitable neighbourhood of $t_0$
		\begin{equation*}
			\begin{split}
				\overline{\nablaAlpha_{\dot{{\gamma}}(t)} \widehat{Z}}\at[\Big]{t}
				&=
				\overline{\nablaAlpha_{{\widetilde{X}}} \widetilde{Z}} \at[\Big]{g(t)} \\
				&=
				\nablaHorAlpha_{\overline{\widetilde{X}}} \overline{\widetilde{Z}} \at[\Big]{g(t)} \\
				&=
				\big( \Lie_{\overline{\widetilde{X}}}  \widetilde{z}^j \big)(g(t)) A_j^L(g(t))
				+  \widetilde{x}^i(g(t)) \widetilde{z}^j(g(t)) 
				\big( \alpha(A_i, A_j)\big)^L(g(t)) \\
				&= 
				\big( \tfrac{\D}{\D t} \widetilde{z}^j(g(t))\big) A_j^L(g(t)) 
				+ \widetilde{x}^i(g(t)) \widetilde{z}^j(g(t)) 
				\big( \alpha(A_i, A_j)\big)^L (g(t)) \\
				&=
				\Big( \big( \tfrac{\D}{\D t} z^j(t)\big) A_i\Big)^L(g(t)) 
				+  
				\Big( \alpha\big( x^i(t) A_i, z^i(t) A_i \big) \Big)^L(g(t)) \\
				&=  \big( \dot{z}(t) \big)^L(g(t)) + 
				\big( \alpha(x(t), z(t))\big)^L(g(t)) .
			\end{split} 
		\end{equation*}
		Applying this argument for each $t_0 \in I$
		yields the desired result.
	\end{proof}
\end{proposition}	
Proposition~\ref{proposition:covariant_derivative_along_curve}
allows for characterizing parallel vector fields along curves.

\begin{corollary}
	\label{corollary:parallel_vectorfields_ODE_on_m}
	Let $G / H$ be a reductive homogeneous space equipped with
	the invariant covariant derivative $\nablaAlpha$.
	Let $\widehat{Z} \colon I \to T(G / H)$ be a vector field along
	${\gamma} \colon I \to G / H$.
	Then $\widehat{Z}$ is parallel along $\gamma$ iff the ODE
	\begin{equation}
		\dot{z}(t) = - \alpha\big(x(t), z(t) \big)
	\end{equation}
	is fulfilled, where $x, z \colon I \to \liealg{m}$ are defined as in
	Proposition~\ref{proposition:covariant_derivative_along_curve}.
	\begin{proof}
		Let $g \in G$. The map $(T_e \ell_g)^{-1} \colon \Hor(G) \to \liealg{m}$ 
		is a linear isomorphism which fulfills
		$(T_e \ell_g)^{-1} \xi^L(g) = \xi$ 
		for all $\xi \in \liealg{m}$ by the definition of left-invariant
		vector fields. Hence
		Proposition~\ref{proposition:covariant_derivative_along_curve}
		yields the desired result due to 
		\begin{equation*}
			0
			=
			\overline{\nablaAlpha_{\dot{\gamma}(t)} \widehat{Z}\at{t}}
			\quad\iff\quad
			0 
			=
			\big(T_e \ell_{g(t)} \big)^{-1 } \overline{\nablaAlpha_{\dot{\gamma}(t)} \widehat{Z}}\at[\Big]{t}
			=
			\dot{z}(t) + \alpha\big(x(t), z(t) \big)
		\end{equation*}
		for $t \in I$.
	\end{proof}
\end{corollary}

\subsection{Geodesics}
\label{subsec:geodesics}

In this short section, we consider geodesics on the reductive homogeneous space $G / H$ with respect to an invariant covariant derivative $\nablaAlpha$.
Recall that a curve $\gamma \colon I \to G / H$ is a geodesic
if the vector field $\dot{\gamma} \colon I \to T (G / H)$
along $\gamma$ is a parallel.
Thus Corollary~\ref{corollary:parallel_vectorfields_ODE_on_m}
can be used to obtain the following characterization of
the geodesics on $G / H$ with respect to $\nablaAlpha$.

\begin{lemma}
	\label{lemma:geodesics_on_G_H_ODE_on_m_characterization}
	Let $G / H$ be a reductive homogeneous space endowed with
	the invariant covariant derivative $\nablaAlpha$
	corresponding to the $\Ad(H)$-invariant bilinear map
	$\alpha \colon \liealg{m} \times \liealg{m} \to \liealg{m}$.
	Let $\gamma \colon I \to G  / H$ be a curve in $G / H$ and 
	$g \colon I \to G$ be a horizontal lift of $\gamma$.
	Define $x \colon I \ni t  \mapsto x(t) = (T_e \ell_{g(t)})^{-1} \dot{g}(t) \in \liealg{m}$.
	Then $\gamma \colon I \to G  / H$ is a geodesic
	with respect to $\nablaAlpha$ iff the ODE
	\begin{equation}
		\dot{x}(t) = - \alpha\big(x(t), x(t) \big)
	\end{equation}
	is satisfied for all $t \in I$.
	\begin{proof}
		The curve $\gamma \colon I \to G / H$ is a geodesic with respect to $\nablaAlpha$ iff
		the vector field $\dot{\gamma}  \colon I \to T (G / H)$
		is parallel along $\gamma \colon I \to G / H$.
		Thus the desired result follows by Corollary~\ref{corollary:parallel_vectorfields_ODE_on_m}.
	\end{proof}
\end{lemma}

We now apply
Lemma~\ref{lemma:geodesics_on_G_H_ODE_on_m_characterization}
to a reductive homogeneous space equipped with the Levi-Civita
covariant derivative
defined by some invariant pseudo-Riemannian metric.
Inspired by
the well known characterization of geodesics on
a Lie group equipped with 
a left-invariant metric given
in~\cite[Ap. B]{arnold:1978},
see also~\cite[Sec. 4]{Elshafei.etal:2023}
for a discussion in the complex setting,
we obtain the next corollary
which generalizes the description
of geodesics on Lie groups
equipped with left-invariant metrics.

\begin{corollary}
	\label{corollary:geodesics_pseudo-riemannian_equation_on_m}
	Let $G / H$ be a reductive homogeneous space and let $\langle \cdot, \cdot \rangle \colon \liealg{m} \times \liealg{m} \to \field{R}$
	be an $\Ad(H)$-invariant scalar product.
	Moreover, let $\nablaLC$ denote the Levi-Civita covariant derivative
	defined by the invariant metric on $G / H$
	corresponding to $\langle \cdot, \cdot \rangle$.
	Let $\gamma \colon I \to G  / H$ be a curve in $G / H$ and 
	$g \colon I \to G$ be a horizontal lift of $\gamma$.
	Define $x \colon I \ni t  \mapsto x(t) = (T_e \ell_{g(t)})^{-1} \dot{g}(t) \in \liealg{m}$.
	Then $\gamma \colon I \to G  / H$ is a geodesic
	with respect to 
	$\nablaLC$ iff the ODE
	\begin{equation}
		\dot{x}(t) 
		=
		(\pr_{\liealg{m}} \circ \ad_{x(t)})^* (x(t))
	\end{equation}
	is satisfied for all $t \in I$.
	Here
	$(\pr_{\liealg{m}} \circ \ad_X)^* \colon \liealg{m} \to \liealg{m}$
	denotes the adjoint
	with respect to $\langle \cdot, \cdot \rangle$
	of the linear map
	defined for fixed $X \in \liealg{m}$ by
	$\pr_{\liealg{m}} \circ \ad_X \colon \liealg{m} \to \liealg{m}$.
	\begin{proof}
		We first recall Proposition~\ref{proposition:Levi_Civita_covariant_derivative}.
		The Levi-Civita covariant derivative on $G / H$ with respect to the invariant metric fulfills $\nablaLC = \nablaAlpha$,
		where $\alpha \colon \liealg{m} \times \liealg{m} \to \liealg{m}$ is given by $\alpha(X, Y) = - \tfrac{1}{2} [X, Y]_{\liealg{m}} + U(X, Y)$
		for all $X, Y \in \liealg{m}$ with
		\begin{equation}
			\label{corollary:geodesics_pseudo-riemannian_equation_on_m_definition_U}
			2 \langle U(X, Y), Z \rangle 
			=
			\langle [Z, X]_{\liealg{m}}, Y \rangle 
			+
			\langle X, [Z, Y]_{\liealg{m}} \rangle
			=
			- \big( 
			\langle [X, Z]_{\liealg{m}}, Y \rangle 
			+
			\langle X, [Y, Z]_{\liealg{m}} \rangle \big)
		\end{equation}
		Obviously,~\eqref{corollary:geodesics_pseudo-riemannian_equation_on_m_definition_U}
		is equivalently to
		\begin{equation}
			\label{corollary:geodesics_pseudo-riemannian_equation_on_m_definition_U_2}
			U(X, Y) 
			=
			- \tfrac{1}{2}
			\big( \langle (\pr_{\liealg{m}} \circ \ad_X)^*(Y) , Z \rangle
			+
			\langle (\pr_{\liealg{m}} \circ \ad_Y)^*(X) , Z \rangle \big)
		\end{equation}
		for all $X, Y, Z \in \liealg{m}$, where
		$(\pr_{\liealg{m}} \circ \ad_X)^*$ and
		$(\pr_{\liealg{m}} \circ \ad_Y)^*$
		denote the adjoints of
		the linear maps
		$(\pr_{\liealg{m}} \circ \ad_Y)$ and
		$(\pr_{\liealg{m}} \circ \ad_Y)$
		with respect to $\langle \cdot, \cdot \rangle$, respectively.
		Since $\langle \cdot, \cdot \rangle$ is non-degenerated,
		we can
		rewrite~\eqref{corollary:geodesics_pseudo-riemannian_equation_on_m_definition_U_2}
		equivalently as
		\begin{equation*}
			U(X, Y)
			=
			- \tfrac{1}{2}
			\big( (\pr_{\liealg{m}} \circ \ad_X)^*(Y) 
			+
			 (\pr_{\liealg{m}} \circ \ad_Y)^*(X) \big) .
		\end{equation*}
		Thus we obtain
		\begin{equation*}
			\alpha(X, X) 
			=
			- \tfrac{1}{2} [X, X]_{\liealg{m}} + U(X, X)
			=
			- (\pr_{\liealg{m}} \circ \ad_X)^* (X),
		\end{equation*}
		for all $X \in \liealg{m}$.
		Now Lemma~\ref{lemma:geodesics_on_G_H_ODE_on_m_characterization}
		yields the desired result.
	\end{proof}
\end{corollary}
As indicated above,
by applying
Corollary~\ref{corollary:geodesics_pseudo-riemannian_equation_on_m}
to a Lie group equipped with a left-invariant pseudo-Riemannian metric
considered as the reductive homogeneous space $G \cong G / \{ e \}$,
one obtains the following corollary
concerning geodesics on $G$.
Its statement is well-known and can be found
in~\cite[Ap. 2]{arnold:1978}.
We also refer to~\cite[Sec. 4]{Elshafei.etal:2023} for a
discussion of this characterization of geodesics
in the complex setting, where it is
named Euler-Arnold Formalism.

\begin{corollary}
	Let $G$ be a Lie group equipped with a left-invariant metric
	defined by the scalar product $\langle \cdot, \cdot \rangle \colon \liealg{g} \times \liealg{g} \to \field{R}$.
	Then $g \colon I \to G$ is a geodesic iff the curve
	$x \colon I \ni t \mapsto x(t) = (T_e \ell_{g(t)})^{-1} \dot{g}(t) \in \liealg{g}$ satisfies
	\begin{equation}
		\dot{x}(t) = (\ad_{x(t)})^* (x(t))
	\end{equation}
	for all $t \in I$.
	Here $(\ad_X)^* \colon \liealg{g} \to \liealg{g}$ denotes
	the adjoint of $\ad_X \colon  \liealg{g} \to \liealg{g}$ with
	respect to $\langle \cdot, \cdot \rangle$,
	where $X \in \liealg{g}$ is fixed.
	\begin{proof}
		Clearly, the Lie group $G$ equipped with the left-invariant metric
		defined by the scalar product $\langle \cdot, \cdot \rangle \colon \liealg{g} \times \liealg{g} \to \field{R}$
		can be viewed as the reductive homogeneous space
		$G / H$ for $H = \{e\}$ with reductive decomposition
		$\liealg{g} = \{0\} \oplus \liealg{g}$
		equipped with the pseudo-Riemannian metric defined by the $\Ad(\{e\})$-invariant scalar product $\langle \cdot , \cdot \rangle$ on $\liealg{g}$.
		Thus the assertion follows by Corollary~\ref{corollary:geodesics_pseudo-riemannian_equation_on_m}
		due to $\pr_{\liealg{m}} = \id_{\liealg{g}}$.
	\end{proof}
\end{corollary}

\subsection{Canonical Invariant Covariant Derivatives}

We now relate two particular
invariant covariant derivatives on $G / H$ to
the canonical affine connections of first and second kind
from~\cite[Sec. 10]{nomizu:1954}.
To this end, we list the two properties concerning
invariant covariant derivatives 
which correspond to the properties of invariant affine connections
from~\cite[Sec. 10, (A1) and (A2)]{nomizu:1954}.
This is the next definition.

\begin{definition}
	\label{definition:invariant_covariant_derivatives_properties}
	Let $\nablaAlpha$ be an invariant covariant derivative on $G / H$
	corresponding to the $\Ad(H)$-invariant bilinear map
	$\alpha \colon \liealg{m} \times \liealg{m} \to \liealg{m}$.
	The following properties of $\nablaAlpha$ are of particular interest:
	\begin{enumerate}
		\item
		\label{item:definition_invariant_covariant_derivatives_properties_1}
		The curves $\gamma_X \colon \field{R} \ni t \mapsto \pr(\exp(t X)) \in G / H$
		are geodesics with respect to $\nablaAlpha$
		for all $X \in \liealg{m}$.
		\item
		\label{item:definition_invariant_covariant_derivatives_properties_2}
		The curves $\gamma_X \colon \field{R} \ni t \mapsto \pr(\exp(t X)) \in G / H$
		are geodesics
		with respect to $\nablaAlpha$ for all $X \in \liealg{m}$
		and the parallel transport of $T_e \pr Z \in T_{\pr(e)} (G / H)$
		along $\gamma_X$ with respect to $\nablaAlpha$ is given by
		$\widehat{Z} \colon \field{R} \ni t \mapsto \big( T_{\exp(t X)} \pr \circ T_e \ell_{\exp(t X)} \big) Z \in T (G / H)$
		for all $Z \in \liealg{m}$.
	\end{enumerate}
\end{definition}
The next lemma is very similar some parts of~\cite[Sec. 10]{nomizu:1954}.

\begin{lemma}
	\label{lemma:invariant_covariant_derivatives_geodesics_parallel_one_parameter}
	Let $G / H$ be a reductive homogeneous space equipped with an invariant covariant derivative $\nablaAlpha$ corresponding to the
	$\Ad(H)$-invariant bilinear map
	$\alpha \colon  \liealg{m} \times \liealg{m} \to \liealg{m}$.
	Then the following assertions are fulfilled:
	\begin{enumerate}
		\item
		\label{item:lemma_invariant_covariant_derivatives_geodesics_parallel_one_parameter_geodesic}
		$\nablaAlpha$ fulfills the property from
		Definition~\ref{definition:invariant_covariant_derivatives_properties},
		Claim~\ref{item:definition_invariant_covariant_derivatives_properties_1}
		iff $\alpha(X, X) = 0$ holds for all $X \in \liealg{m}$.
		\item
		\label{item:lemma_invariant_covariant_derivatives_geodesics_parallel_one_parameter_parallel}
		$\nablaAlpha$ fulfills the property from
		Definition~\ref{definition:invariant_covariant_derivatives_properties},
		Claim~\ref{item:definition_invariant_covariant_derivatives_properties_2}
		iff $\alpha(X, Y) = 0$ 
		is fulfilled for all $X, Y \in \liealg{m}$.
	\end{enumerate}
	\begin{proof}
		Let $X \in \liealg{m}$ be arbitrary. We define the curve
		$\gamma_X \colon \field{R} \ni t \mapsto \pr(\exp(t X)) \in G / H$.
		Obviously, the curve $\field{R} \ni t \mapsto \exp(t X) \in G$
		is a horizontal lift of $\gamma$.
		Define $x \colon \field{R} \to \liealg{m}$ by
		$x(t) = (T_e \ell_{\exp(t X)})^{-1} \circ (T_e \pr\at{\liealg{m}})^{-1} \dot{\gamma}(t)$.
		Clearly, $x(t) = X$ holds for all $t \in \field{R}$.
		By Lemma~\ref{lemma:geodesics_on_G_H_ODE_on_m_characterization},
		the curve
		$\gamma \colon I \to G / H$ is a geodesic with
		respect to $\nablaAlpha$ iff $\alpha(X, X) = 0$ holds,
		i.e. Claim~\ref{item:lemma_invariant_covariant_derivatives_geodesics_parallel_one_parameter_geodesic} is shown.
		
		It remains to prove Claim~\ref{item:lemma_invariant_covariant_derivatives_geodesics_parallel_one_parameter_parallel}.
		To this end, let $Z \in \liealg{m}$ be arbitrary.
		We now define the vector field
		$\widehat{Z} \colon \field{R} 
		\ni t \mapsto \big(T_{\exp(t X)} \pr \circ T_e \ell_{\exp(t X)} \big) Z
		\in T(G / H)$
		along the curve $\gamma_X \colon \field{R} \ni t \mapsto \pr(\exp(t X)) \in G / H$.
		Next we consider the curve $z \colon \field{R} \to \liealg{m}$
		given by
		$z(t) 
		=
		\big(T_e \ell_{\exp(t X)} \big)^{-1} \circ \big(T_e \pr\at{\liealg{m}} \big)^{-1} \widehat{Z}(t) 
		=
		Z$.
		According to
		Corollary~\ref{corollary:parallel_vectorfields_ODE_on_m},
		the vector field $\widehat{Z} \colon I \to T(G / H)$
		is parallel along $\gamma$
		iff $\alpha(x(t), z(t)) = \alpha(X, Z) = 0$
		holds for
		all $t \in \field{R}$.
		This yields the desired result.
	\end{proof}
\end{lemma}
The next proposition can be viewed as a reformulation
of~\cite[Thm. 10.1]{nomizu:1954} and~\cite[Thm. 10.2]{nomizu:1954}.

\begin{proposition}
	\label{proposition:canonical_invariant_derivatives}
	Let $G / H$ be a reductive homogeneous space. 
	\begin{enumerate}
		\item
		\label{item:proposition_canonical_invariant_derivatives_first_kind}
		Define the $\Ad(H)$-invariant bilinear map 
		\begin{equation}
			\alpha
			\colon \liealg{m} \times \liealg{m}  \to \liealg{m} ,
			\quad
			(X, Y) \mapsto
			\alpha(X, Y) = \tfrac{1}{2}[X, Y]_{\liealg{m}} .
		\end{equation}
		The corresponding invariant covariant derivative
		$\nablaAlpha$ is the unique invariant covariant derivative 
		on $G / H$
		which is torsion free and satisfies
		Definition~\ref{definition:invariant_covariant_derivatives_properties},
		Claim~\ref{item:definition_invariant_covariant_derivatives_properties_1}.
		\item
		\label{item:proposition_canonical_invariant_derivatives_first_second_kind}
		Define the $\Ad(H)$-invariant bilinear map 
		\begin{equation}
			\alpha
			\colon \liealg{m} \times \liealg{m}  \to \liealg{m} ,
			\quad
			(X, Y) \mapsto
			\alpha(X, Y) = 0 .
		\end{equation}
		The corresponding invariant covariant derivative
		$\nablaAlpha$ is the unique invariant covariant derivative
		which satisfies 
		Definition~\ref{definition:invariant_covariant_derivatives_properties},
		Claim~\ref{item:definition_invariant_covariant_derivatives_properties_1}
		and
		Claim~\ref{item:definition_invariant_covariant_derivatives_properties_2}.
	\end{enumerate}
	\begin{proof}
		Claim~\ref{item:proposition_canonical_invariant_derivatives_first_second_kind}
		is an immediate consequence of
		Lemma~\ref{lemma:invariant_covariant_derivatives_geodesics_parallel_one_parameter},
		Claim~\ref{item:lemma_invariant_covariant_derivatives_geodesics_parallel_one_parameter_parallel}.
		
		It remains to proof Claim~\ref{item:proposition_canonical_invariant_derivatives_first_second_kind}.
		Obviously, $\nablaAlpha$ is torsion free for $\alpha(X, Y)= \tfrac{1}{2} [X, Y]_{\liealg{m}}$
		by Lemma~\ref{lemma:invariant_covariant_derivative_torsion}.
		Moreover, $\nablaAlpha$ fulfills
		Definition~\ref{definition:invariant_covariant_derivatives_properties},
		Claim~\ref{item:definition_invariant_covariant_derivatives_properties_1} 
		by
		Lemma~\ref{lemma:invariant_covariant_derivatives_geodesics_parallel_one_parameter},
		Claim~\ref{item:lemma_invariant_covariant_derivatives_geodesics_parallel_one_parameter_geodesic}
		because of $\alpha(X, X) = 0$ for all $X \in \liealg{m}$.
		It remains to prove the uniqueness of $\alpha$.
		To this end, let
		$\beta \colon \liealg{m} \times \liealg{m} \to \liealg{m}$
		be an $\Ad(H)$-invariant bilinear map
		and assume that
		the $\Ad(H)$-invariant bilinear map
		\begin{equation*}
			\gamma = \alpha + \beta \colon \liealg{m} \times \liealg{m} \to \liealg{m},
			\quad (X, Y) \mapsto \tfrac{1}{2} [X, Y]_{\liealg{m}} + \beta(X, Y)
		\end{equation*}
		fulfills $\gamma(X, X) = 0$ for all $X \in \liealg{m}$
		such that $\gamma$ defines the torsion-free invariant covariant
		derivative $\nabla^{\gamma}$ on $G / H$, i.e.
		$\gamma(X, Y) - \gamma(Y, X) = [X, Y]_{\liealg{m}}$
		holds
		for all $X, Y \in \liealg{m}$
		by
		Lemma~\ref{lemma:invariant_covariant_derivative_torsion}.
		This yields
		\begin{equation}
			\label{equation:proposition_canonical_invariant_derivatives_beta_symmetric}
			\begin{split}
				\gamma(X, Y) -  \gamma(Y, X)
				&=
				\tfrac{1}{2} [X, Y]_{\liealg{m}} + \beta(X, Y)
				-
				\big( \tfrac{1}{2} [Y, X]_{\liealg{m}} + \beta(Y, X)\big) \\
				&=
				[X, Y]_{\liealg{m}} + \beta(X, Y) - \beta(Y, X) \\
				&= 
				[X, Y]_{\liealg{m}} .
			\end{split}
		\end{equation}
		By~\eqref{equation:proposition_canonical_invariant_derivatives_beta_symmetric},
		one obtains
		$\beta(X, Y) - \beta(Y, X) = 0$
		for all $X, Y \in \liealg{m}$, i.e.
		$\beta(X, Y) = \beta(Y, X)$ is symmetric.
		Moreover, we have $\beta(X, X)  = 0$
		for all $X \in \liealg{m}$ due to
		\begin{equation*}
			0 = \gamma(X, X) = \tfrac{1}{2}[X, X]_{\liealg{m}}
			+ \beta(X,X) = \beta(X, X) .
		\end{equation*}
		Thus $\beta \colon \liealg{m} \times \liealg{m} \to \liealg{m}$
		is a symmetric bilinear map that fulfills $\beta(X, X) = 0$
		for all $X \in \liealg{m}$.
		By polarization,
		we obtain $\beta(X, Y) = 0$ for all $X, Y \in \liealg{m}$.
		Hence $\gamma = \alpha + \beta = \alpha$ holds, i.e.
		$\alpha \colon \liealg{m} \times \liealg{m} \ni (X, Y) \mapsto \tfrac{1}{2} [X, Y]_{\liealg{m}} \in \liealg{m}$
		is the unique $\Ad(H)$-invariant bilinear map
		that satisfies
		$\alpha(X, Y) - \alpha(Y, X) = [X, Y]_{\liealg{m}}$
		and $\alpha(X, X) =0$ for all $X, Y \in \liealg{m}$.
		This yields the desired result. 
	\end{proof}
\end{proposition}

\begin{definition}
	\label{definition:canonical_invariant_derivatives}
	Let $G / H$ be a reductive homogeneous space.
	\begin{enumerate}
		\item
		The invariant covariant derivative defined by
		$\alpha(X, Y) =  \tfrac{1}{2}[X, Y]_{\liealg{m}}$
		for all $X, Y \in \liealg{m}$ is called
		the
		canonical invariant covariant derivative of first kind.
		It is denoted by $\nablaCan$.
		\item
		The invariant covariant derivative defined by
		$\alpha(X, Y) =  0$
		for all $X, Y \in \liealg{m}$ is called the
		canonical invariant covariant derivative of second kind.
		It is denoted by $\nablaCanSecond$.
	\end{enumerate}
\end{definition}

\begin{remark}
	By
	Proposition~\ref{proposition:one-to-one-correspondence-invariant-affine-connection-covariant-derivatives},
	the canonical covariant derivatives
	of first kind $\nablaCan$ and of second kind $\nablaCanSecond$ from
	Definition~\ref{definition:canonical_invariant_derivatives}
	correspond to the canonical affine
	connections of first and second kind
	form~\cite[Sec. 10]{nomizu:1954},
	respectively.
\end{remark}

\begin{remark}
	\label{remark:horizontal_lift_covariant_derivative_naturally_redutive_and_symmetric}
	Assume that $G / H$ is a naturally reductive homogeneous space.
	Then the Levi-Civita covariant derivative coincides with the canonical
	covariant derivative of first kind
	by
	Remark~\ref{remark:invariant_covariant_derivatives_levi_civita}, i.e $\nablaLC = \nablaCan$ holds.
	This is has already been proven
	in~\cite[Thm. 13.1 and Eq. (13.2)]{nomizu:1954}.
\end{remark}

\begin{remark}
	Let $G / H$ be equipped with an invariant pseudo-Riemannian metric.
	Then $\nablaCanSecond$ is an invariant metric covariant derivative
	on $G / H$ by
	Proposition~\ref{proposition:invariant_covariant_derivative_compatible_with_structure_skew_adjoint}.
\end{remark}

We briefly comment on the canonical covariant derivatives on
symmetric homogeneous spaces in the next remark
following~\cite[Thm. 15.1]{nomizu:1954}.

\begin{remark}
	\label{remark:symmetric_spaces_invariant_covariant_derivatives}
	Let $(G, H, \sigma)$ be a symmetric pair and let
	$G / H$ be the corresponding
	symmetric homogeneous space.
	Let
	$\liealg{g} = \liealg{h} \oplus \liealg{m}$
	denote the canonical reductive decomposition.
	Then $[X, Y] \in \liealg{h}$ holds for all
	$X, Y \in \liealg{m}$
	by Lemma~\ref{lemma:symmetric_pair_canonical_reductive_decomposition}.
	Therefore $\tfrac{1}{2}[X, Y]_{\liealg{m}} = 0$ 
	is fulfilled for all $X, Y \in \liealg{m}$.
	Hence $\nablaCan = \nablaCanSecond$ holds
	by Proposition~\ref{proposition:canonical_invariant_derivatives}.
\end{remark}
Moreover, for pseudo-Riemannian symmetric spaces, we obtain the 
following remark whose
statement can be found in~\cite[Thm. 15.6]{nomizu:1954}.

\begin{remark}
	Let $G / H$ be a pseudo-Riemannian symmetric homogeneous
	space.
	Then one has $\nablaLC = \nablaCan = \nablaCanSecond$
	by Remark~\ref{remark:horizontal_lift_covariant_derivative_naturally_redutive_and_symmetric}
	combined with
	Remark~\ref{remark:pseudo-Riemannian-symmetric-naturally-reductive}.
\end{remark}
We end this section by specializing Corollary~\ref{corollary:parallel_vectorfields_ODE_on_m}
on parallel vector fields along curves
to the
canonical covariant derivatives $\nablaCan$ and
$\nablaCanSecond$.

\begin{corollary}
	\label{corollary:parallel_vectorfields_ODE_on_m_canonical_covariant_derivatives}
	Let $\widehat{Z} \colon I \to T (G / H)$ be a vector field along
	the curve $\gamma \colon I \to G / H$.
	Using the notation of
	Corollary~\ref{corollary:parallel_vectorfields_ODE_on_m},
	the following assertions are fulfilled:
	\begin{enumerate}
		\item
		\label{corollary:parallel_vectorfields_ODE_on_m_canonical_covariant_derivatives_first_kind}
		$\widehat{Z}$ is parallel along $\gamma$ with respect to $\nablaCan$
		iff
		\begin{equation}
			\dot{z}(t) = - \tfrac{1}{2} [x(t), z(t)]_{\liealg{m}}
		\end{equation}
		holds for all $t \in I$.
		\item
		\label{corollary:parallel_vectorfields_ODE_on_m_canonical_covariant_derivatives_second_kind}
		$\widehat{Z}$ is parallel along $\gamma$
		with respect to $\nablaCanSecond$ iff
		\begin{equation}
			\dot{z}(t) = 0
		\end{equation}
		is fulfilled for all $t \in I$.
	\end{enumerate}
\end{corollary}

\begin{remark}
	A similar description of parallel vector fields
	as in 
	Corollary~\ref{corollary:parallel_vectorfields_ODE_on_m_canonical_covariant_derivatives},
	Claim~\ref{corollary:parallel_vectorfields_ODE_on_m_canonical_covariant_derivatives_first_kind}
	has
	already appeared in \cite[Prop. 2.12]{smith:1993} for the
	special case, where $G / H$ is a normal naturally reductive space,
	see Remark~\ref{remark:normal_naturally_recductive}
	for this notion,
	and $\gamma \colon I \to G / H$ is a geodesic,
	i.e. for $x \colon I \to \liealg{m}$ being constant.
\end{remark}

\section{Conclusion}

We considered invariant covariant derivatives on
a reductive homogeneous space in detail.
We proved that they are uniquely 
determined by evaluating them on fundamental vector fields.
Moreover, we provided a new proof for their existence by expressing
them in terms of horizontally lifted vector fields.
By this result, a characterization of parallel vector fields
along curves
in a reductive homogeneous space equipped
with an invariant covariant derivative is obtained.
In addition, the so-called canonical covariant derivatives of first and
second kind corresponding the canonical affine connections of first and
second kind from~\cite{nomizu:1954} are considered.

\section*{Acknowledgments}
This work has been supported by the 
German Federal Ministry of Education and Research
(BMBF-Projekt 05M20WWA: Verbundprojekt 05M2020 - DyCA).


\begin{thebibliography}{10}
	
	\bibitem {arnold:1978}
	\chairxauthorbibfont{Arnold, V.~I.:}\newblock \chairxtitlebibfont{Mathematical
		methods of classical mechanics}, vol.~60 in \chairxseriesbibfont{Graduate
		Texts in Mathematics}.
	\newblock Springer-Verlag, New York-Heidelberg, 1978.
	\newblock Translated from the Russian by K. Vogtmann and A. Weinstein.
	\par\csname arnold:1978chairxnote\endcsname
	
	\bibitem {berestovskii.nikonorov:2020}
	\chairxauthorbibfont{Berestovskii, V., Nikonorov, Y.:}\newblock
	\chairxtitlebibfont{Riemannian manifolds and homogeneous geodesics}.
	\newblock \chairxseriesbibfont{Springer Monographs in Mathematics}.
	\newblock Springer, Cham, Cham, 2020.
	\par\csname berestovskii.nikonorov:2020chairxnote\endcsname
	
	\bibitem {chevalley:1946}
	\chairxauthorbibfont{Chevalley, C.:}\newblock \chairxtitlebibfont{Theory of Lie
		groups I}.
	\newblock Princeton University Press, Princeton, NJ, 1946.
	\par\csname chevalley:1946chairxnote\endcsname
	
	\bibitem {Elshafei.etal:2023}
	\chairxauthorbibfont{Elshafei, A., Ferreira, A.~C., Reis, H.:}\newblock
	\chairxtitlebibfont{Geodesic completeness of pseudo and
		holomorphic-{R}iemannian metrics on {L}ie groups}.
	\newblock Nonlinear Anal.  \textbf{232} (2023), Paper No. 113252, 37.
	\par\csname Elshafei.etal:2023chairxnote\endcsname
	
	\bibitem {gallier.quaintance:2020}
	\chairxauthorbibfont{Gallier, J., Quaintance, J.:}\newblock
	\chairxtitlebibfont{Differential geometry and {L}ie groups---a computational
		perspective}, vol.~12 in \chairxseriesbibfont{Geometry and Computing}.
	\newblock Springer, Cham, Cham, 2020.
	\par\csname gallier.quaintance:2020chairxnote\endcsname
	
	\bibitem {helgason:1978}
	\chairxauthorbibfont{Helgason, S.:}\newblock \chairxtitlebibfont{Differential
		geometry, {L}ie groups, and symmetric spaces}, vol.~80 in
	\chairxseriesbibfont{Pure and Applied Mathematics}.
	\newblock Academic Press, Inc., New York-London, 1978.
	\par\csname helgason:1978chairxnote\endcsname
	
	\bibitem {jurdjevic.markina.leite:2023}
	\chairxauthorbibfont{Jurdjevic, V., Markina, I., Silva~Leite, F.:}\newblock
	\chairxtitlebibfont{Symmetric spaces rolling on flat spaces}.
	\newblock J. Geom. Anal.  \textbf{33}.3 (2023), Paper No. 94, 33.
	\par\csname jurdjevic.markina.leite:2023chairxnote\endcsname
	
	\bibitem {knapp:2002}
	\chairxauthorbibfont{Knapp, A.~W.:}\newblock \chairxtitlebibfont{Lie groups
		beyond an introduction}, vol. 140 in \chairxseriesbibfont{Progress in
		Mathematics}.
	\newblock Birkh\"{a}user Boston, Inc., Boston, MA, second. edition, 2002.
	\par\csname knapp:2002chairxnote\endcsname
	
	\bibitem {kobayashi.nomizu:1963a}
	\chairxauthorbibfont{Kobayashi, S., Nomizu, K.:}\newblock
	\chairxtitlebibfont{Foundations of Differential Geometry {I}}.
	\newblock \chairxseriesbibfont{Interscience Tracts in Pure and Applied
		Mathematics} no. \textbf{15}.
	\newblock John Wiley {\&} Sons, New York, London, Sydney, 1963.
	\par\csname kobayashi.nomizu:1963achairxnote\endcsname
	
	\bibitem {lee:2013}
	\chairxauthorbibfont{Lee, J.~M.:}\newblock \chairxtitlebibfont{Introduction to
		smooth manifolds}, vol. 218 in \chairxseriesbibfont{Graduate Texts in
		Mathematics}.
	\newblock Springer, New York, New York, second. edition, 2013.
	\par\csname lee:2013chairxnote\endcsname
	
	\bibitem {lee:2018}
	\chairxauthorbibfont{Lee, J.~M.:}\newblock \chairxtitlebibfont{Introduction to
		{R}iemannian manifolds}, vol. 176 in \chairxseriesbibfont{Graduate Texts in
		Mathematics}.
	\newblock Springer, Cham, Cham, 2018.
	\par\csname lee:2018chairxnote\endcsname
	
	\bibitem {michor:2008}
	\chairxauthorbibfont{Michor, P.~W.:}\newblock \chairxtitlebibfont{Topics in
		differential geometry}, vol.~93 in \chairxseriesbibfont{Graduate Studies in
		Mathematics}.
	\newblock American Mathematical Society, Providence, RI, 2008.
	\par\csname michor:2008chairxnote\endcsname
	
	\bibitem {nomizu:1954}
	\chairxauthorbibfont{Nomizu, K.:}\newblock \chairxtitlebibfont{Invariant affine
		connections on homogeneous spaces}.
	\newblock Amer. J. Math.  \textbf{76} (1954), 33--65.
	\par\csname nomizu:1954chairxnote\endcsname
	
	\bibitem {oneill:1983}
	\chairxauthorbibfont{O'Neill, B.:}\newblock
	\chairxtitlebibfont{Semi-{R}iemannian geometry}, vol. 103 in
	\chairxseriesbibfont{Pure and Applied Mathematics}.
	\newblock Academic Press, Inc., New York, 1983.
	\newblock With applications to relativity.
	\par\csname oneill:1983chairxnote\endcsname
	
	\bibitem {rabenoro.pennec:2023}
	\chairxauthorbibfont{Rabenoro, D., Pennec, X.:}\newblock
	\chairxtitlebibfont{The geometry of Riemannian submersions from compact Lie
		groups. Application to flag manifolds}, 2023.
	\par\csname rabenoro.pennec:2023chairxnote\endcsname
	
	\bibitem {rudolph.schmidt:2017}
	\chairxauthorbibfont{Rudolph, G., Schmidt, M.:}\newblock
	\chairxtitlebibfont{Differential geometry and mathematical physics. {P}art
		{II}}.
	\newblock \chairxseriesbibfont{Theoretical and Mathematical Physics}.
	\newblock Springer, Dordrecht, Dordrecht, 2017.
	\newblock Fibre bundles, topology and gauge fields.
	\par\csname rudolph.schmidt:2017chairxnote\endcsname
	
	\bibitem {smith:1993}
	\chairxauthorbibfont{Smith, S.~T.:}\newblock \chairxtitlebibfont{Geometric
		Optimization Methods for Adaptive Filtering}.
	\newblock PhD thesis, Harvard University, Cambridge, 1993.
	\par\csname smith:1993chairxnote\endcsname
	
	\bibitem {xu:2022}
	\chairxauthorbibfont{Xu, M.:}\newblock \chairxtitlebibfont{Submersion and
		homogeneous spray geometry}.
	\newblock J. Geom. Anal.  \textbf{32}.6 (2022), Paper No. 172, 43.
	\par\csname xu:2022chairxnote\endcsname
	
\end{thebibliography}
\end{document}